\documentclass[final,3p]{elsarticle}

\usepackage{lineno}
\usepackage[colorlinks=true,breaklinks=true,pdftex]{hyperref}

\modulolinenumbers[1]

\biboptions{authoryear}

\usepackage{dblfloatfix}
\usepackage{float}
\usepackage{graphicx}
\usepackage{lineno,hyperref}
\usepackage{times}
\usepackage{epsfig}
\usepackage{bm}
\usepackage{breakurl}
\usepackage{graphicx}
\usepackage{amsmath,amssymb, amsfonts, euscript, mathrsfs}
\usepackage{subeqnarray}
\usepackage{cases}
\usepackage{caption}
\usepackage{fancybox}
\usepackage{picins}
\usepackage{picinpar} 
\usepackage[table]{xcolor}
\usepackage{wrapfig}
\definecolor{lightgray}{gray}{0.9}
\definecolor{lightblue}{rgb}{0.98,0.98,1.0}
\usepackage{subfigure}
\usepackage{multirow}
\usepackage{booktabs}
\usepackage{ulem}
\usepackage{comment}
\usepackage{lipsum}
\usepackage{amsmath}
\usepackage{tabularx}
\usepackage{overpic}
\usepackage{adjustbox}
\usepackage{natbib}
\usepackage{doi}

\usepackage{algorithm}
\usepackage{algorithmicx}
\usepackage{algpseudocode}

\usepackage{multirow} 
\usepackage{xcolor}

\newcommand{\T}{\mathsf{T}}

\makeatletter
\DeclareRobustCommand\bigop[1]{%
	\mathop{\vphantom{\sum}\mathpalette\bigop@{#1}}\slimits@
}
\newcommand{\bigop@}[2]{%
	\vcenter{%
		\sbox\z@{$#1\sum$}%
		\hbox{\resizebox{\ifx#1\displaystyle.9\fi\dimexpr\ht\z@+\dp\z@}{!}{$\m@th#2$}}%
	}%
}
\makeatother

\newcommand{\stkout}[1]{\ifmmode\text{\sout{\ensuremath{#1}}}\else\sout{#1}\fi}

\begin{document}

\begin{frontmatter}

\title{Density-based isogeometric topology optimization of shell structures}



\author[first]{Qiong Pan}
\ead{PQ2019@mail,ustc.edu.cn}
\author[first]{Xiaoya Zhai\corref{cor}}
\ead{xiaoya93@ustc.edu.cn}
\author[first]{Falai Chen}
\ead{chenfl@ustc.edu.cn}
\cortext[cor]{Corresponding author}
\address[first]{School of Mathematical Sciences, University of Science and Technology of China, Hefei, China}

\begin{abstract}
Shell structures with a high stiffness-to-weight ratio are desirable in various engineering applications. In such scenarios, topology optimization serves as a popular and effective tool for shell structures design. Among the topology optimization methods, solid isotropic material with penalization method(SIMP) is often chosen due to its simplicity and convenience. However, SIMP method is typically integrated with conventional finite element analysis(FEA) which has limitations in computational accuracy. Achieving high accuracy with FEA needs a substantial number of elements, leading to computational burdens. In addition, the discrete representation of the material distribution may result in rough boundaries and checkerboard structures. To overcome these challenges, this paper proposes an isogeometric analysis(IGA) based SIMP method for optimizing the topology of shell structures based on Reissner-Mindlin theory. We use NURBS to represent both the shell structure and the material distribution function with the same basis functions, allowing for higher accuracy and smoother boundaries. The optimization model takes compliance as the objective function with a volume fraction constraint and the coefficients of the density function as design variables. The Method of Moving Asymptotes is employed to solve the optimization problem, resulting in an optimized shell structure defined by the material distribution function. To obtain fairing boundaries in the optimized shell structure, further process is conducted by fitting the boundaries with fair B-spline curves automatically. Furthermore, the IGA-SIMP framework is applied to generate porous shell structures by imposing different local volume fraction constraints. Numerical examples are provided to demonstrate the feasibility and efficiency of the IGA-SIMP method, showing that it outperforms the FEA-SIMP method and produces smoother boundaries.
\end{abstract}


\end{frontmatter}


\section{Introduction}

\label{sec:Intro} 
	Shell structures occupy a crucial position in various engineering fields, including architectural, mechanical, aeronautical, and biomechanical engineering \citep{rotter1998shell, ventsel2002thin}. 
	These structures are favored for their exceptionally high stiffness-to-weight ratios and well-performing force transfer properties. In particular, porous shell structures possess immense aesthetic value in architectural and furniture designs and play a crucial role in the development of medical implants and orthopedic devices \citep{berlinberg2022minimum, chen2022porous, hu2023parametric,shichman2022outcomes}. And for the design of shell structures, topology optimization is considered an effective tool \citep{maute1997adaptive}.

	Curved elements in shell structure analysis were first introduced by ~\cite{ahmad1970analysis}, allowing for the conversion of topology optimization of shell structures into a degenerate surface topology optimization problem on the mid-surfaces of the shells. Since then, various topology optimization methods have been applied to optimize shell structures, such as the homogenization method
	\citep{bendsoe1988generating},  
	Solid Isotropic Material with Penalization (SIMP) method
	\citep{bendsoe1989optimal},  
	level set method (LSM)
	\citep{osher1988fronts},  
	and moving morphable components (MMC) method
	\citep{guo2014doing}.  
	\cite{ansola2002integrated} firstly applied the homogenization method to perform shape and topology optimization alternatively for shell structures using the second rank microstructures proposed 
	by ~\cite{soto1993modelling}. 
	Later on, \cite{hassani2013simultaneous} improved the results by optimizing shape and topology simultaneously under the SIMP framework.
	The level set method
	\citep{townsend2019level,ye2019topology}  
	has also been employed to study the topology optimization of shell structures, whose boundaries can be tracked through the update of the level set function. 
	Additionally, the MMC method has recently been developed for explicit shell structure optimization.
	~\cite{huo2022topology} utilized conformal mapping theory to transform the optimization of shell structures into a distribution problem of basic components controlled by a few variables in the parameter domain, thereby reducing the degrees of freedom in the optimization. 
	Among all these methods for shell structure optimization, SIMP is the most commonly used one because it is not only straightforward and intuitive but also easy to implement.

	The aforementioned methods for shell structure optimization are based on Finite Element Analysis (FEA) 
	which allows for  
	applying forces and boundary conditions directly.
	However, the accuracy of analysis results heavily depends on the resolution of finite element meshes. To obtain high-quality results, refinement of meshes  is inevitable, which leads to high computation costs and low efficiency. To address these issues, isogeometric analysis (IGA) provides a promising alternative due to its efficiency, accuracy and representation consistency between geometry and analysis models.

	IGA was introduced by~\cite{hughes2005isogeometric} to integrate two related disciplines--computer-aided design and computer-aided engineering. In recent years, it has become an effective alternative tool to FEA in a variety of engineering fields 
	\citep{cottrell2009isogeometric},  
	and has been applied in topology optimization of shell structures. The first attempt was proposed by \cite{kang2016isogeometric} based on trimmed surface analysis (TSA).
	However, this method may create irregular trimming curves with self-intersection. To address this issue, 
	\cite{zhang2020explicit} and \cite{jiang2023explicit} 
	combined trimming-based IGA with the moving morphable void (MMV) method 
	\citep{Zhang2018CMS}  
	to design shell structures. They introduce a group of moving morphable voids with convex control polygon as the building blocks and the optimal structure is achieved by optimizing the parameters describing the locations, morphologies and layout of the voids in the design domain. 
	Subsequently, \cite{zhang2020stress} further developed an IGA-MMV approach to address stress-related optimization problems which results in a design with explicit and clear boundaries.
	IGA-based approaches have several advantages, including exact geometry representation, mature refinement schemes, and higher order of continuity, which overcome the limitations of FEA-based approaches. 
	However, the numerical implementation of the above TSA-related methods is complicated, and computational costs may increase considerably when a large number of trimming curves or voids are used. 
	To address these drawbacks, we propose a straightforward approach to the topology optimization problem of shell structures through IGA-SIMP based on Reissner-Mindlin theory
	\citep{bischoff2004models}.  

	The basic idea of our approach is as follows. First, we represent the mid-surface of a shell structure using NURBS basis functions and construct a density function that is represented by the same NURBS basis functions as those of the shell structure. 
	Then we take the control coefficients of the density function as design variables and minimize the compliance of the shell structure under a volume constraint. 
	A multi-level system
	\citep{nagy2013isogeometric}  
	is utilized for the 
	design and the optimization problem is solved by the Method of Moving Asymptotes(MMA)  
	\citep{svanberg1987method}.  
	Finally, a post-processing step is performed to obtain fair boundaries by fitting the boundaries with fair B-spline curves. 
	In addition, by imposing the local volume constraint~\citep{wu2017infill}, our method can also be applied to design porous shell structures conveniently.  
	Several representative shell structures are provided with different load and boundary conditions to verify the effectiveness and robustness of the proposed method.

	In summary, our contributions are listed as follows.
	\begin{itemize}
		\item [(1)]
		We propose an IGA-based SIMP method for optimizing shell structures. In this method, the density function and shell model are represented using NURBS with the same B-spline basis functions.
		\item [(2)]
		We perform a post-processing step to obtain explicit and fair boundaries automatically.
	\end{itemize}
	
	The paper is organized as follows. 
	Section ~\ref{sec:pre} presents the theoretical foundations of NURBS functions and Reissner-Mindlin shell theory. 
	Section ~\ref{sec:method} proposes formulations for topology optimization of (porous) shell structures and provides some details of numerical implementation. 
	In Section~\ref{sec:post-processing}, we put forward a post-processing algorithm to obtain fair boundaries automatically . 
	Section~\ref{sec:results} illustrates numerical examples of typical shell structures using the proposed IGA-SIMP approach and compares the results with those obtained using the FEA-SIMP method. 
	Additionally, experiments on generating porous shell structures and results of boundary fairing 
	are presented, followed by a discussion on parameter selection and a comparison between the IGA-SIMP method and other IGA-based methods.
	Finally, section ~\ref{sec:conclu} closes the paper with some concluding remarks.



\section{Preliminaries}
	\label{sec:pre}
	In this section, we present some basic knowledge about NURBS(Non-uniform rational B-splines) and Reissner-Mindlin theory, 
	which lay the groundwork for our IGA-SIMP method.
	
	\subsection{NURBS representations}
	\label{sec:NURBS}
	In computer aided design and geometric modeling, NURBS is the standard representation of geometric objects.
	Given a degree $p$ and a knot vector
	$\mathcal{S} = \left\{s_0,s_1,\cdots,s_{m+p+1}\right\}$ $(s_0 \leq s_1 \leq \cdots \leq s_{m+p+1})$, 
	a B-spline basis function is defined recursively as
	\begin{equation} 
		\begin{aligned}
			&N_{i}^0 (s) = 
			\begin{cases}
				1, \quad &\text{if} \ s_i \leq s \leq s_{i+1}; \\
				0, \quad & \text{otherwise};
			\end{cases} \qquad p = 0 \\
			&N_{i}^p (s) = \frac{s-s_i}{s_{i+p}-s_i} N_{i}^{p-1} (s) + \frac{s_{i+p+1}-s}{s_{i+p+1}-s_{i+1}} N_{i+1}^{p-1} (s), \, p \geq 1
		\end{aligned}
	\end{equation}
	where $0/0$ is defined as zero.

	Given a set of control points ${\mathbf P}_i\in {\mathbb R}^2$ and corresponding weights $w_i>0$, $i=0,1,\ldots, m$, a NURBS curve is defined by
	\begin{equation}
		\label{bcurve}
		{\mathbf P}(s)=\sum_{i=0}^m R_i^p(s){\mathbf P}_i, \qquad s\in [s_p,s_{m+1}],   
	\end{equation}
	where
	\begin{equation}
		R_{i}^p(s) = \frac{w_i \cdot N_{i}^p(s)}{\sum\limits_{j=0}^{m} w_j  N_{j}^p(s)}, \quad i=0,1,\ldots,m
	\end{equation}
	are the univariate NURBS basis functions. 
	
\begin{figure*}[t]
		\begin{center}
			\begin{overpic}[width=0.8\textwidth]{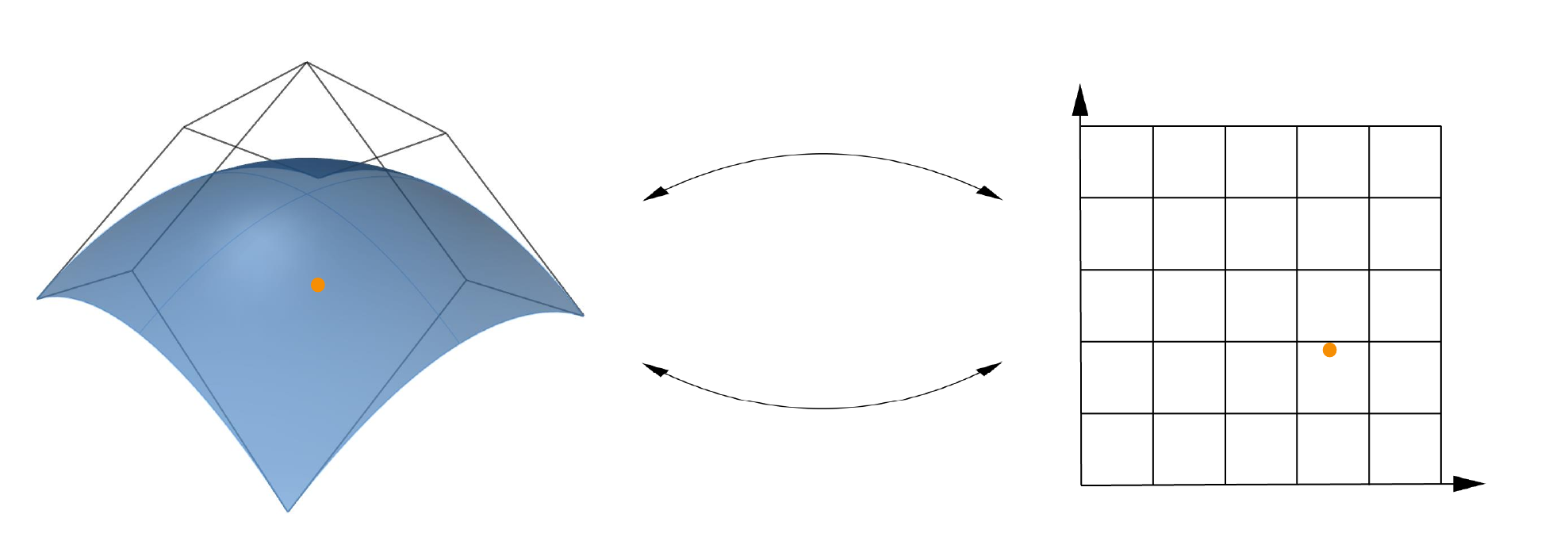}
				\put(0,16){$\boldsymbol{P}_{00}$}
				\put(7,19){$\boldsymbol{P}_{10}$}
				\put(8,28){$\boldsymbol{P}_{01}$}
				\put(17.5,3.5){$\boldsymbol{P}_{20}$}
				\put(18,26){$\boldsymbol{P}_{02}$}
				\put(18,32.5){$\boldsymbol{P}_{11}$}
				\put(28,27.5){$\boldsymbol{P}_{01}$}
				\put(29,18.5){$\boldsymbol{P}_{01}$}
				\put(37,14){$\boldsymbol{P}_{22}$}
				\put(40,28){$\boldsymbol{S}(s,t) = \sum\limits_{ij}R_{ij}^{pq}(s,t) \boldsymbol{P}_{ij}$}
				\put(49,17){\color{orange}{bijection}}
				\put(40,7){$\tilde{\rho}(s,t) = \sum\limits_{ij}R_{ij}^{pq}(s,t) \rho_{ij}$}
				\put(68.5,31){$t$}
				\put(95,4){$s$}
				\put(12,0){NURBS surface}
				\put(74,0){Parametric space}
			\end{overpic}
			\caption{A bi-quadratic NURBS surface and its parametric space.}
			\label{fig:NURBSSurface}
		\end{center}
	\end{figure*}
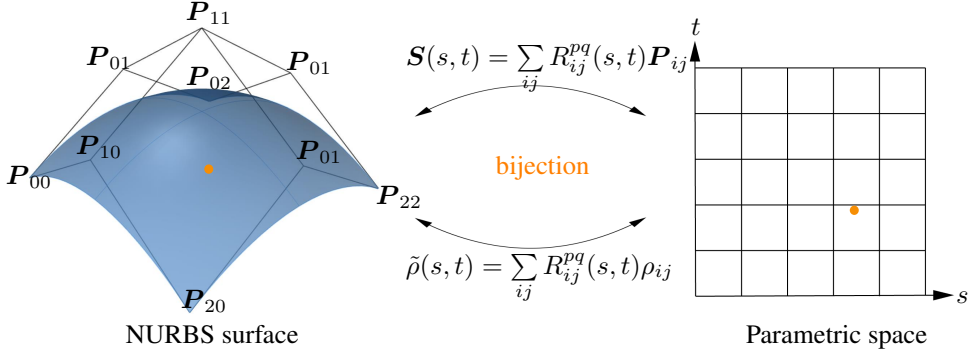

	Similarly, given a set of points ${\mathbf P}_{ij}\in {\mathbb R}^3$ and the corresponding weights, $i=0,1,\ldots,m$, $j=0,1,\ldots, n$,
	a tensor product NURBS surface of bi-degree $(p,q)$ is defined by
	\begin{equation}
		\label{eq:explicit}
		\boldsymbol{S}(s,t) = \sum\limits_{i=0}^{m}\sum\limits_{j=0}^{n} R_{ij}^{pq}(s,t) \boldsymbol{P}_{ij},
	\end{equation}
	where $s\in [s_p, s_{m+1}], t\in [t_q,t_{n+1}]$, and
	\begin{equation}
		R_{ij}^{pq}(s,t) = \frac{ w_{ij}  N_{i}^{p}(s) N_{j}^{q}(t)  }{\sum\limits_{k=0}^{m}\sum\limits_{l=0}^{n}  w_{kl}  N_{k}^{p}(s) N_{l}^{q}(t)}, \end{equation}
	$ i=0,1,\ldots,m, \, j=0,1,\ldots, n
	$ are bivariate NURBS basis functions defined over the knot vectors
	$\mathcal{S} = \left\{s_0,s_1,\cdots,s_{m+p+1}\right\}$ and $\mathcal{T} = \left\{t_0,t_1,\cdots,t_{n+q+1}\right\}$ with weights $\left\{w_{ij}\right\}_{i=0,j=0}^{m,n}$. 
	Fig.~\ref{fig:NURBSSurface} illustrates a NURBS surface.

	\medskip
	
	In this paper, we will represent the mid-surface, the material distribution function and the displacement vector of a shell structure in the NURBS forms.

	\subsection{Reissner–Mindlin theory}
\begin{figure*}[ht]
		\begin{center}
			\begin{overpic}[width=0.4\textwidth]{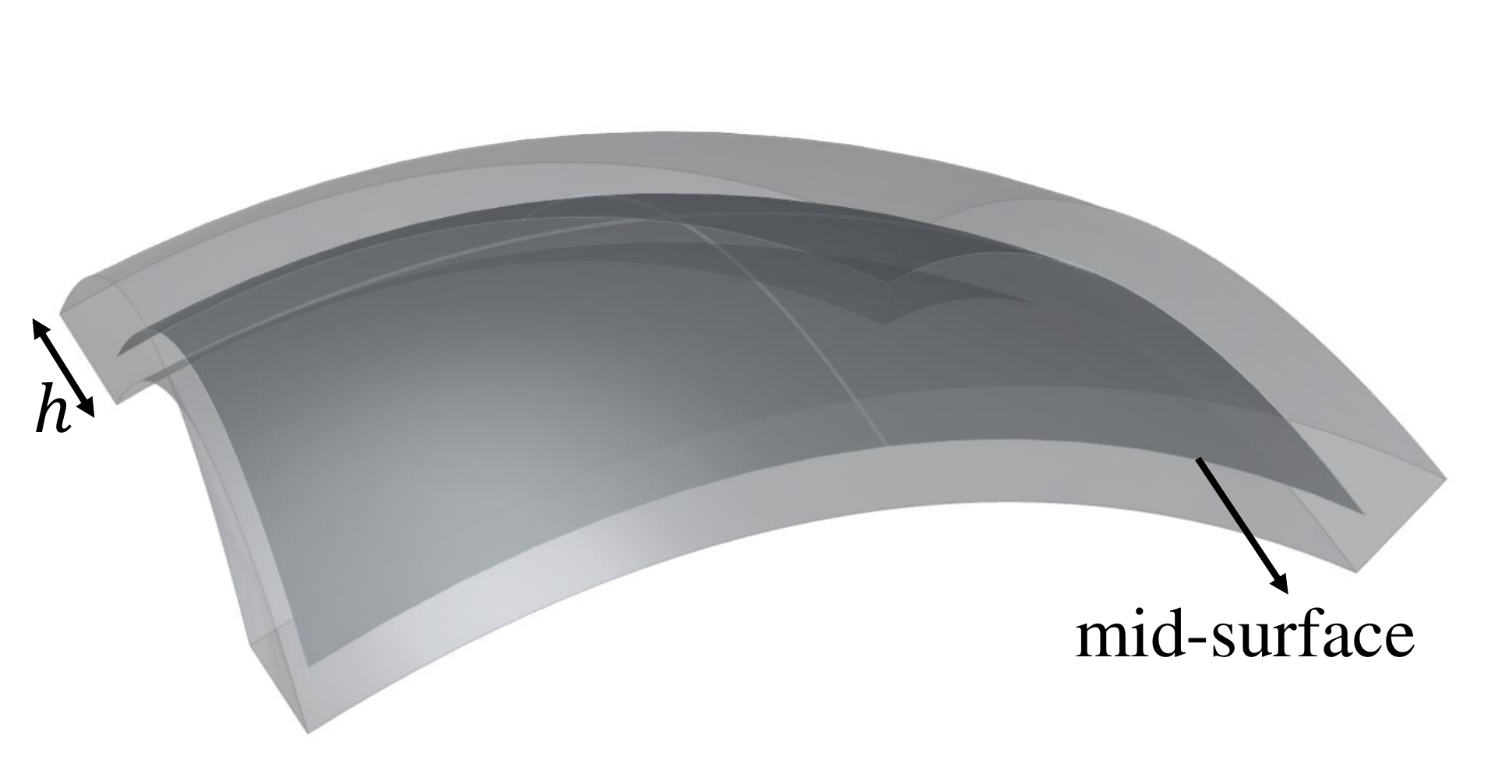}
				\put(20,19){$\boldsymbol{X}(s,t)$}
				\put(35,0){(a) Shell structure}
			\end{overpic}
			\qquad
			\begin{overpic}[width=0.4\textwidth]{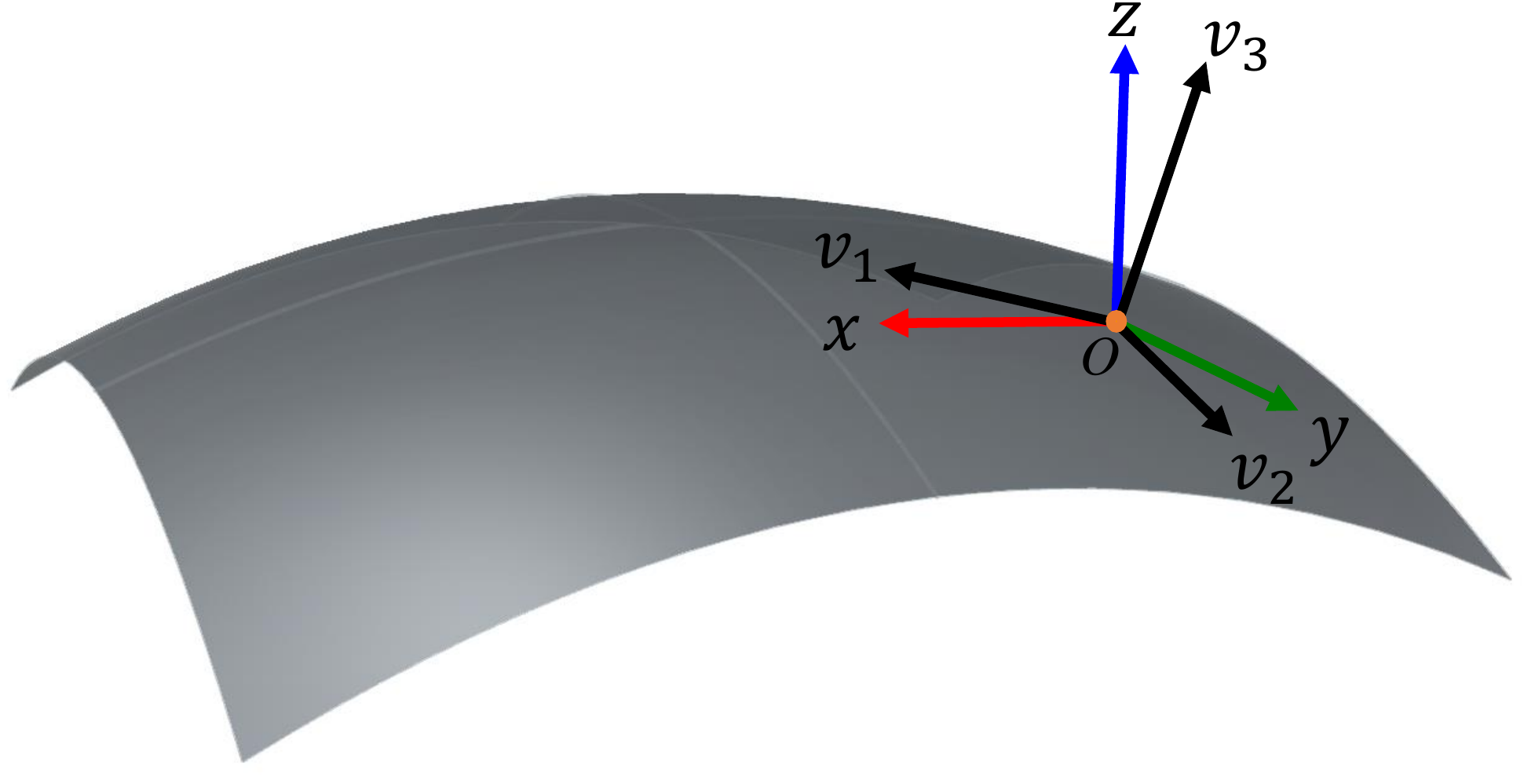}
				\put(25,19){$\boldsymbol{S}(s,t)$}
				\put(35,0){(b) Mid-surface}
			\end{overpic}
			\caption{A shell structure (a) and its mid-surface (b). There is a local orthonormal coordinate system $\{\boldsymbol{v}_1,\boldsymbol{v}_2,\boldsymbol{v}_3\}$ at the point $\boldsymbol{O}$. Here $\boldsymbol{v}_3$ is the normal of the point $\boldsymbol{O}$ on the mid-surface and $\boldsymbol{v}_1, \boldsymbol{v}_2$ are two tangent vectors of the mid-surface.} 
			\label{fig: Shell structure and its mid-surface}
		\end{center}
	\end{figure*}

	The Reissner–Mindlin theory, based on zero normal stress assumption, is one of the most widely used among various theories of shell structures developed by \cite{simo1989stress}. 
	A shell structure $\boldsymbol{X}(s,t)$ is described by its mid-surface $\boldsymbol{S}(s,t)$ together with its thickness (as shown in Fig.~\ref{fig: Shell structure and its mid-surface}(a)). For any point $\boldsymbol{O}$ on the mid-surface, the local coordinates system $ \{\boldsymbol{v_1},\boldsymbol{v_2},\boldsymbol{v_3}\}$ is constructed as in Fig.~\ref{fig: Shell structure and its mid-surface}(b), where
	$\boldsymbol{v}_3 = [l_3,m_3,n_3]^\T = \frac{\boldsymbol{S}_s\times \boldsymbol{S}_t}{\lvert \boldsymbol{S}_s\times \boldsymbol{S}_t \rvert}$ is the unit normal vector of the point $\boldsymbol{O}$ on the mid-surface.
	According to Reissner-Mindlin theory, the stress alone the direction of $v_3$ is assumed to be zero, i.e.,  $\sigma_{\boldsymbol{v}_3,\boldsymbol{v}_3} = 0$. 
	Once $\boldsymbol{v}_3$ is determined, $\boldsymbol{v}_1$ can be chosen as the outer product of $\boldsymbol{v}_3$ and one of the three axis vectors in the global system $\{\boldsymbol{x},\boldsymbol{y},\boldsymbol{z}\}$, as illustrated in Fig.~\ref{fig: Shell structure and its mid-surface}(b). Subsequently, $\boldsymbol{v}_2$ can be taken as the outer product of $\boldsymbol{v}_3$ and $\boldsymbol{v}_1$
	\citep{kang2015isogeometric}.  

	The shell structure is represented in a NURBS representation as
	\begin{equation} 
		\label{eq: shell}
		\boldsymbol{X}(s,t)=  \boldsymbol{S}(s,t)+ \zeta \boldsymbol{v_3}(s,t) \sum\limits_{i=0}^m\sum\limits_{j=0}^n R_{ij}^{pq}(s,t) \frac{h}{2},
	\end{equation}
	where the mid-surface
	\begin{equation}
		\label{eqn: mid-surface}
		{\boldsymbol S}(s,t)=
		\sum\limits_{i=0}^m\sum\limits_{j=0}^n R_{ij}^{pq}(s,t) \boldsymbol{P}_{ij},  \end{equation}
	$h$ denotes the thickness of the shell structure which is a constant, and $\zeta\in[-1,1]$ is the parameter along with the thickness direction.

	The displacement of a shell structure can be described with $5$ degrees of freedom per control point, consisting of three translation variables and two rotation variables with respect to two tangential directors, denoted by $\boldsymbol{v_1}$ and $\boldsymbol{v_2}$. The displacement vector $\boldsymbol{u}(s,t)$ is represented using the same NURBS basis in Eq.\eqref{eq: shell}
	\begin{equation}
		\label{eq: displacement}
			\boldsymbol{u}(s,t) 
   = \begin{pmatrix}
				u\\v\\w
			\end{pmatrix}
			= \sum\limits_{i=0}^m\sum\limits_{j=0}^n R_{ij}^{pq}(s,t) \begin{pmatrix}
				u_{ij}\\v_{ij}\\w_{ij} \end{pmatrix}  
   + \zeta \boldsymbol{\mu}(s,t) \sum\limits_{i=0}^m\sum\limits_{j=0}^n R_{ij}^{pq}(s,t) \frac{h}{2} \begin{pmatrix}
				\alpha_{ij}\\\beta_{ij}
			\end{pmatrix},
	\end{equation}
	where $(u_{ij}$,$v_{ij}$,$w_{ij})$ are the translations along $\boldsymbol{x}$, $\boldsymbol{y}$, and $\boldsymbol{z}$ directions, respectively. The variables $\alpha_{ij}$ and $\beta_{ij}$ denote the rotation coefficients around the $\boldsymbol{v}_1$-axis and $\boldsymbol{v}_2$-axis, respectively.
	$\boldsymbol{\mu}(s,t) = \left(-\boldsymbol{v_2}(s,t),\boldsymbol{v_1}(s,t)\right)$.
	
	\medskip
	
	\subsection{Material density function}
	Since there is a bijective correspondence between  the mid-surface and  the parametric domain (as shown in Fig.\ref{fig:NURBSSurface}), the material density function of the shell structure 
	Eq.~\eqref{eq: shell}
	can be also defined in NURBS form on 2D parametric domain:
	\begin{equation}
		\label{eq:implicit}
		{\rho}(s,t) = \sum\limits_{i=0}^{m}\sum\limits_{j=0}^{n} R_{ij}^{pq}(s,t) \rho_{ij},
	\end{equation}
	where $\rho_{ij}$ is density control coefficients whose values are in the range of $[0,1]$.
	For each point $(x,y,z)$ on the mid-surface, there is a corresponding point $(s,t) = {\boldsymbol S}^{-1}(x,y,z)$ on the parametric space. 
	Thus the material density function defined on the mid-surface is 
	$ \rho \circ {{\boldsymbol S}}^{-1}$. 
	This is why we call our method IGA based.

	\medskip

	In order to reduce intermediate density, we can get a well-defined material distribution of the structure by the following \textit{\textbf{implicit}} function:
		\begin{equation}
			\left\{
			\begin{aligned}
				{ \rho}(s,t)& = \kappa; \quad boundary\\
				{ \rho}(s,t)& < \kappa; \quad void \\
				{ \rho}(s,t)& > \kappa; \quad solid\\
			\end{aligned}
			\right.
			\label{eq:sur_rep}
		\end{equation}
		To facilitate the numerical implementation, we apply the heaviside function $H(\cdot)$ as discussed  by \cite{wang2011projection} 
		to $\rho(s,t)$ for an approximation of the aforementioned implicit definition:
		\begin{equation}
			\label{eq: heaviside function}
			{\tilde \rho}(s,t) = H({\rho}(s,t)) = \frac{\tanh{\frac{\tau}{2}} + \tanh{\left(\tau ({\rho} - \kappa)\right)}}{2\tanh{\frac{\tau}{2}}},
		\end{equation}
		where $\tau$ is a parameter that controls the sharpness of the projection function. 
		${\tilde \rho}(s,t)$ is regarded as the projected density function of the mid-surface $S(s,t)$.  
	Based on our experience, we typically recommend $\kappa$ to be within the range of $[0.25,0.75]$. Within this range, the optimized structures exhibit similar physical properties.
	In this paper, we set $\kappa = 0.5$.
	
	
        \section{Methods}
	\label{sec:method}
	In this section, we first introduce the formulation of shell structure optimization problem based on IGA-SIMP.
	Then we put forward a model for generating porous shell structures. The numerical method to solve the two models is also presented. 
	
	\subsection{ Formulation of shell structure optimization}
	The objective of this study is to minimize the compliance of a shell structure which satisfies a predefined volume fraction. Based on the NURBS representations of the shell structure (Eq.\eqref{eq: shell}), the displacement vector (Eq.\eqref{eq: displacement}) and the density function (Eq.\eqref{eq:implicit}), the problem is formulated as
	
	\begin{equation}
		\begin{aligned}
			\label{eq: formulation 1}
			(\mathcal{P}) \quad   \min\limits_{\rho_{ij}} &\quad C = \boldsymbol{U^\T F} = \boldsymbol{U^\T K U}, \\
			s.t. & \quad \boldsymbol{KU} = \boldsymbol{F}, \\
			& \quad V \leq V^*, \\
			& \quad  0 \leq \rho_{ij} \leq 1,
		\end{aligned}
	\end{equation}
	where $\{\rho_{ij}\}$ are the density coefficients in Eq.\eqref{eq:implicit} and serve as the design variables, 
	$C$ is the compliance of the shell structure,  $\boldsymbol{K}$ is the stiffness matrix, $\boldsymbol{U}$ is the displacement control vector whose elements are the control coefficients 
	$\{u_{ij},v_{ij},w_{ij},\alpha_{ij},\beta_{ij}\}$ of
	$\boldsymbol{u}(s,t)$ defined in Eq.\eqref{eq: displacement}, 
	$\boldsymbol{F}$ is the equivalent force vector 
	corresponding to the  shell load,
	$\boldsymbol{KU} = \boldsymbol{F}$ is the equilibrium equation, 
	$V$ is the volume of the shell structure, and $V^*$ is a prescribed threshold. 
	
	Notice that the above formulation is a discrete form of the original minimum compliance problem 
	by \cite{bendsoe2003topology} 
	based on the NURBS representations of the shell structure, the displacement vector and the material density function. It has the similar form as its FEA counterpart but is different in the symbols $\boldsymbol{U}$, $\boldsymbol{F}$, $\boldsymbol{K}$
	and the design variables. The details can be found in ~\ref{subsec:analysis}.
	
	The optimization problem ($\mathcal P$) is solved by the Method of Moving Asymptotes (MMA) 
	that will be described in detail in Section~\ref{sec: analysis method}.

	\subsection{Formulation of porous shell structures }
	\begin{figure}[htbp]
		\centering
		\begin{overpic}[width=0.6\textwidth]{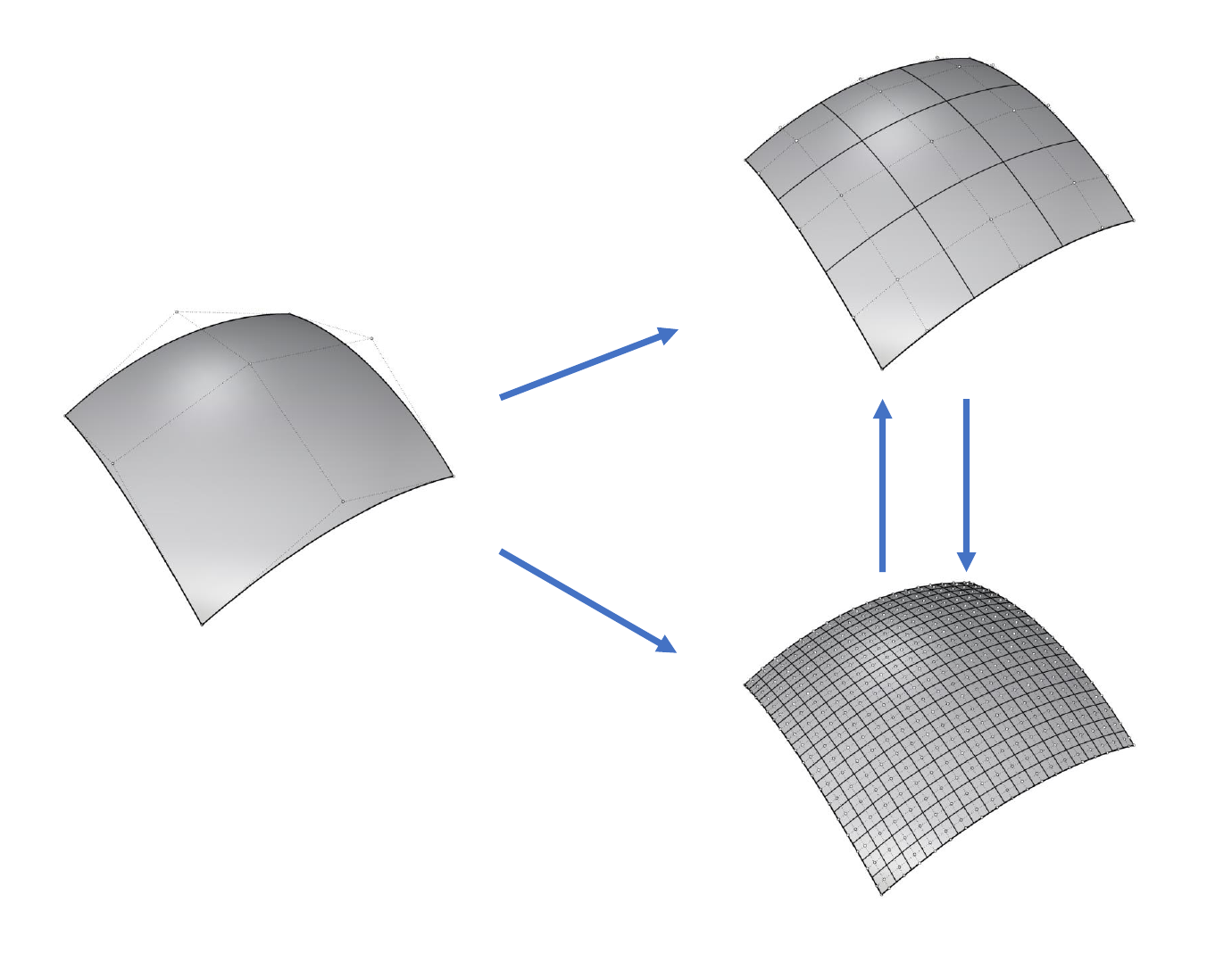}
			\put(10,57){\textbf{CAD model}}
			\put(74,48){\textbf{Design model}}
			\put(44,11){\textbf{Analysis model}}
			\put(0,20){$\boldsymbol{S}(s,t) = \sum\limits_{i = 0}^{m}\sum\limits_{j = 0}^{n} R_{ij}^{pq}(s,t) \boldsymbol{P}_{ij}$}
			\put(37,76){${\rho}(s,t) = \sum\limits_{i = 0}^{m_1}\sum\limits_{j = 0}^{n_1} \bar{R}_{ij}^{pq}(s,t) \bar{\rho}_{ij}$}
			\put(40,2){$\boldsymbol{u_S}(s,t) 
    					= \sum\limits_{i = 0}^{m_2} \sum\limits_{j = 0}^{n_2} \tilde{R}_{ij}^{pq}(s,t) 
                            \begin{pmatrix}
    						\tilde{u}_{ij} \\ \tilde{v}_{ij} \\ \tilde{w}_{ij} 
    					\end{pmatrix}  $
			}
			\put(39,53){\small{{refinement}}}
			\put(80,39){\small{{refinement}}}
			\put(39,25){\small{{refinement}}}
			\put(55,41){\small{coarsen}}
			\put(50,37){\small{(knot removal)}}        
		\end{overpic}
		\caption{Relationship between three models. The mid-surface $\boldsymbol{S}(s,t)$ of a shell structure $\boldsymbol{X}(s,t)$ is refined for different purposes. A coarse refinement of the CAD model is used to define the density function $\boldsymbol{{\rho}}(s,t)$, and a finer refinement is used to analyze the shell displacement $\boldsymbol{u}(s,t)$. $\boldsymbol{u_S}(s,t)$ is a part of $\boldsymbol{u}(s,t)$ which denotes the displacement of the shell mid-surface.}
		\label{fig: Analysis-Design}
	\end{figure}
	\begin{figure*}[htbp]
		\centering
		\begin{overpic}[width=0.8    \textwidth]{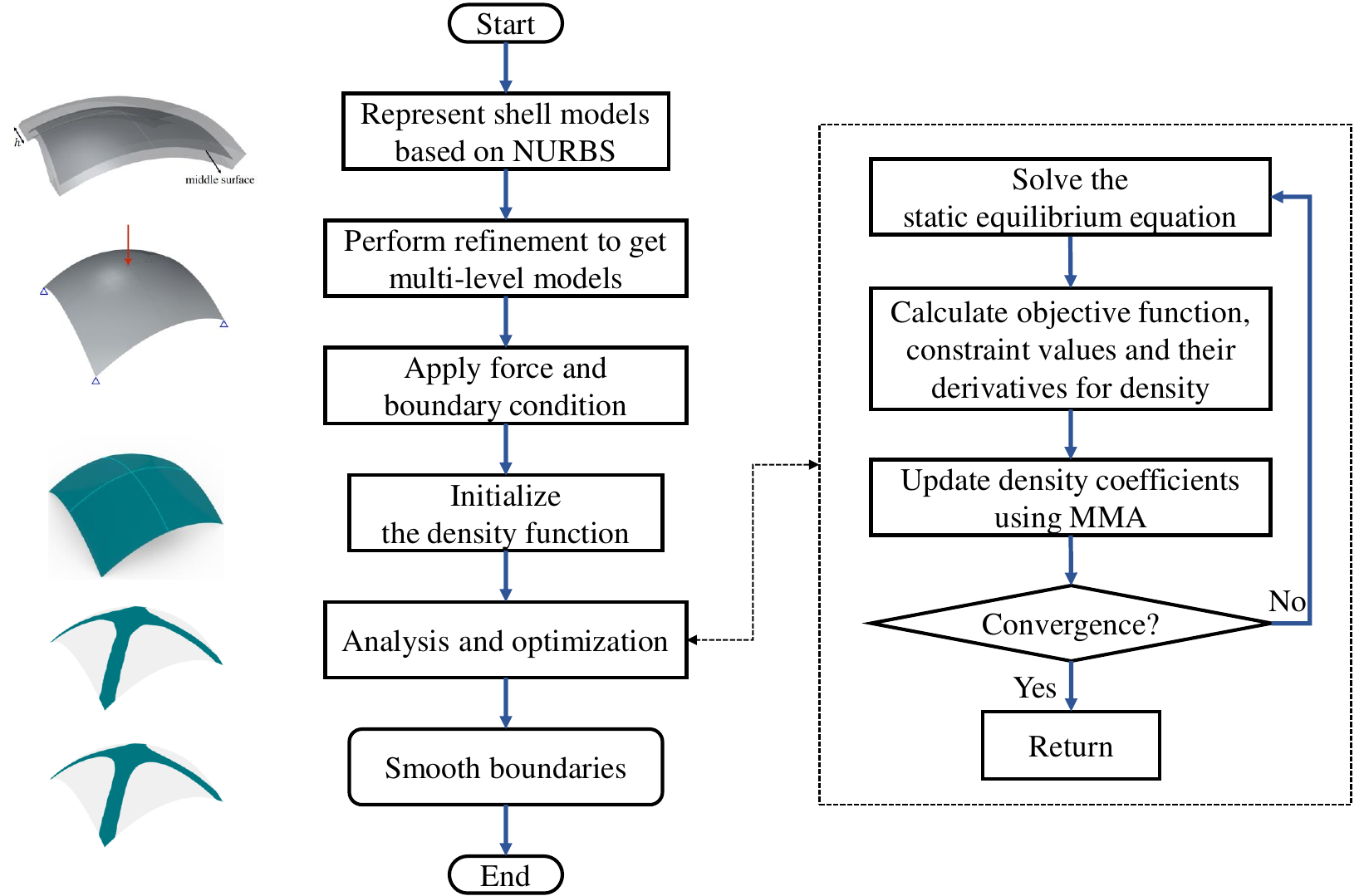}
			\put(10,48){\tiny{$G = 100$}}
		\end{overpic}
		\caption{The pipeline of our algorithm.}
		\label{fig:flowchart}
	\end{figure*}
	
	Porous structures offer some desirable properties such as being lightweight, aesthetic pleasing and high robustness against load perturbations
	\citep{clausen2017topology,dou2020projection}.  
	They can be used not only in medical science and biomimetic technology for repairing or replacing damaged bones
	\citep{balla2010porous,wang2016topological}, 
	but also in daily life to fabricate heat dissipation structures, 
	furniture and adsorbents 
	\citep{davis2002ordered,das2020multi,hu2023parametric}.  
	Further studies on porous structures can be found in various literature sources 
	\citep{wang2013cost, alexandersen2015topology,wu2017infill,liu2017additive,liu2018porous, Khabazi2011}.  
	
	By replacing the global volume constraint in the problem (Eq.~\eqref{eq: formulation 1}) with a set of local volume constraints, 
	we propose a new model to generate porous shell structures with good properties. To be specific, the problem is formulated as 
	\begin{equation}
		\begin{aligned}
			\label{eq: formulation 2}
			(\mathcal{Q}) \quad  \min\limits_{\rho_{ij}} &\quad C = \boldsymbol{U^\T F} = \boldsymbol{U^\T K U}, \\
			s.t.&\quad \boldsymbol{KU} = \boldsymbol{F}, \\
			&\quad  \Bar{\rho}_e \leq {\alpha}, \quad e=1,2,\ldots, N_e \\
			& \quad 0 \leq \rho_{ij} \leq 1, 
		\end{aligned}
	\end{equation}
	where $C$, $\boldsymbol{U, F, K}$, $\rho_{ij}$ are the same symbols as in the previous formulation~\eqref{eq: formulation 1},  
	$\alpha \in (0,1)$ is a prescribed upper bound, and 
	\begin{equation}
		\label{eq: local volume}
		\Bar{\rho}_e = \frac{\sum\limits_{f \in \mathcal{N}_e} \tilde\rho_f V_f }{\sum\limits_{f \in \mathcal{N}_e} V_f} 
	\end{equation}
	measures the local material density of an element $e$ which is the average of the material densities of elements in the neighborhood 
	$\mathcal{N}_e = \{f| d(x_e,x_f) \leq \mathcal{R} \cdot \delta \}$ of the element $e$. $\mathcal{R}$ is an integer and $\delta$ is the average element length. 
	$V_f$ denotes the volume of the element $f$, and $N_E$ is the number of elements of the shell structure. 
	It is worth noting that the element $e$ represents a solid shell element enclosed by the offset surfaces of the curve elements in the mid-surface. The curved elements on the mid-surface are obtained from isoparametric curves under different knots.
	
	The element-wise local volume constraints $\{\Bar{\rho}_e \leq {\alpha}\}$ are almost impossible to satisfy, so we relax them with one global constraint 
	\begin{equation}
		\label{gbcons}
		\Bar{V} = \left(\frac{1}{N_E}\sum_{e=1}^{N_E} (\Bar{\rho}_e)^{\gamma}\right)^{\frac{1}{\gamma}} \le \alpha,
	\end{equation}
	which is an aggregation of individual constraints for each element and is an approximation of the following maximum constraint
	\begin{equation*} 
		\max\limits_{1\le e \le N_E} \Bar{\rho}_e
		\leq {\alpha}. 
	\end{equation*}
	Again the optimization problem $(\mathcal{Q})$ with the global constraint \eqref{gbcons} 
	is solved by the MMA, and the sensitivity analysis is described in ~\ref{subsec:sensitivity analysis}.

	\subsection{Numerical method}
	\label{sec: analysis method}
	The optimization problems ($\mathcal P$) and ($\mathcal Q$) are solved using MMA.
	MMA converts the original problem into successive convex optimization subproblems. After initializing the design variables
	$\rho_{ij}$ (in our implementation, we set $\rho_{ij}=V^*/V_s$, where $V_s$ represents the volume of the design domain for shell structures), in each iteration, a convex quadratic programming problem is solved, where the key is to solve the updated equilibrium equation $\boldsymbol{K}\boldsymbol{U}=\boldsymbol{F}$.  
	In our IGA-SIMP framework, the shell structure can be divided into curved elements according to the knot spans of NURBS.
	For an element $e$, the displacement control vector, the equivalent force vector, and the stiffness matrix are denoted as $\boldsymbol{U}_e$, $\boldsymbol{F}_e$, and $\boldsymbol{K}_e$, respectively. 
	Accordingly, the compliance $C$ is assembled as
	\begin{equation}
		\label{eq: compliance}
		C = \sum\limits_{e=1}^{N_E} \boldsymbol{U}_e^\T \boldsymbol{F}_e  = \sum\limits_{e=1}^{N_E} \boldsymbol{U}_e^\T \boldsymbol{K}_e \boldsymbol{U}_e,
	\end{equation}
	and the volume is computed as
	\begin{equation}
		\label{eq: volume}
		V = \sum\limits_{e=1}^{N_E} V_e = \sum\limits_{e=1}^{N_E} V_e^0 \cdot \tilde\rho_e,
	\end{equation}
	where $V_e^0$ denotes the solid volume of element $e$, and $\tilde\rho_e$ 
	is the material density at the center of element $e$. The sensitivity analysis for the gradient computation of $C$ and $V$ with respect to $\rho_{ij}$ is provided in ~\ref{subsec:sensitivity analysis}.

	The stiffness matrix $\boldsymbol{K}$ is assembled by element-wise $\boldsymbol{K}_e$
	which is calculated as an integral over the parametric domain.
	Details are presented in ~\ref{subsec:analysis}. 
	During each iteration, the force vector $\boldsymbol{F}$ and the stiffness matrix $\boldsymbol{K}_e^0$ for a solid element remain fixed and can be computed in advance. The density-related stiffness matrix $\boldsymbol{K}_e$ is recomputed according to the update of the design variables $\rho_{ij}$.
	
	\medskip
	
	For a given shell structure represented by Eq.\eqref{eq: shell}, it has to be refined for design and analysis purposes. To represent the material distribution function $\rho(s,t)$ accurately, we refine the NURBS representation of the mid-surface Eq.\eqref{eqn: mid-surface} by knot insertion, and consequently the material density function $\rho(s,t)$ is refined accordingly:
	\begin{equation}
		\rho(s,t)=\sum\limits_{i=0}^{m_1}\sum\limits_{j=0}^{n_1} {\bar R}_{ij}^{pq}(s,t) {\bar\rho}_{ij}.
	\end{equation}
	Similarly, to solve the equilibrium equation with high accuracy, 
	the NURBS representation of the shell structure is refined, and accordingly
	the displacement $\boldsymbol{u}$ is refined as
	\begin{equation}
		\begin{aligned}
			\boldsymbol{u}(s,t) 
			&= \begin{pmatrix}
				u\\v\\w
			\end{pmatrix}
			= \sum\limits_{i=0}^{m_2}\sum\limits_{j=0}^{n_2} {\tilde R}_{ij}^{pq}(s,t) \begin{pmatrix}
				{\tilde u}_{ij}\\{\tilde v}_{ij}\\{\tilde w}_{ij} \end{pmatrix}  \\
			&+ \zeta \boldsymbol{\mu}(s,t) \sum\limits_{i=0}^{m_2}\sum\limits_{j=0}^{n_2} {\tilde R}_{ij}^{pq}(s,t) \frac{h}{2} \begin{pmatrix}
				{\tilde \alpha}_{ij}\\ {\tilde \beta}_{ij}
			\end{pmatrix}.
		\end{aligned}
	\end{equation}
	${\bar R}_{ij}^{pq}(s,t)$ and ${\tilde R}_{ij}^{pq}(s,t)$ are the corresponding refined NURBS basis functions.
	
	Therefore, the shell structure has three levels of models--the original model, the design model that is used to define the density function of the material, and the analysis model that is applied to solve the displacement vector. 
	In practice, we generally have $m_2>m_1>m$ and $n_2>n_1>n$, that is, the design model is generally coarse, while the analysis model is fine. 
	Fig.~\ref{fig: Analysis-Design} illustrates the relationships between the three different models. The complete flowchart of our optimization approach is shown in Fig.\ref{fig:flowchart}.
	
	
        \section{Boundary fairing}
	\label{sec:post-processing}
	
	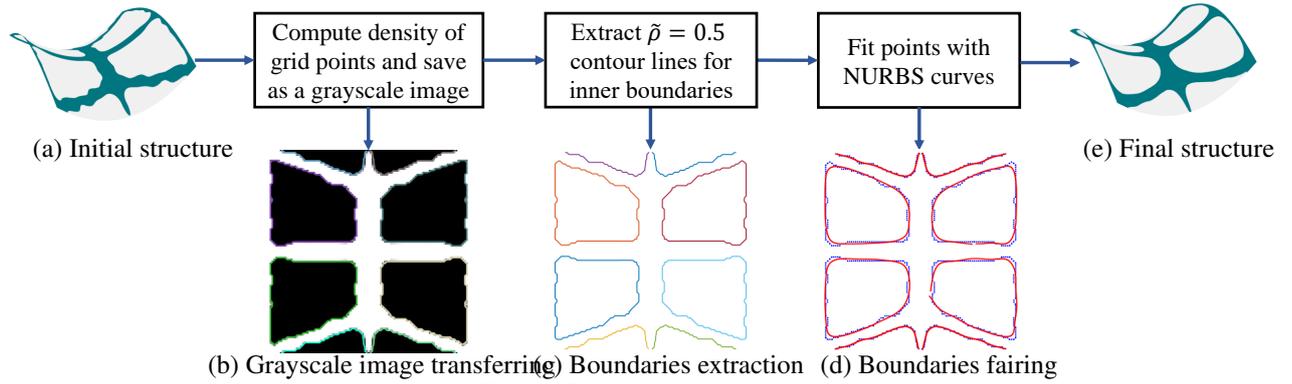
\begin{figure*}[htbp]
		\centering
		\begin{overpic}[width=1\textwidth]{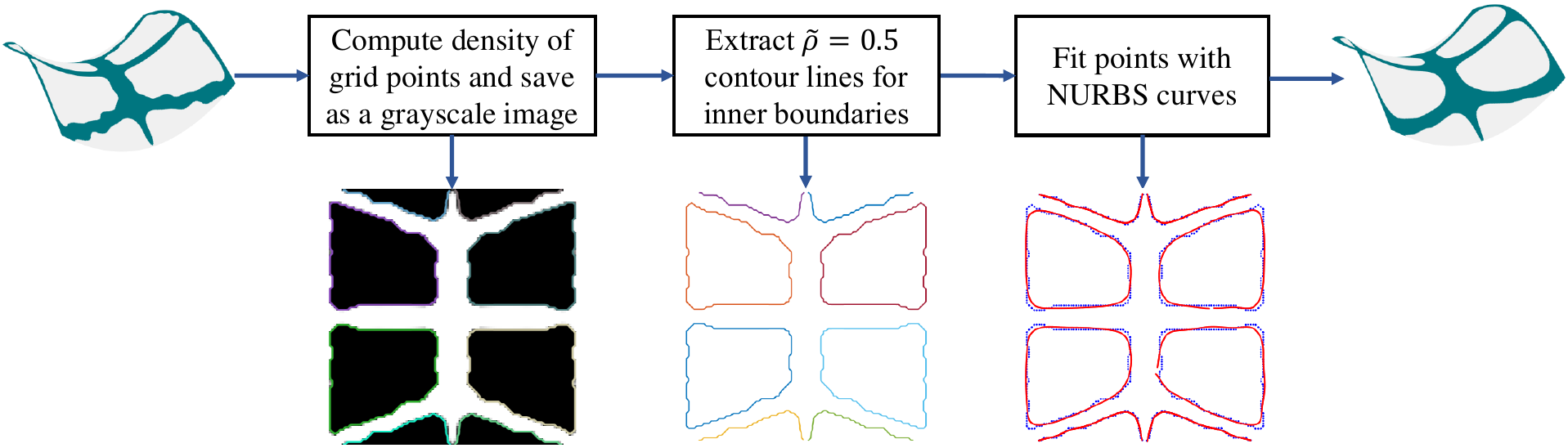}
			\put(2,16){(a) Initial structure}
			\put(86,16){(e) Final structure}
			\put(16,-1.5){(b) Grayscale image transferring}
			\put(42,-1.5){(c) Boundaries extraction}
			\put(65,-1.5){(d) Boundaries fairing}
		\end{overpic}
		\caption{Boundary fairing process.} 
		\label{fig:smoothing process}
	\end{figure*}

	Though the density function $\tilde{\rho}(s,t)$ in Eq.\eqref{eq: heaviside function} is represented by a smooth function,  
	the borders of the holes are generally wavy due to significant differences in density coefficients between neighboring indices.  To generate smooth and fair boundaries, which not only facilitate the manufacture~\citep{sarioz2006optimization} but also enhance the aesthetic appearance of the structures~\citep{farin1987fairing,farin1989curvature,sapidis1990automatic}, a post-processing step becomes essential.
	We propose an automatic approach to obtain fairing boundaries by fitting B-spline curves, which is a crucial operation in Computer-Aided Geometric Design (CAGD). It aims to remove irregularities and noise from curves while preserving their essential shape characteristics.
	
	We first divide the parametric domain into a rectangular grid and calculate density values at nodes. These values are then stored as a matrix, visualized as a grayscale image shown in Fig.~\ref{fig:smoothing process} (b).
	The Marching Squares algorithm \citep{maple2003geometric} is employed to extract the contour lines representing the boundaries in the parametric domain shown in Fig.~\ref{fig:smoothing process} (c). Specifically, the contour lines are identified with their value equal to 0.5. 
	Then, we obtain multiple sets of ordered points along with their corresponding indices from the entire grid point set. Each set represents an inner boundary curve.
	Finally, for each inner boundary, we extract the boundary points $\left\{\boldsymbol{Q_l}\right\}_{l=0}^{c}$ and fit them with fair B-spline curves~\citep{piegl1996nurbs, sarioz2006optimization} in Fig.~\ref{fig:smoothing process} (d).

	The fitting curve $\boldsymbol{Z}(\xi)$ is a NURBS curve of degree 3 defined on uniform knots $\boldsymbol{\Xi} = \{\xi_0,\xi_1, \cdots, \xi_{b+4}\}$ with $b+1$ control points: 
	\begin{equation}
		\label{eq: boundary curve}
		\boldsymbol{Z}(\xi) = \sum\limits_{k=0}^{b} N_{k}^{3}(\xi) \boldsymbol{P}_{k},\quad \xi\in [0,1]. 
	\end{equation}
	
	The point set $\left\{\boldsymbol{Q}_l\right\}_{l=0}^{c}$ is approximated by $\boldsymbol{Z}(\xi)$ in the least squares sense with an additional fairing term in integral form
	\begin{equation} 
		\label{eq: fitting}
		\begin{aligned}
			& \min_{\left\{\boldsymbol{P}_k\right\}_{k=0}^{b}}\  \left\{ \sum\limits_{l=0}^{c} [\boldsymbol{Z}(\eta_l) - \boldsymbol{Q}_l]^2 + \lambda \int\limits_{0}^1 \left| \frac{d^{2}\boldsymbol{Z}(\xi)}{d \xi^2}  \right|^2 d\xi \right\} \\
			= 
			& \min_{\left\{\boldsymbol{P}_k\right\}_{k=0}^{b}}\  \left\{ \sum\limits_{l=0}^{c} [\boldsymbol{Z}(\eta_l) - \boldsymbol{Q}_l]^2 + \lambda \sum\limits_{k=3}^{b}\int\limits_{\xi_{l}}^{\xi_{l+1}} \left|\frac{d^{2}\boldsymbol{Z}(\xi)}{d \xi^2} \right|^2 d\xi \right\},
		\end{aligned}
	\end{equation}
	where $\lambda$ is a positive weight that affects the fairness of the curve. The larger the $\lambda$, the more fair the curve. 
	$\left\{\eta_i\right\}_{l=0}^{c}$ are precomputed parameter values of $\left\{\boldsymbol{Q}_l\right\}_{l=0}^{c}$ obtained using centripetal parameterization \citep{lee1989choosing}.  
	Note that the knot vectors and control points for open and closed curves have different selection criteria, for which additional details can be found in \cite{hoschek1993fundamentals}. 
	
	The fairing term in Eq.\eqref{eq: fitting} can be approximated by two-point Gaussian integration and simplified as follows 
	\begin{equation}
		\sum\limits_{k=3}^{b} \left\{\frac{\xi_k+\xi_{k+1}}{2} \left[ \omega_1  \frac{d^{2}\boldsymbol{Z}(\xi_{k_1})}{d \xi^2} + \omega_2  \frac{d^{2}\boldsymbol{Z}(\xi_{k_2})}{d \xi^2} \right] \right\}^2,
	\end{equation}
	where $\omega_1 = \omega_2 = 1$ and $\xi_{k_1},\xi_{k_2}$ are Gauss weights and Gauss points in $[\xi_k,\xi_{k+1}]$ respectively.
	
	By setting the partial derivatives of the objective function with respect to the control points $\{\boldsymbol{P}_{k} \}$ to be zero, the optimization problem Eq.\eqref{eq: fitting} can be transformed into a linear system
	\begin{equation} 
		\label{eq: linear system}
		\begin{aligned}
			\sum\limits_{l=0}^{c} N_j^3(\eta_l) \boldsymbol{Q}_l 
			& = \lambda \sum\limits_{i=0}^{b} \left( \sum\limits_{k=3}^{b} \frac{(\xi_k+\xi_{k+1})^2}{4} \cdot \left( \frac{d^2 N_{i}^3 (\xi_{k_1})}{d \xi^2} \right. \right. \\
			& + \left. \left. \frac{d^2 N_{i}^3 (\xi_{k_2})}{d \xi^2} \right) \left( \frac{d^2 N_j^3 (\xi_{k_1})}{d \xi^2} +  \frac{d^2 N_j^3 (\xi_{k_2})}{d \xi^2} \right)  \right) \boldsymbol{P}_{i} \\
			& + \sum\limits_{i=0}^{b} \left( \sum\limits_{k=0}^{b} N_{i}^3(\eta_k) N_j^3(\eta_k) \right) \boldsymbol{P}_{i},
		\end{aligned}
	\end{equation}
	for $j = 0,1,\cdots,b$. 
	$\lambda = 0.01$.
	Then, an explicit representation of a fairing boundary curve is obtained. 
	The integration of the fairing procedure with shell structures enables efficient editing to meet users' requirements and manufacturing standards. The entire process is summarized in algorithm~\ref{algorithm 1}.

	\begin{algorithm}
		\caption{Boundary fairing algorithm}
		\label{algorithm 1}
		\begin{algorithmic}[1]
			\Require 
			$\tilde\rho(s,t)$: Projected density  function as Eqs.\eqref{eq:implicit}-\eqref{eq: heaviside function}; \newline
			$b$: number of control points for curves \newline
			$\lambda$: fairing factor;
			\Ensure  
			$\boldsymbol{Z}(\xi)$;
			\State Grid the parametric domain and calculate the density values on grid points by $\tilde\rho(s,t)$;
			\State Extract $\tilde\rho=0.5$ for inner boundaries with ordered points through Marching Squares algorithm;
			\Repeat 
			\State parameterize for $\{\boldsymbol{Q}_l\}_{l=1}^c$;  
			\State Solve Eq.\eqref{eq: linear system} for $\{\boldsymbol{P}_k\}_{k=0}^b$;
			\Until{fitting is finished for all point clusters} 
			\State Compute $\rho(s,t)$  in physical domain via $S(s,t)$\eqref{eq:explicit} and all $\boldsymbol{Z}(\xi)$ in parametric space.
		\end{algorithmic}
	\end{algorithm}

	
        \section{Numerical experiments and discussion}
	\label{sec:results}
	This section presents a series of numerical examples to showcase the effectiveness of our proposed method. 
	Firstly, we conduct a comparison between FEA-SIMP and IGA-SIMP to establish the validity and superiority of our approach. 
	Subsequently, we apply the IGA-SIMP model to generate porous shells.
	We then present the fairing results of the optimized shell structures.
	The selection of the number of design variables is discussed and we compare our method with several other IGA-based methods in the end.

	The parameter settings for the simulation are as follows: the Young's modulus is $E = 2100$, and the Poisson's ratio is $\nu = 0.3$. The thickness of the shell considered in Eq.~\eqref{eq: shell} is set to $h=5$. For the iterative process, the value of $\tau$ in Eq.~\eqref{eq: heaviside function} begins at $2$ and is doubled every $25$ iterations until it reaches $64$. 
	The optimization process terminates if either the maximum change of nodal material densities is less than 0.005 in five consecutive iterations or 200 iteration steps have been completed. The default values of parameters in MMA are $move = 0.1, asyinit = 0.1, asyincr = 1.1, asydecr = 0.7$.

	\subsection{Comparison with FEA-SIMP}
	\begin{figure*}[t]
		\subfigure[\label{fig:initial} initial]
		{
			\begin{minipage}[t]{0.12\linewidth}
				\centering
				\begin{overpic}[scale = 0.08]{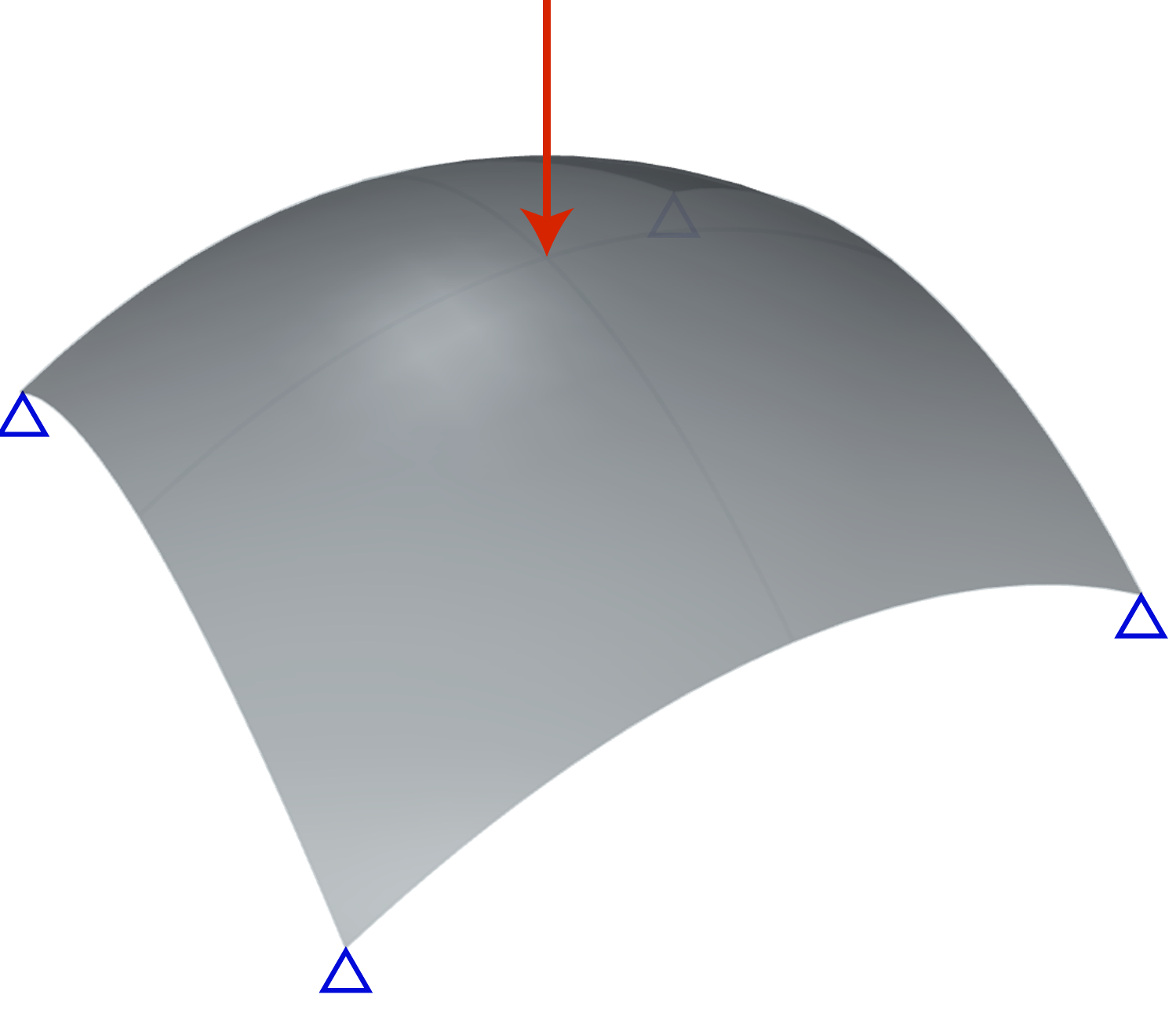}
					\put(50,75){\tiny{$G = 100$}}
				\end{overpic}
				\\
				\vspace{0.1cm}
				\begin{overpic}[scale = 0.08]{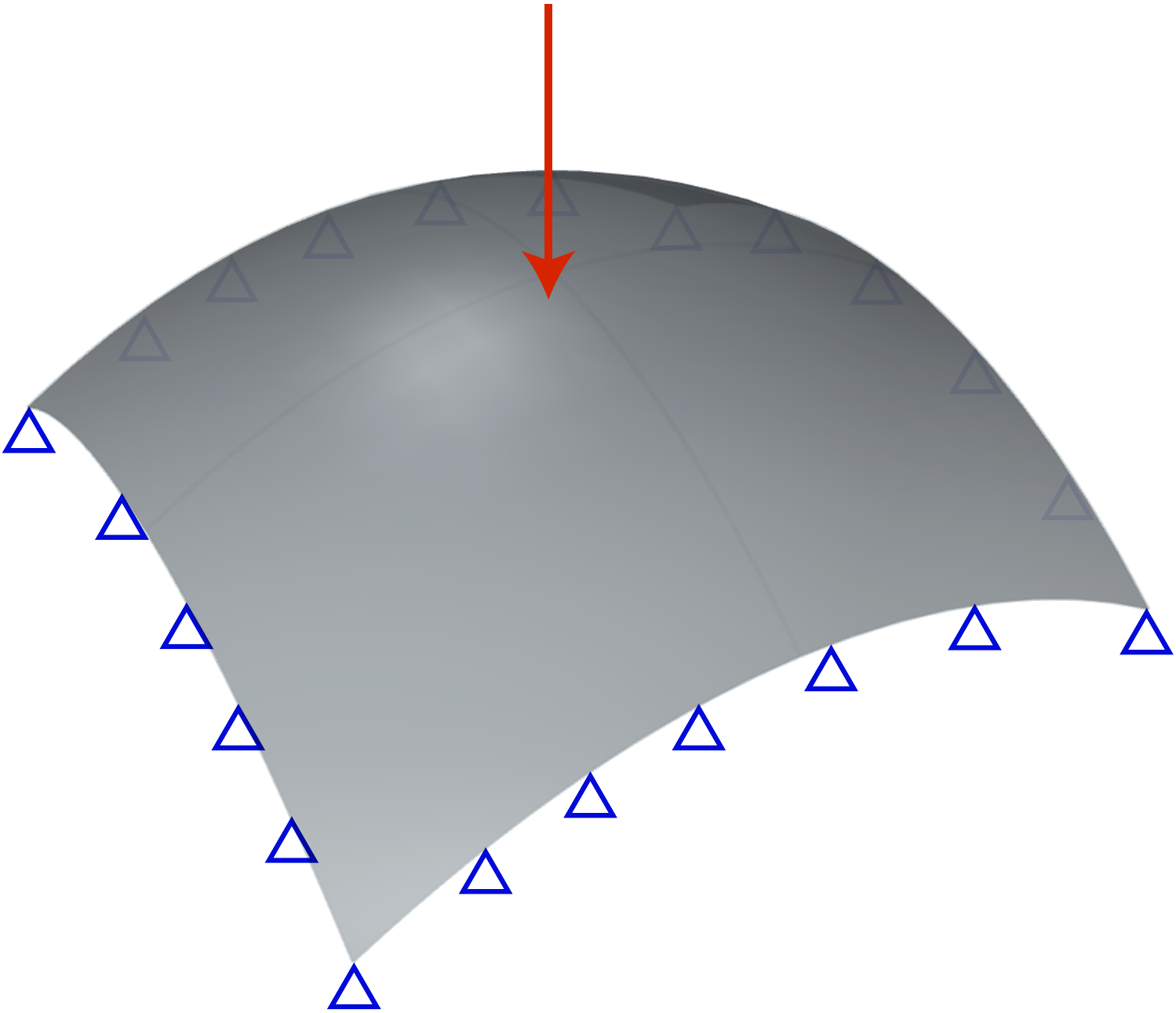}
					\put(50,75){\tiny{$G = 100$}}
				\end{overpic}
				\\
				\vspace{0.1cm}
				\begin{overpic}[scale = 0.08]{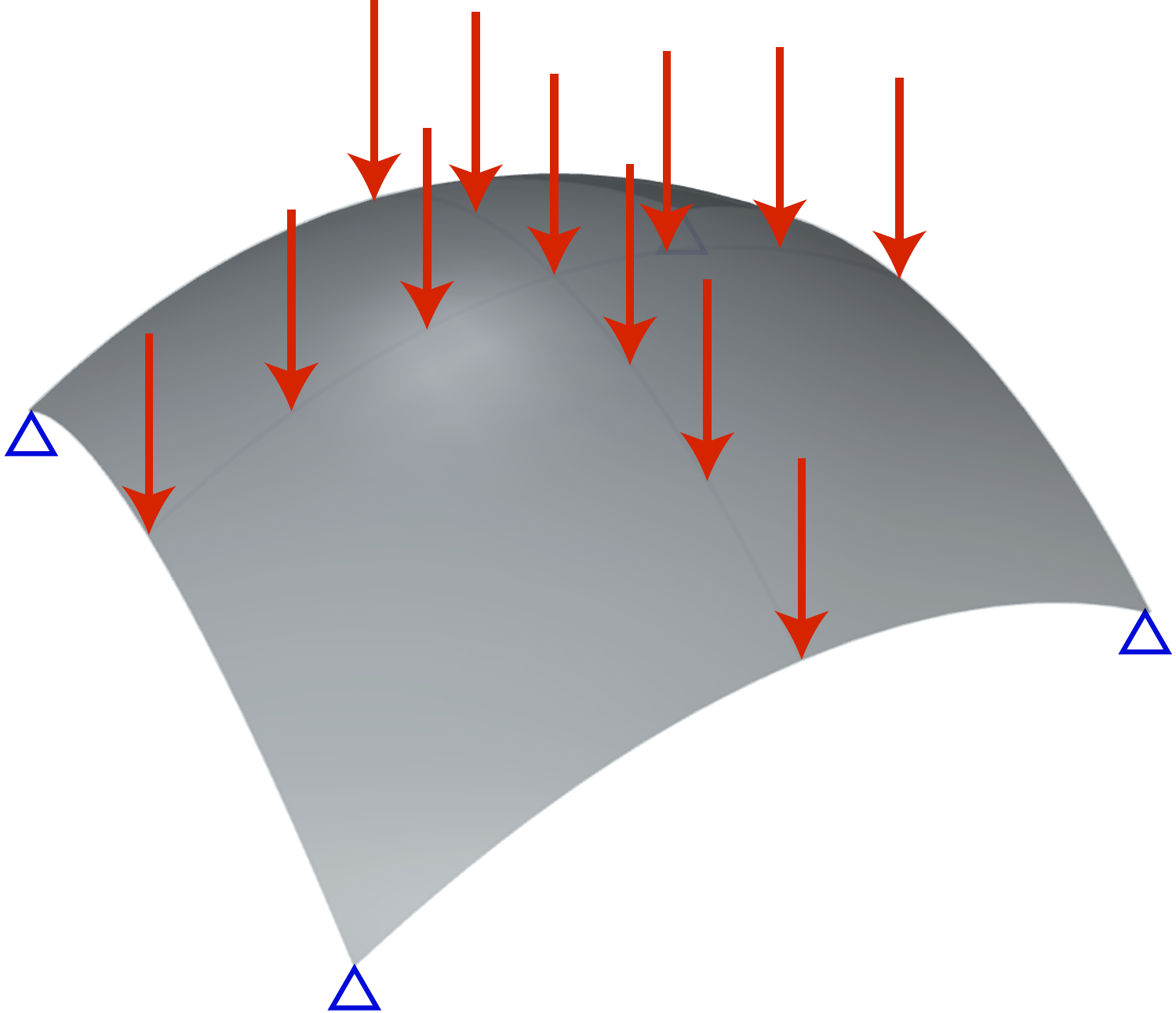}
					\put(38,83){\tiny{$G = 10$}}
				\end{overpic}
				\\        
				\vspace{0.2cm}
			\end{minipage}
		}
		\subfigure[\label{fig:IGA_50 front} {IGA50}]{
			\begin{minipage}[t]{0.125\linewidth}            
				\centering
				\includegraphics[scale = 0.1]{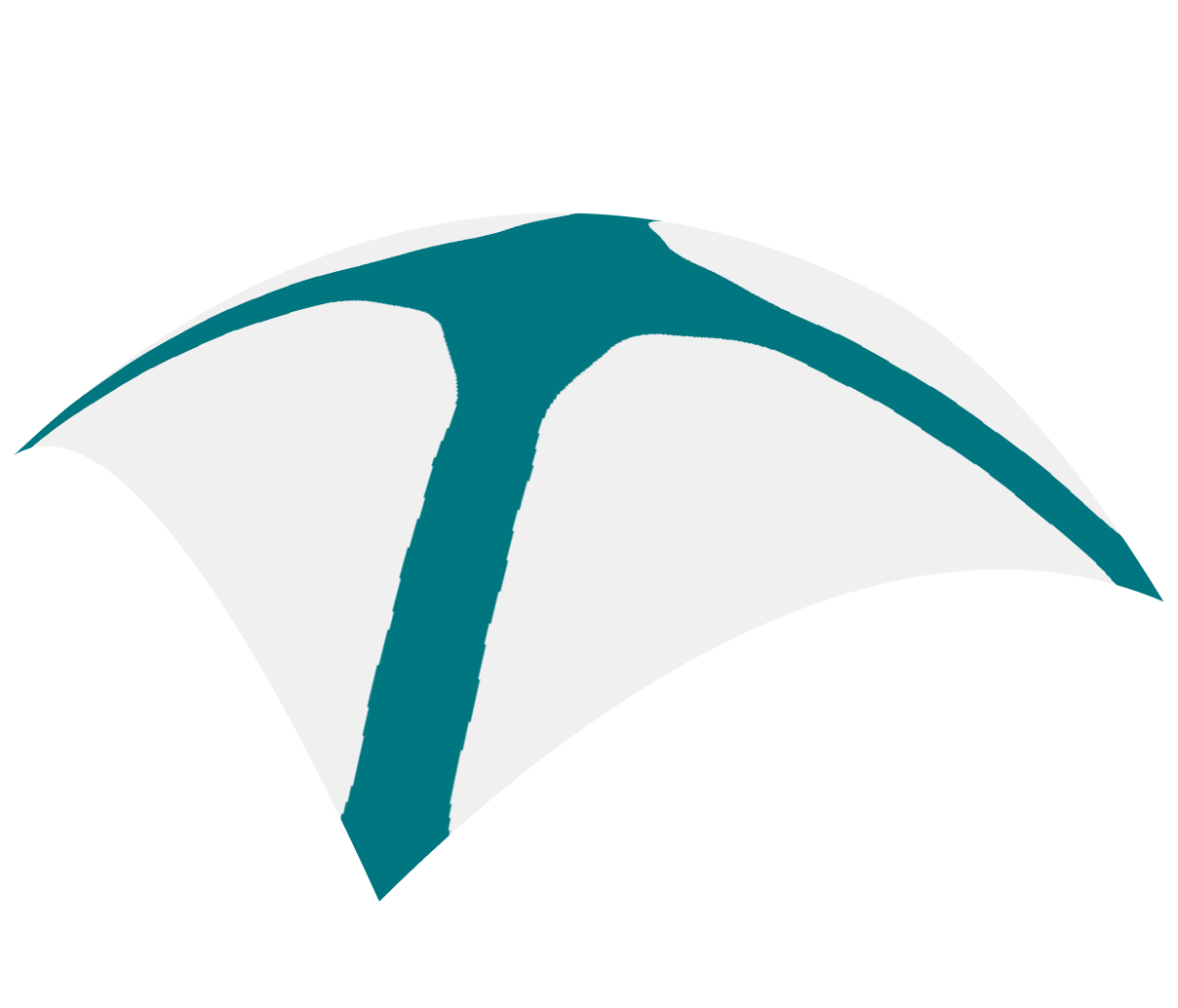}\\
				\vspace{0.1cm}
				\includegraphics[scale = 0.1]{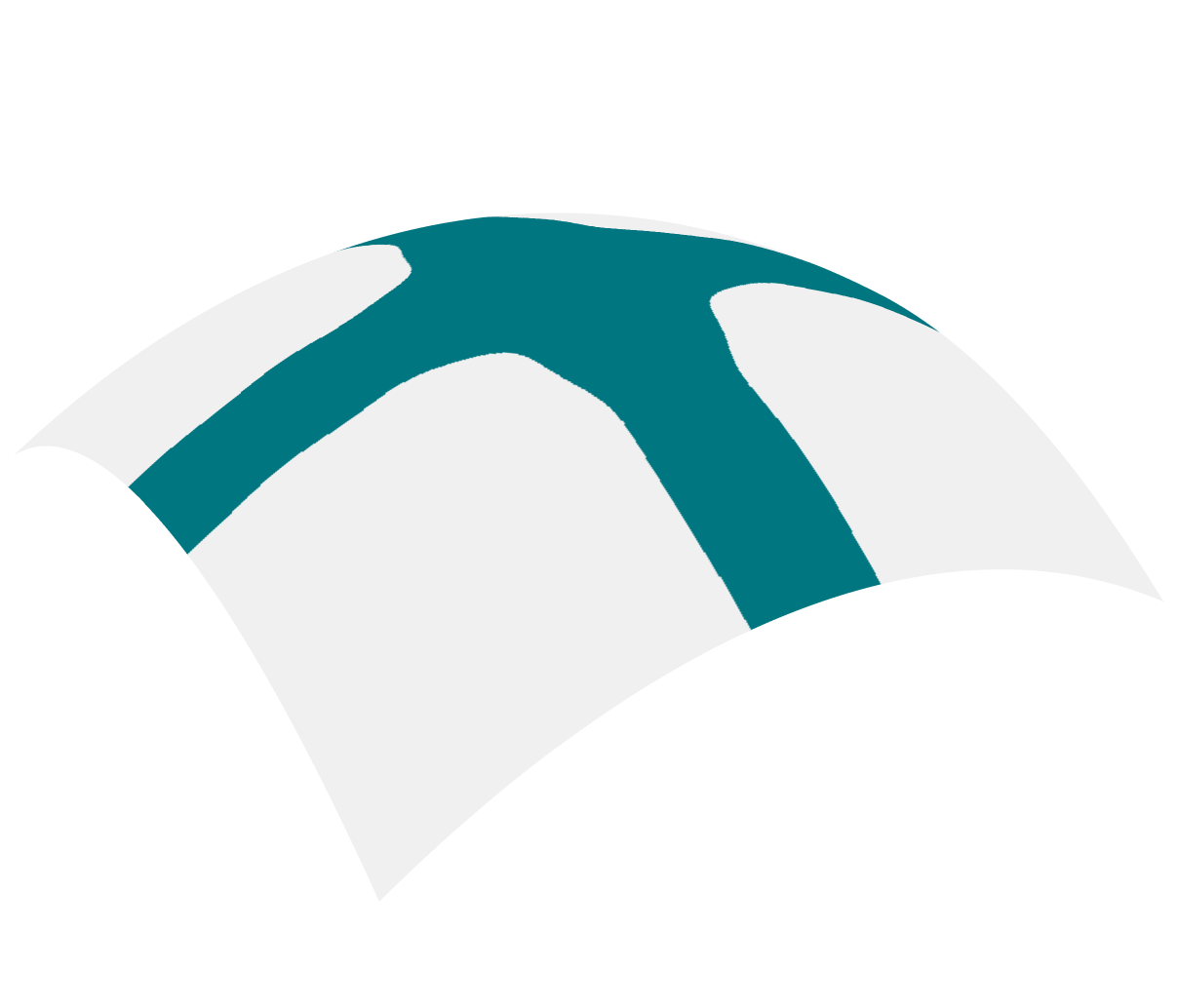}\\
				\vspace{0.1cm}
				\includegraphics[scale = 0.1]{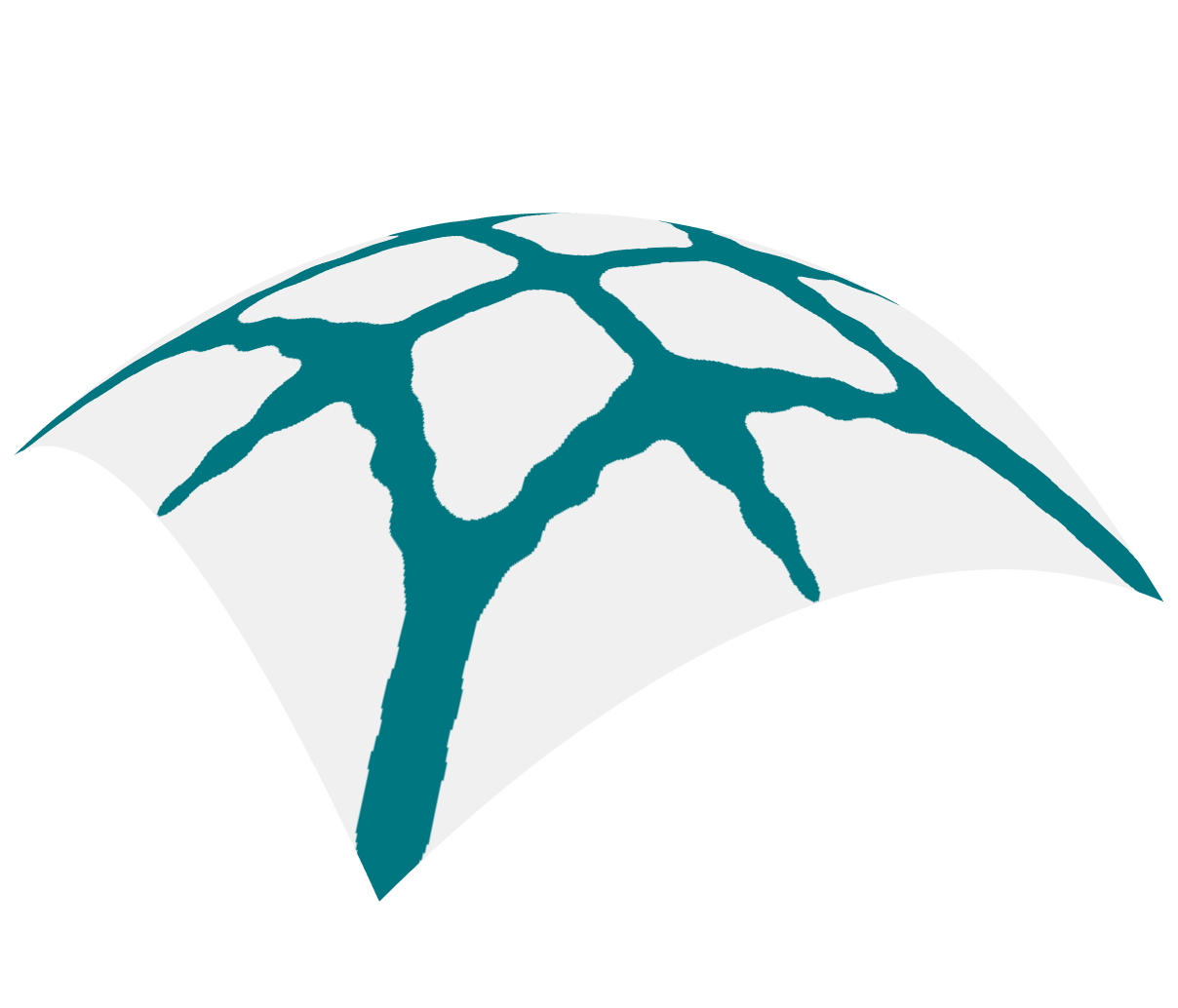}\\
				\vspace{0.2cm}
			\end{minipage}
		}
		\subfigure[\label{fig:IGA_50 top} top view]{
			\begin{minipage}[t]{0.125\linewidth}
				\centering
				\includegraphics[scale = 0.1]{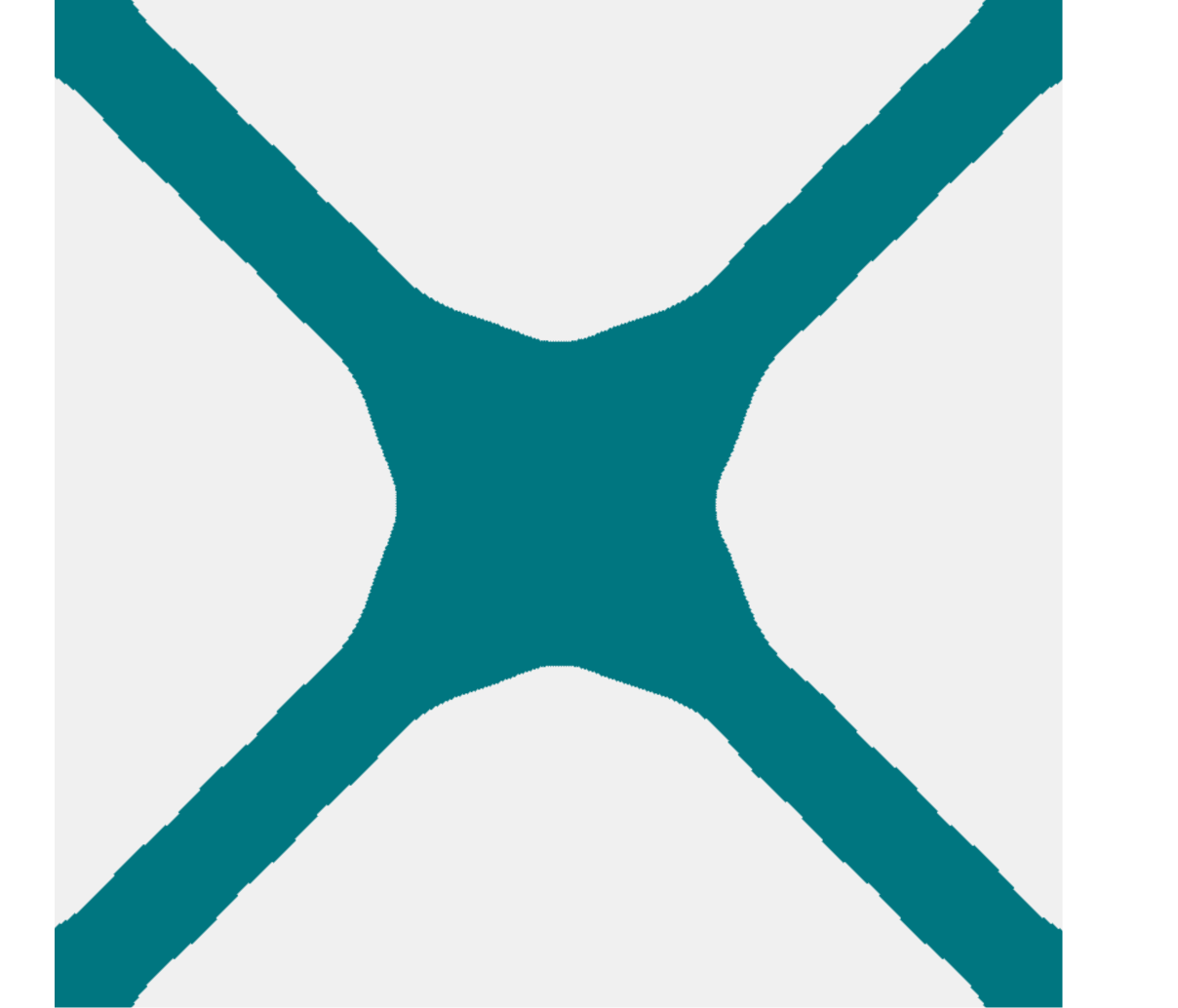}\\
				\vspace{0.1cm}
				\includegraphics[scale = 0.1]{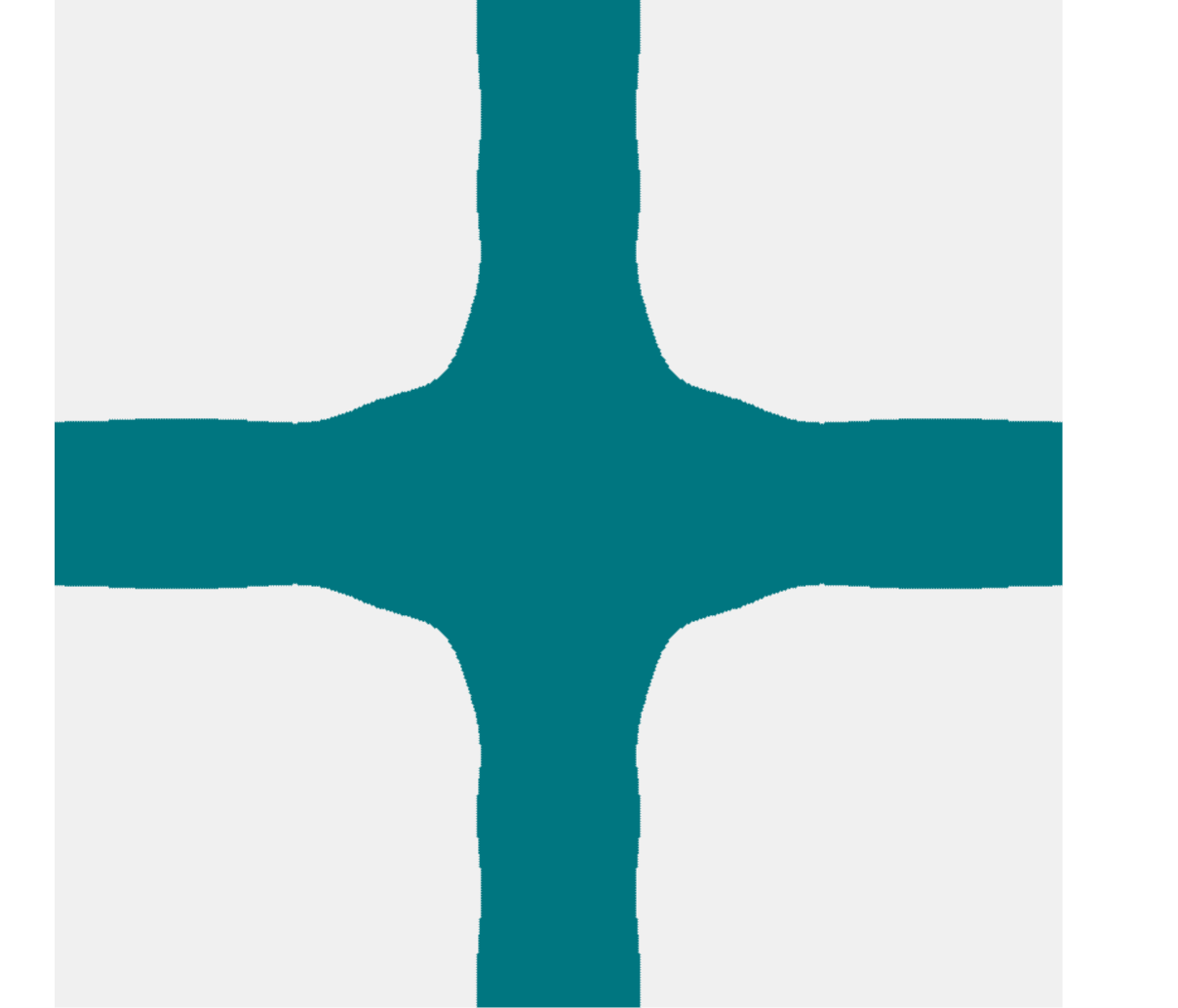}\\
				\vspace{0.1cm}
				\includegraphics[scale = 0.1]{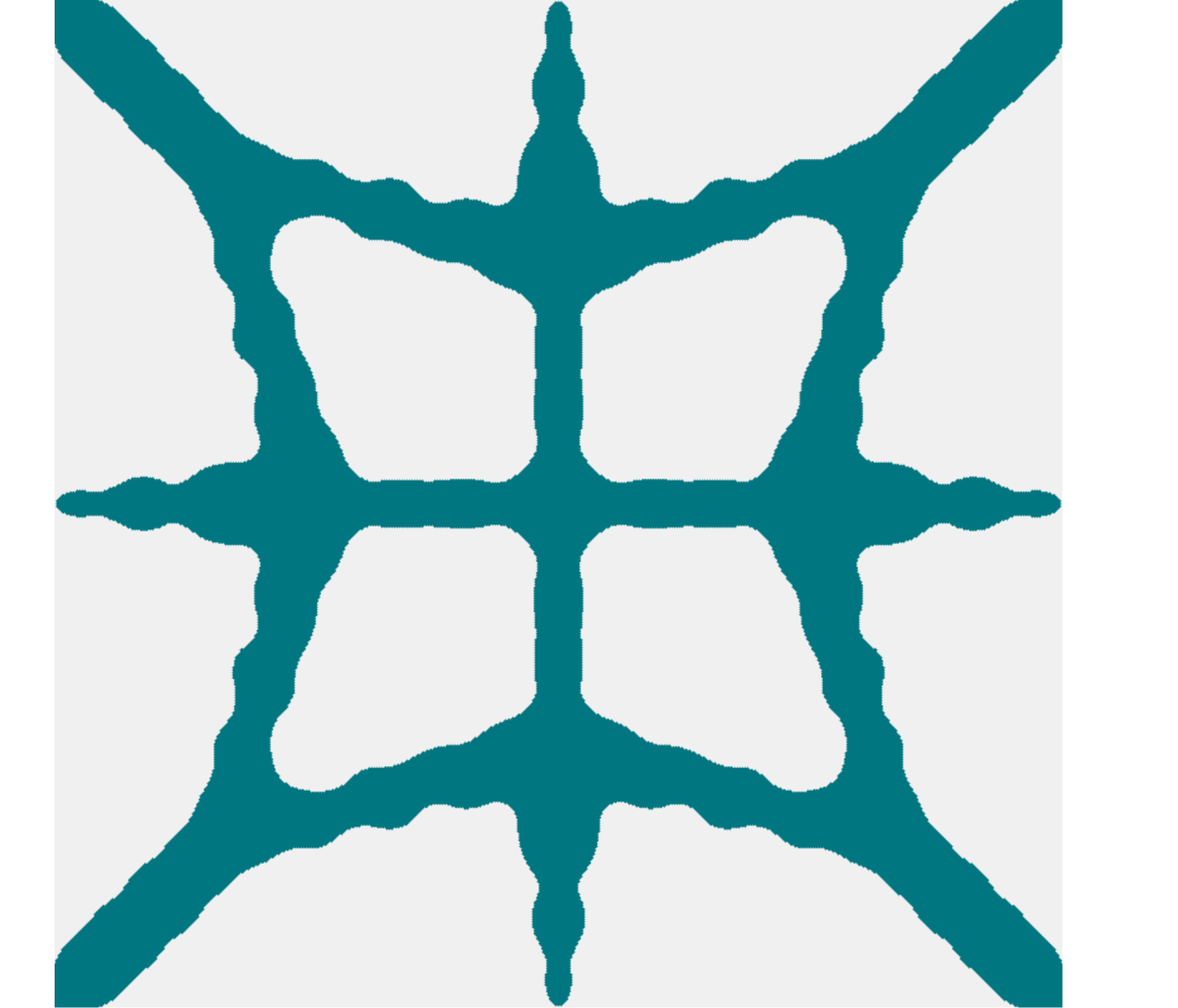}\\
				\vspace{0.2cm}
			\end{minipage}
		}
		\subfigure[\label{fig:IGA_100 front} IGA100]{
			\begin{minipage}[t]{0.125\linewidth}
				\centering
				\includegraphics[scale = 0.1]{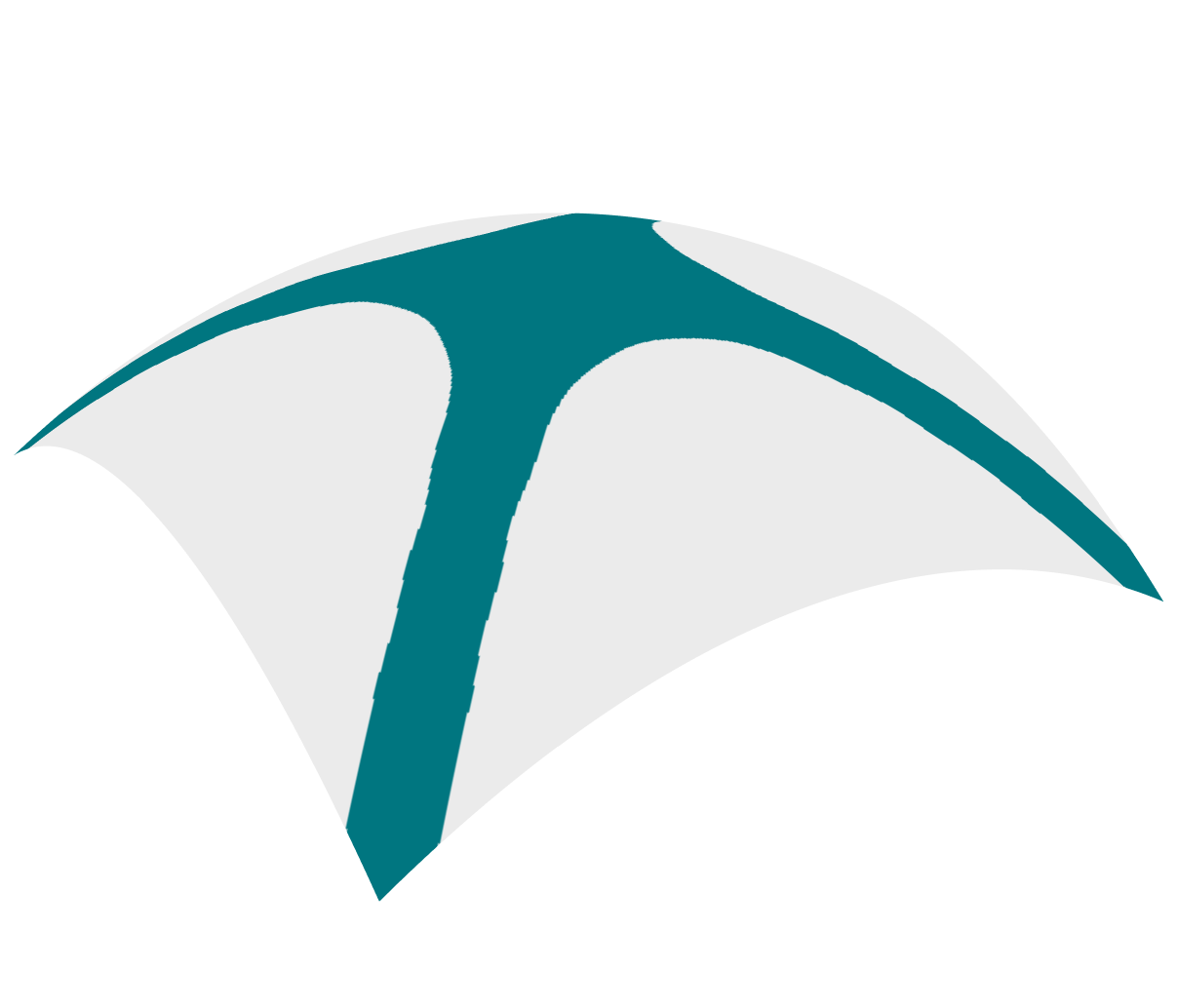}\\
				\vspace{0.1cm}
				\includegraphics[scale = 0.1]{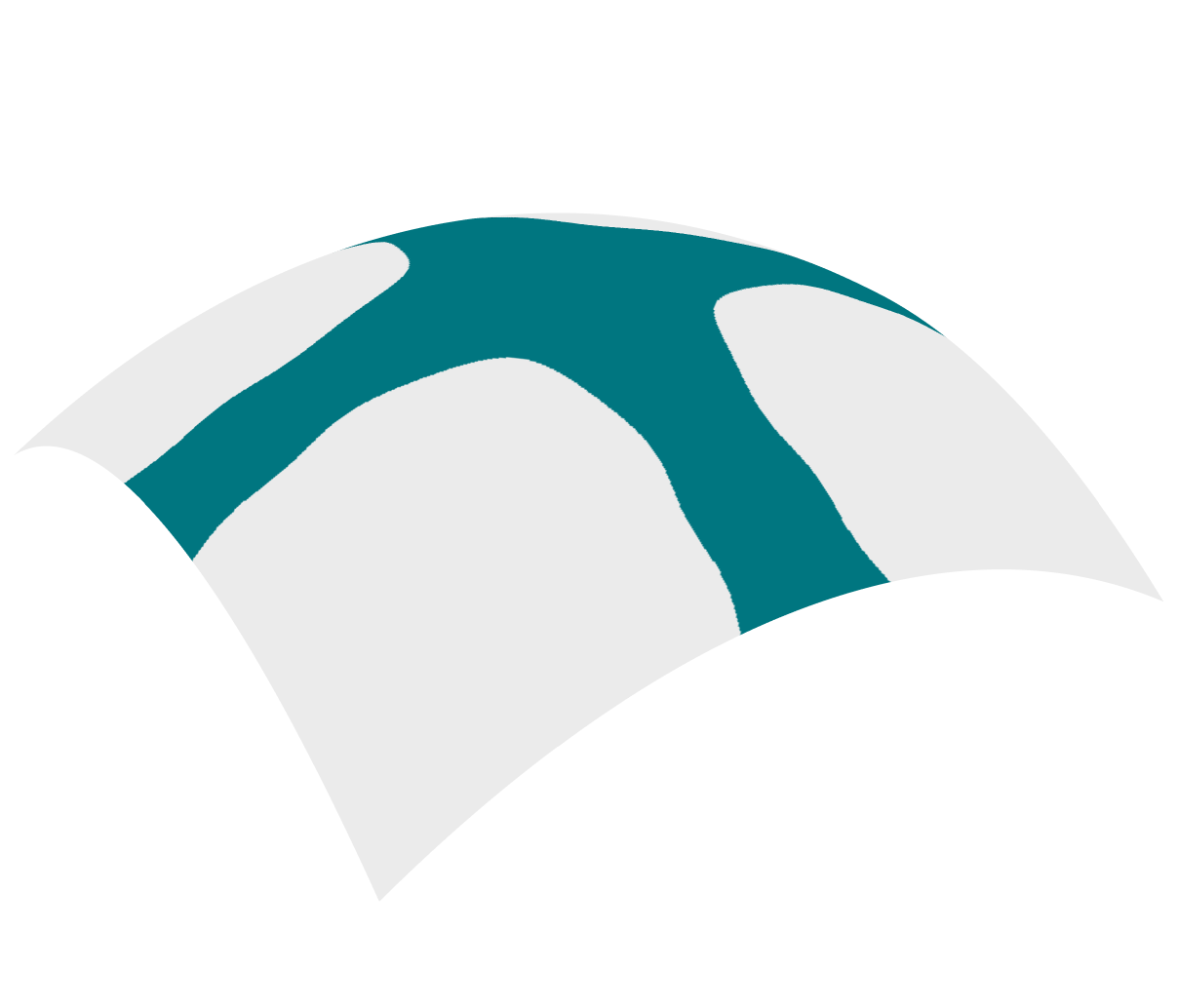}\\
				\vspace{0.1cm}
				\includegraphics[scale = 0.1]{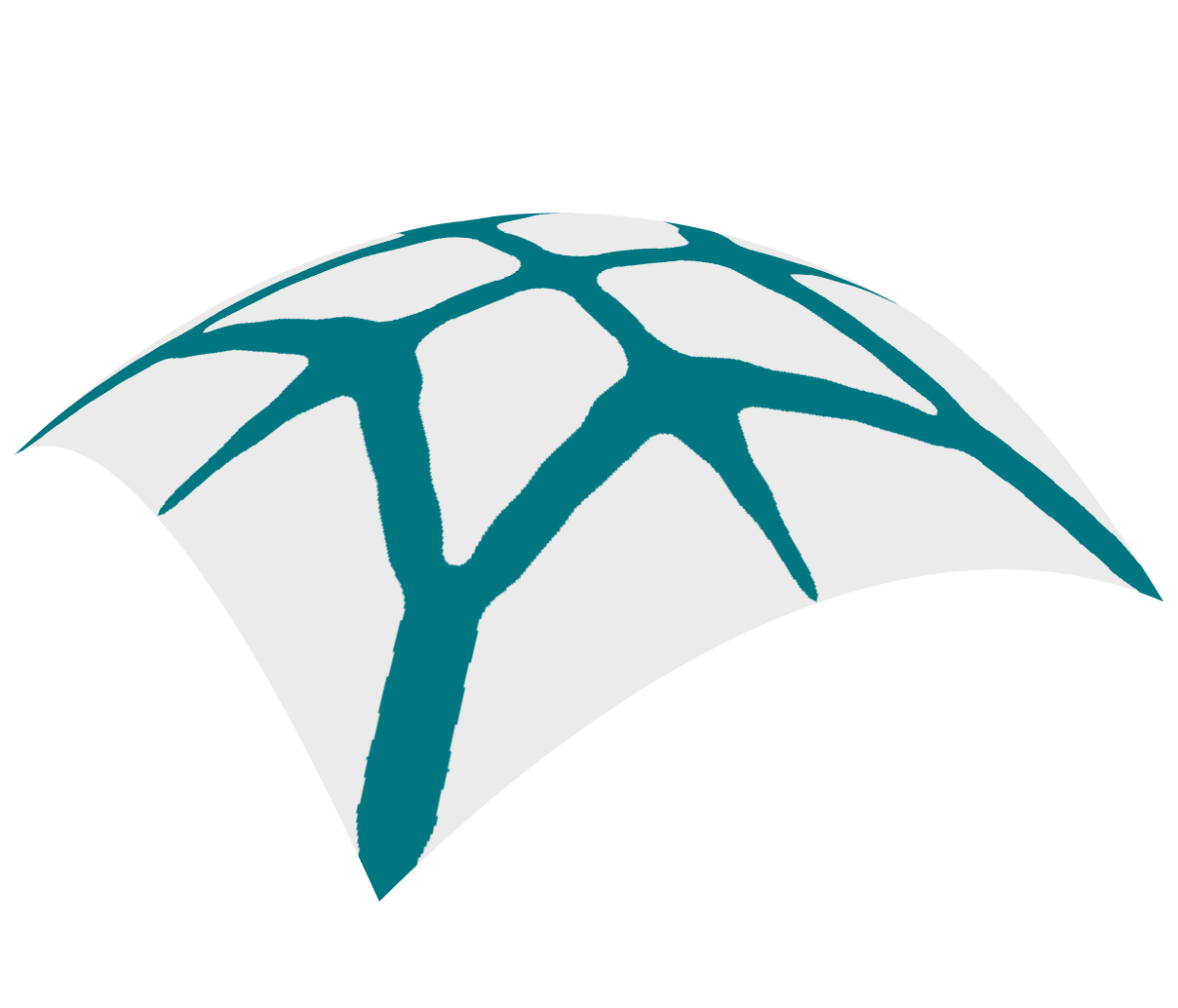}\\
				\vspace{0.2cm}
			\end{minipage}
		}
		\subfigure[\label{fig:IGA_100 top} top view]{
			\begin{minipage}[t]{0.125\linewidth}
				\centering
				\includegraphics[scale = 0.1]{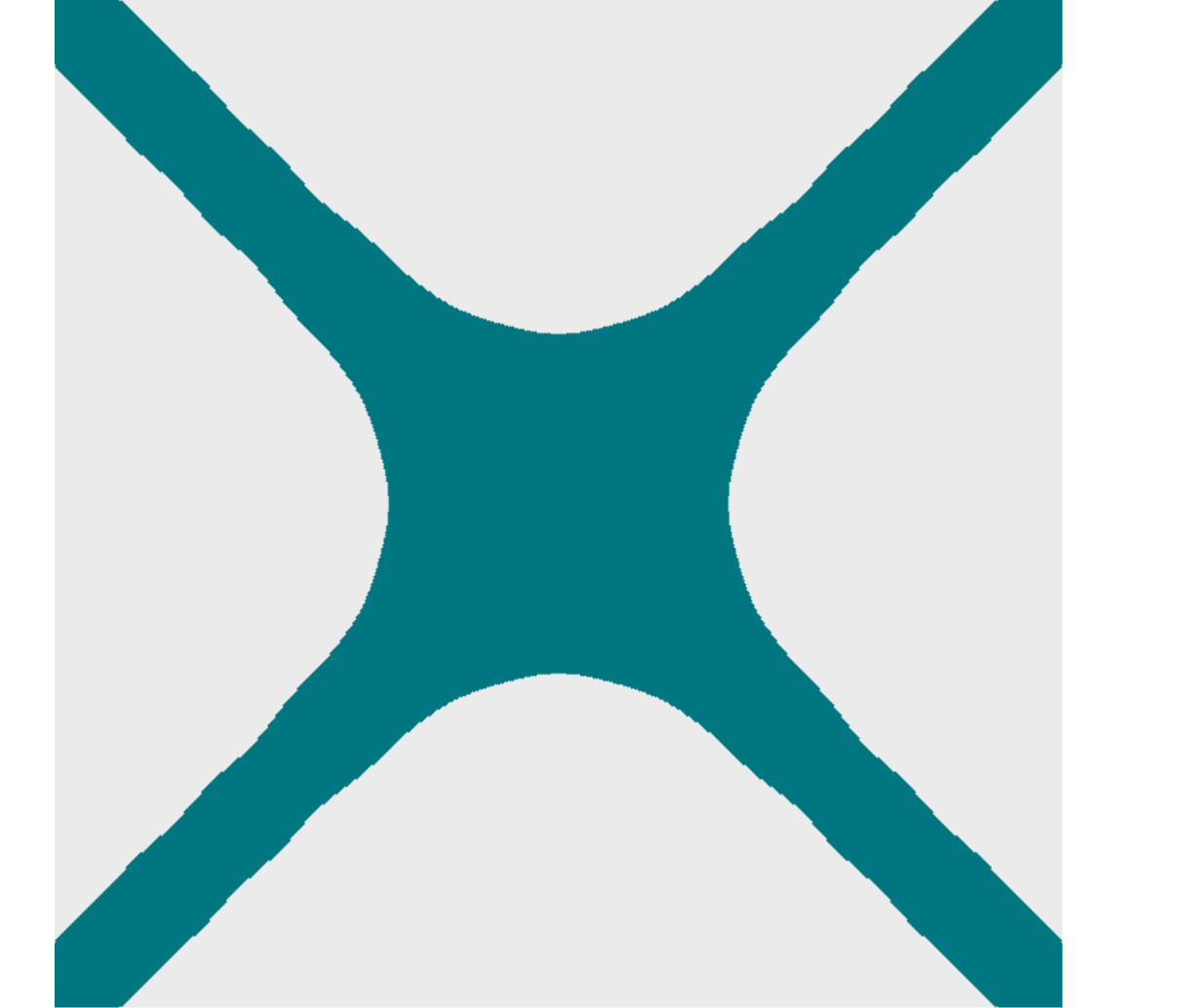}\\
				\vspace{0.1cm}
				\includegraphics[scale = 0.1]{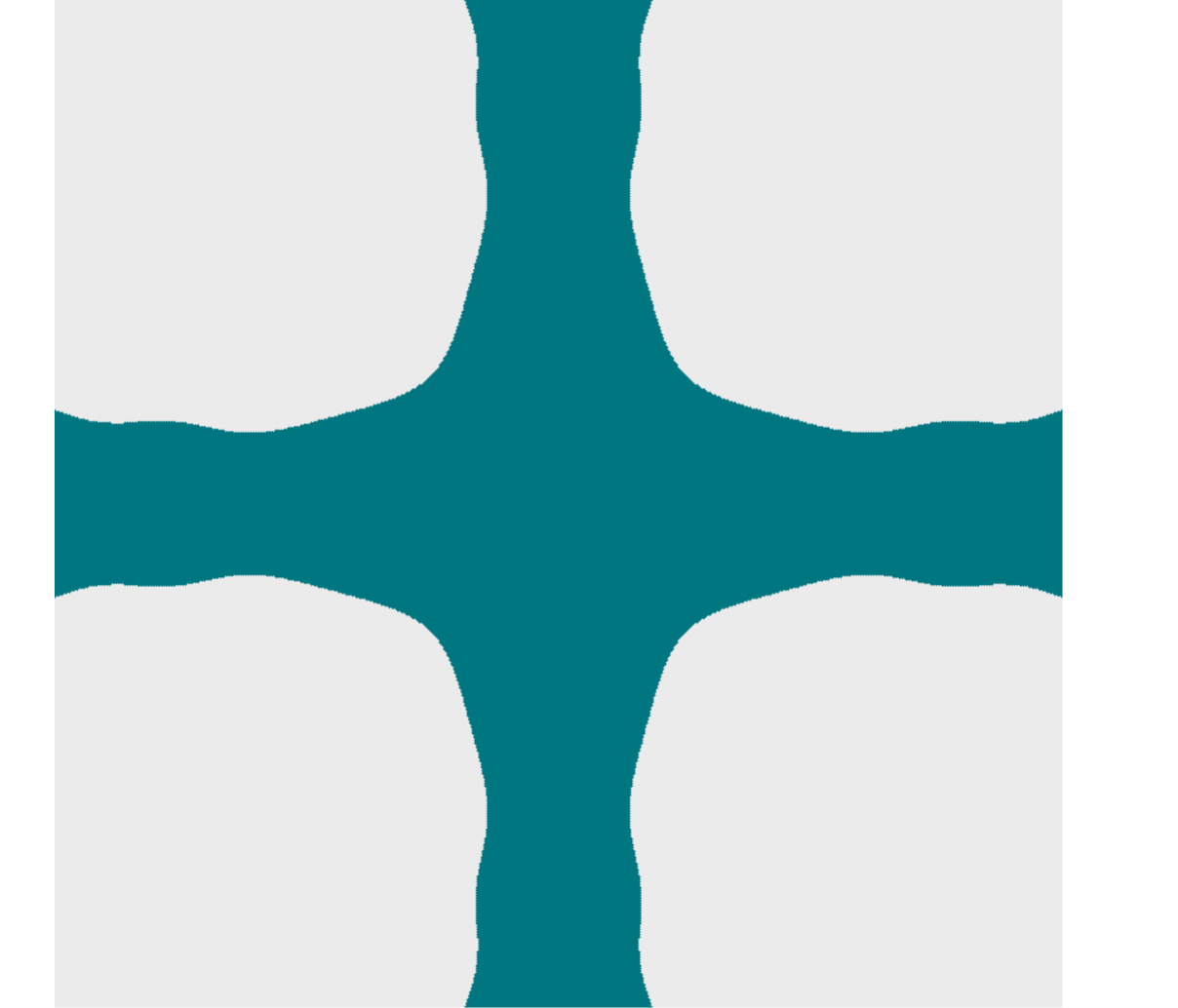}\\
				\vspace{0.1cm}
				\includegraphics[scale = 0.1]{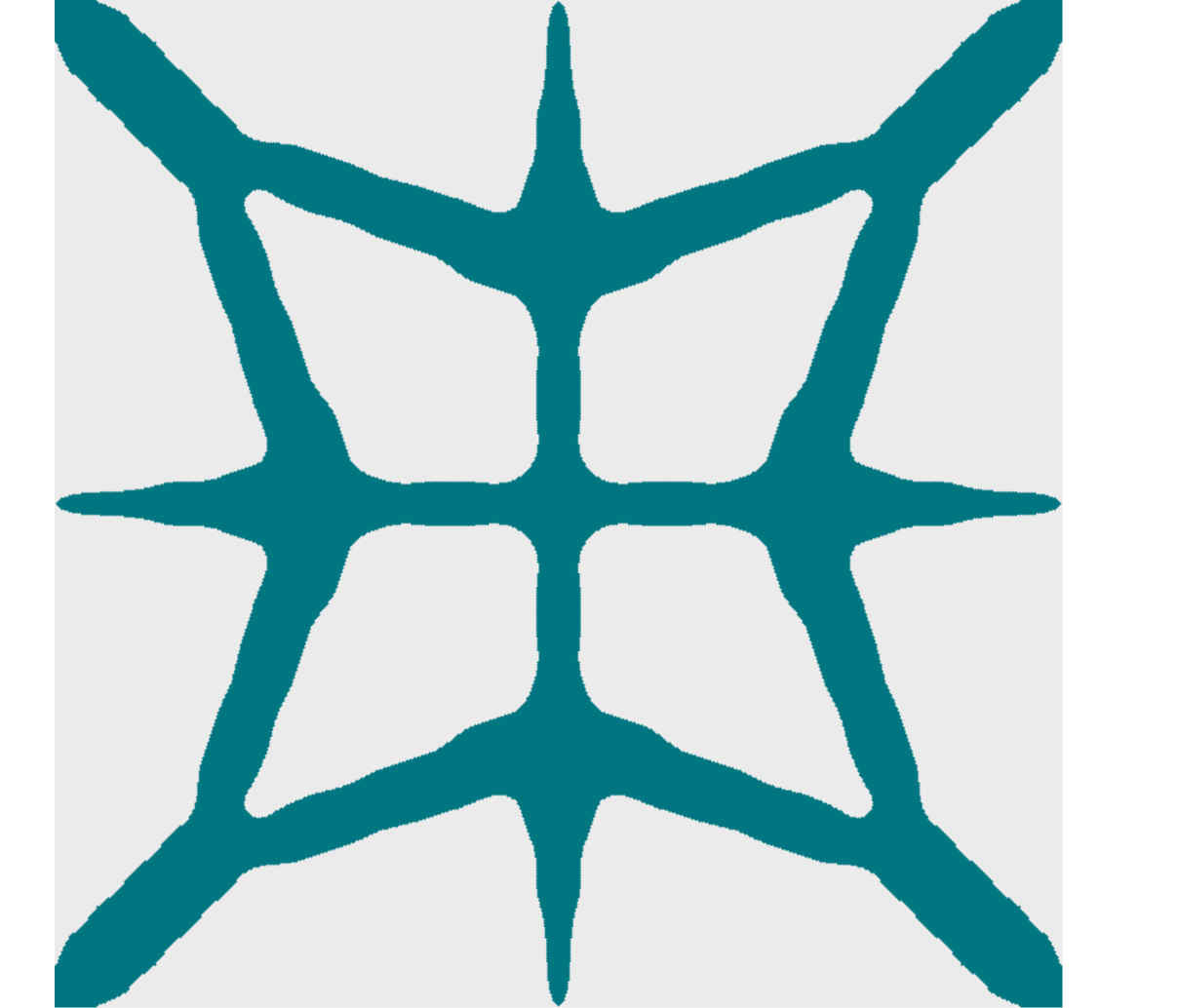}\\
				\vspace{0.2cm}
			\end{minipage}
		}
		\subfigure[\label{fig:FEA_100 front} FEA100]{
			\begin{minipage}[t]{0.125\linewidth}
				\centering
				\includegraphics[scale = 0.1]{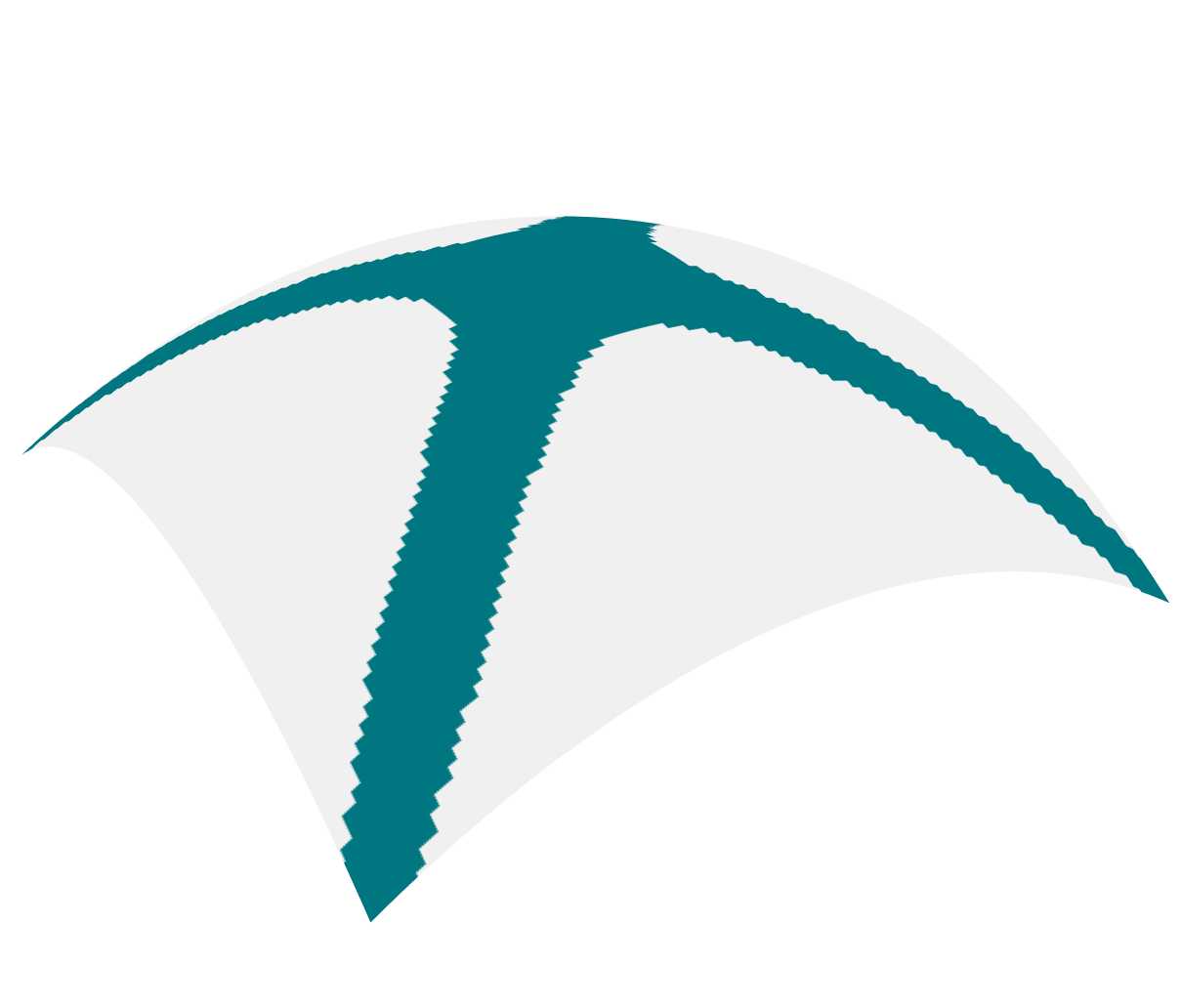}\\
				\vspace{0.1cm}
				\includegraphics[scale = 0.1]{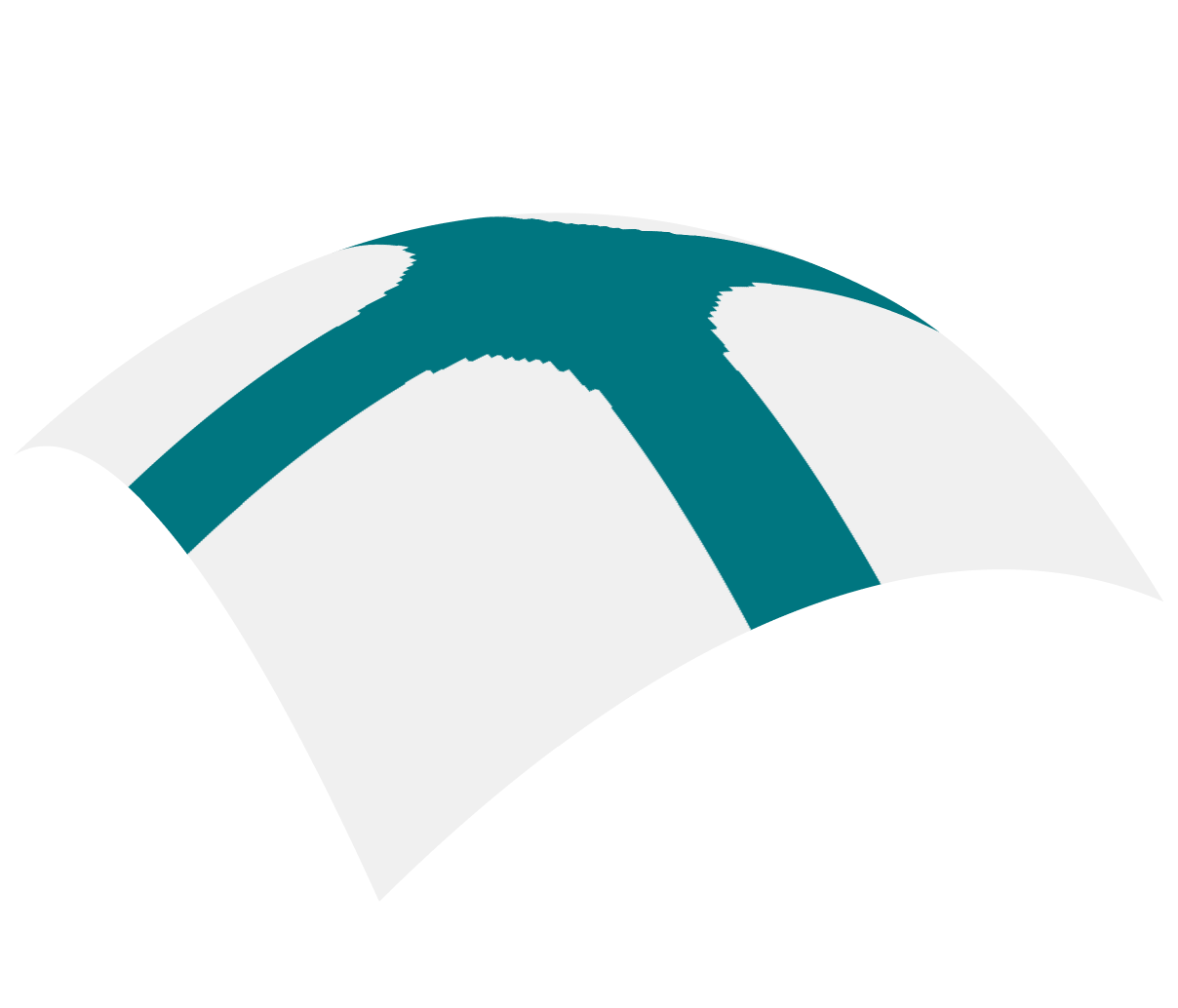}\\
				\vspace{0.1cm}
				\includegraphics[scale = 0.1]{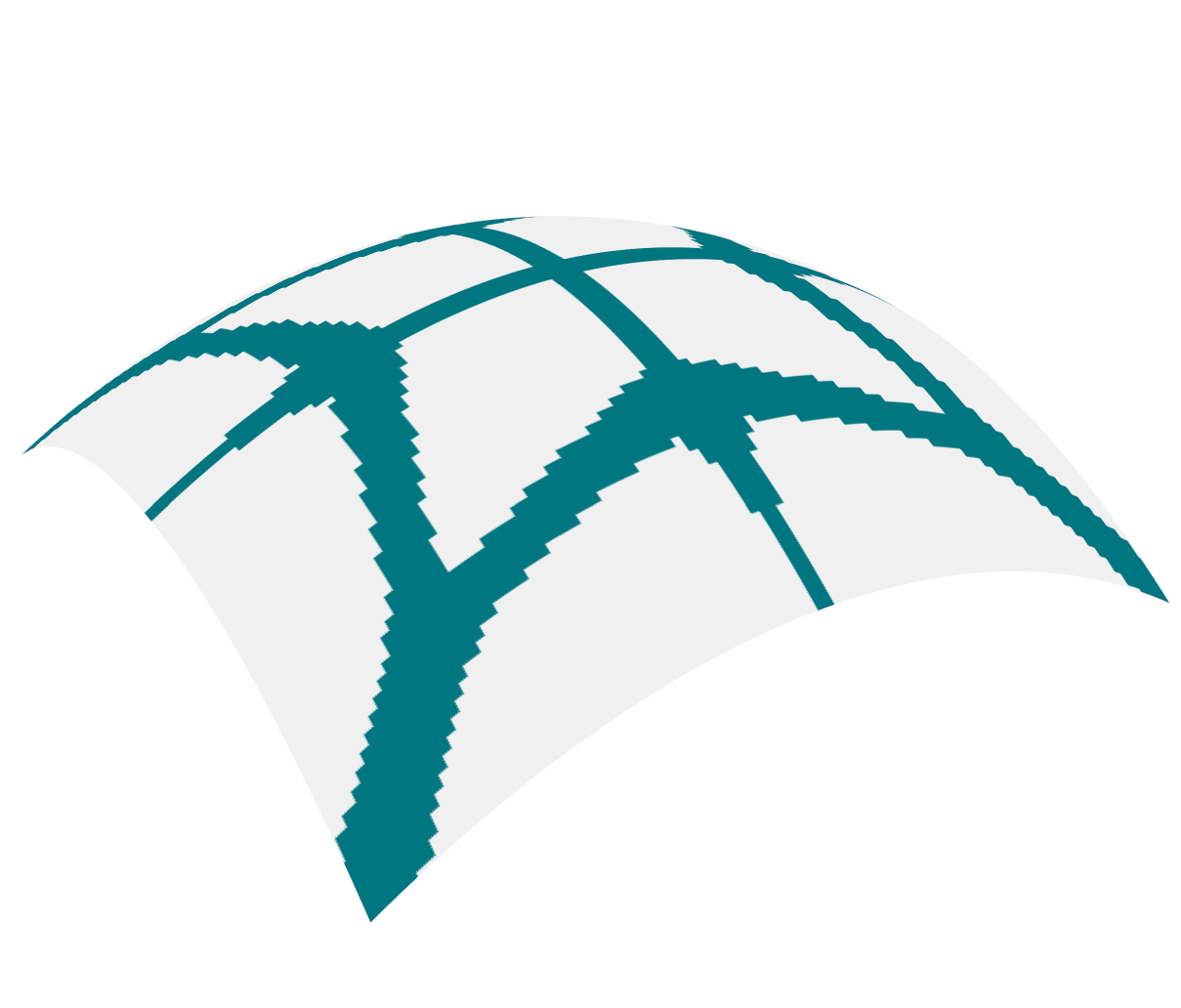}\\
				\vspace{0.2cm}
			\end{minipage}
		}
		\subfigure[\label{fig:FEA_100 top} top view ]{
			\begin{minipage}[t]{0.125\linewidth}
				\centering
				\includegraphics[scale = 0.1]{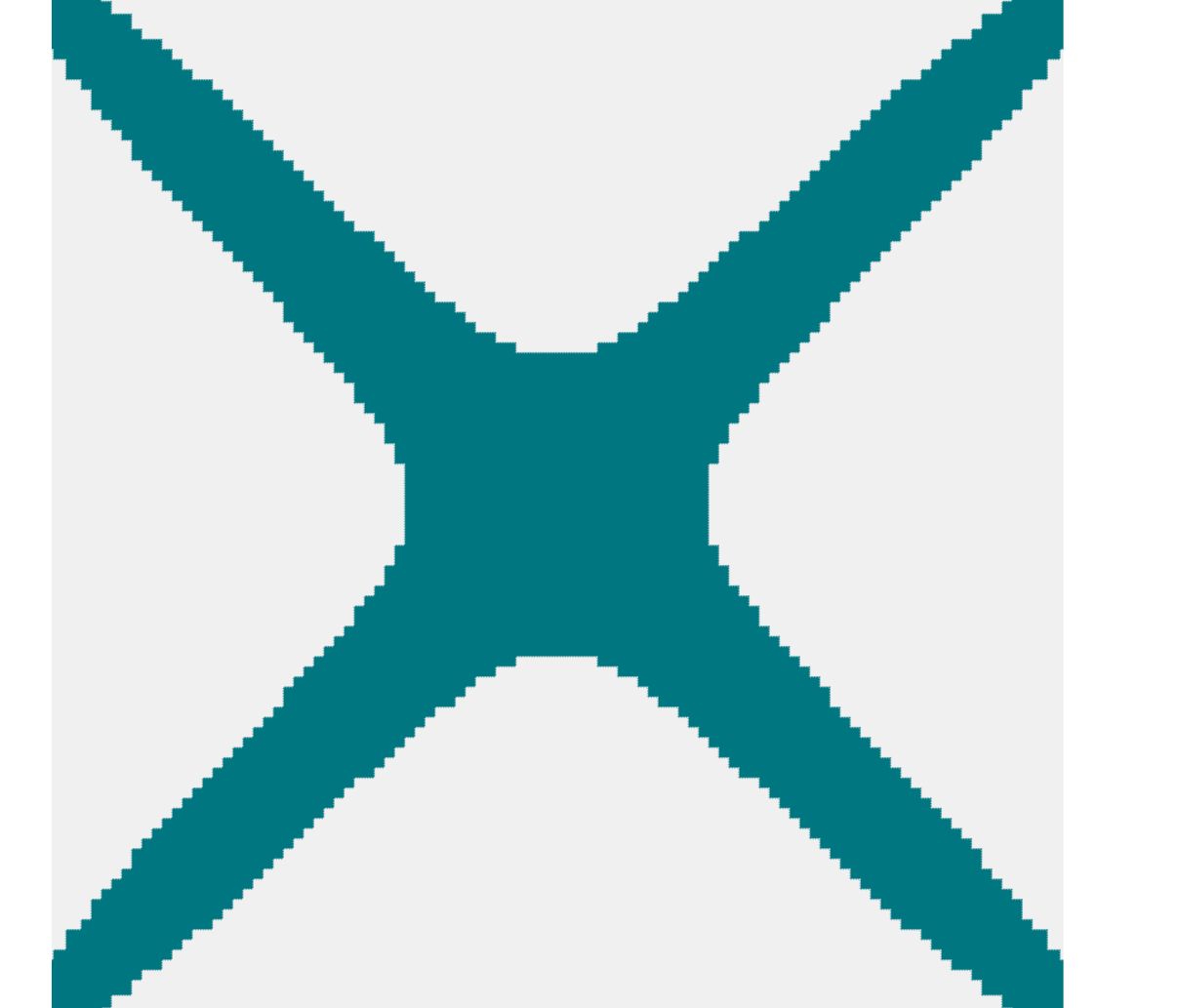}\\
				\vspace{0.1cm}
				\includegraphics[scale = 0.1]{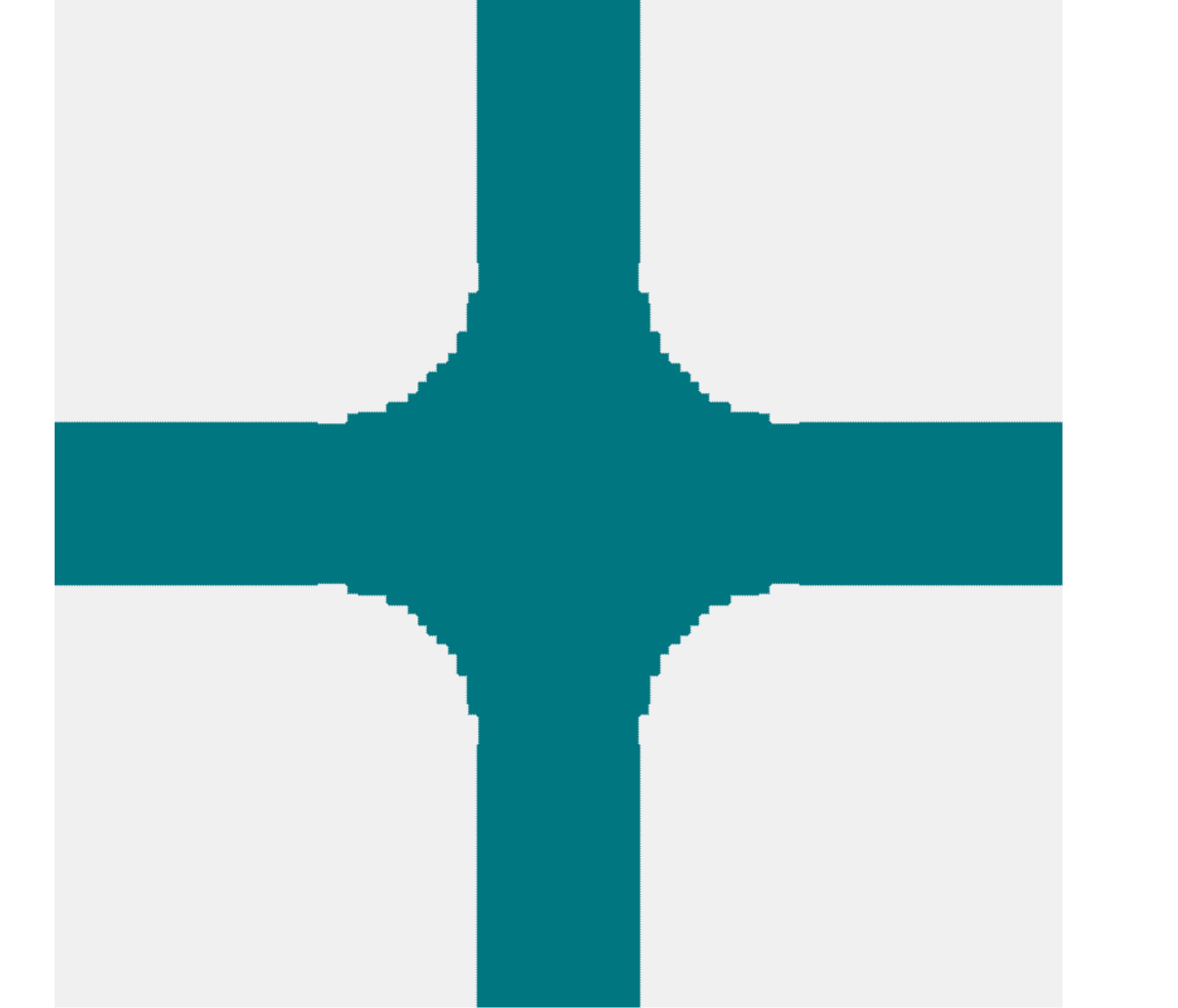}\\
				\vspace{0.1cm}
				\includegraphics[scale = 0.1]{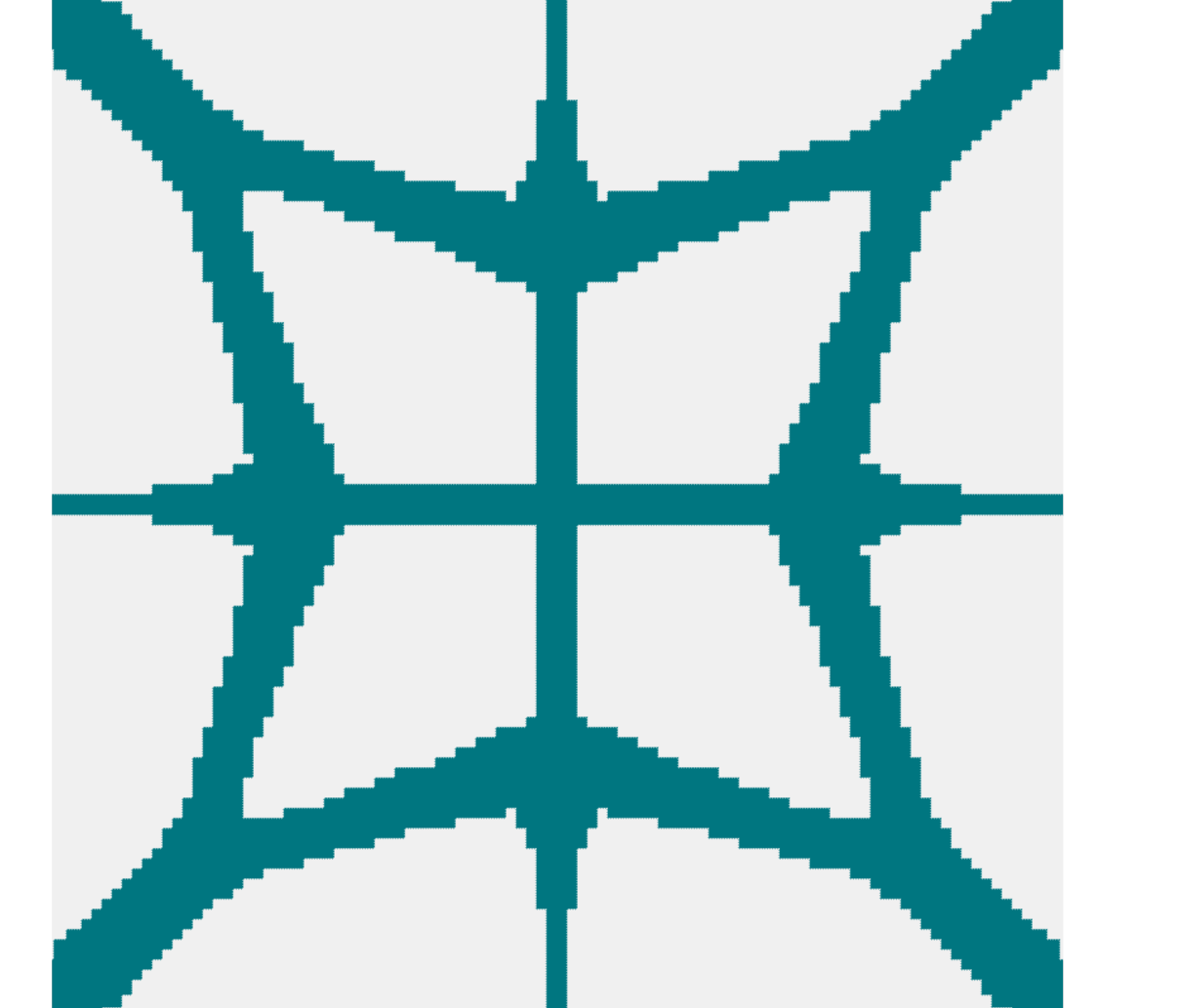}\\
				\vspace{0.2cm}
			\end{minipage}
		}
		\caption{Optimized structures comparisons between IGA-SIMP and FEA-SIMP methods.  The leftmost column describes models with different loads and boundary conditions, followed by the results using the IGA-SIMP method with $50\times 50$ elements and $100\times 100$ elements in the analysis model. The last two columns show the results using the FEA-SIMP method with $100\times 100$ elements. Top view results are shown after the front view results.}
		\label{fig: geometrical comparison - IGA-SIMP VS FEM-SIMP}
	\end{figure*}
	\begin{table*}[!ht]
		\setlength{\belowcaptionskip}{0.35cm}
		\centering
		\caption{Statistics on IGA-SIMP method vs. FEA-SIMP method under different parameters for 3 cases.}
		\begin{tabular*}{\hsize}{@{}@{\extracolsep{\fill}}|c|ccccccc|@{}}
			\hline
			model & method & $DOF$ & $N_E$ & $N_{DV}$ & C& $\lVert \rho-\rho_D \rVert_{L_1}$ & Time (s)\\ 
			\hline 
			\multirow{3}{*}{case 1} & \multirow{2}{*}{IGA-SIMP} & $52\times 52$ & $50\times 50$  & $15\times 15$ & 9.49 & 0.024 & 129 \\
			~ & ~ & $102\times 102$ & $100\times 100$ & $15\times 15$ & 9.55 & 0.012 & 705 \\
			~ & FEA-SIMP & $201\times 201$ & $100\times 100$ & $201\times 201$ & 9.32 & 0.028 & 1337 \\
			\hline
			\multirow{3}{*}{case 2} & \multirow{2}{*}{IGA-SIMP} & $52\times 52$ & $50\times 50$  & $15\times 15$ & 3.27 & 0.015 & 131 \\
			~ & ~ & $102\times 102$ & $100\times 100$ & $15\times 15$ & 3.31 & 0.010 & 673 \\
			~ & FEA-SIMP & $201\times 201$ & $100\times 100$ & $201\times 201$ & 3.36 & 0.017 & 1246 \\
			\hline    	
			\multirow{3}{*}{case 3} & \multirow{2}{*}{IGA-SIMP} & $52\times 52$ & $50\times 50$  & $30\times 30$ & 237.86 & 0.044 & 134 \\
			~ & ~ & $102\times 102$ & $100\times 100$ & $30\times 30$ & 239.40 & 0.020 & 710 \\
			~ & FEA-SIMP & $201\times 201$ & $100\times 100$ & $201\times 201$ & 232.16 & 0.038 & 1343 \\
			\hline
		\end{tabular*}
		\label{Table: IGA-SIMP VS FEM-SIMP}
	\end{table*}
	\begin{figure*}[htpb]
		\subfigure[case 1]
		{
			\begin{overpic}[width= 0.32\textwidth]{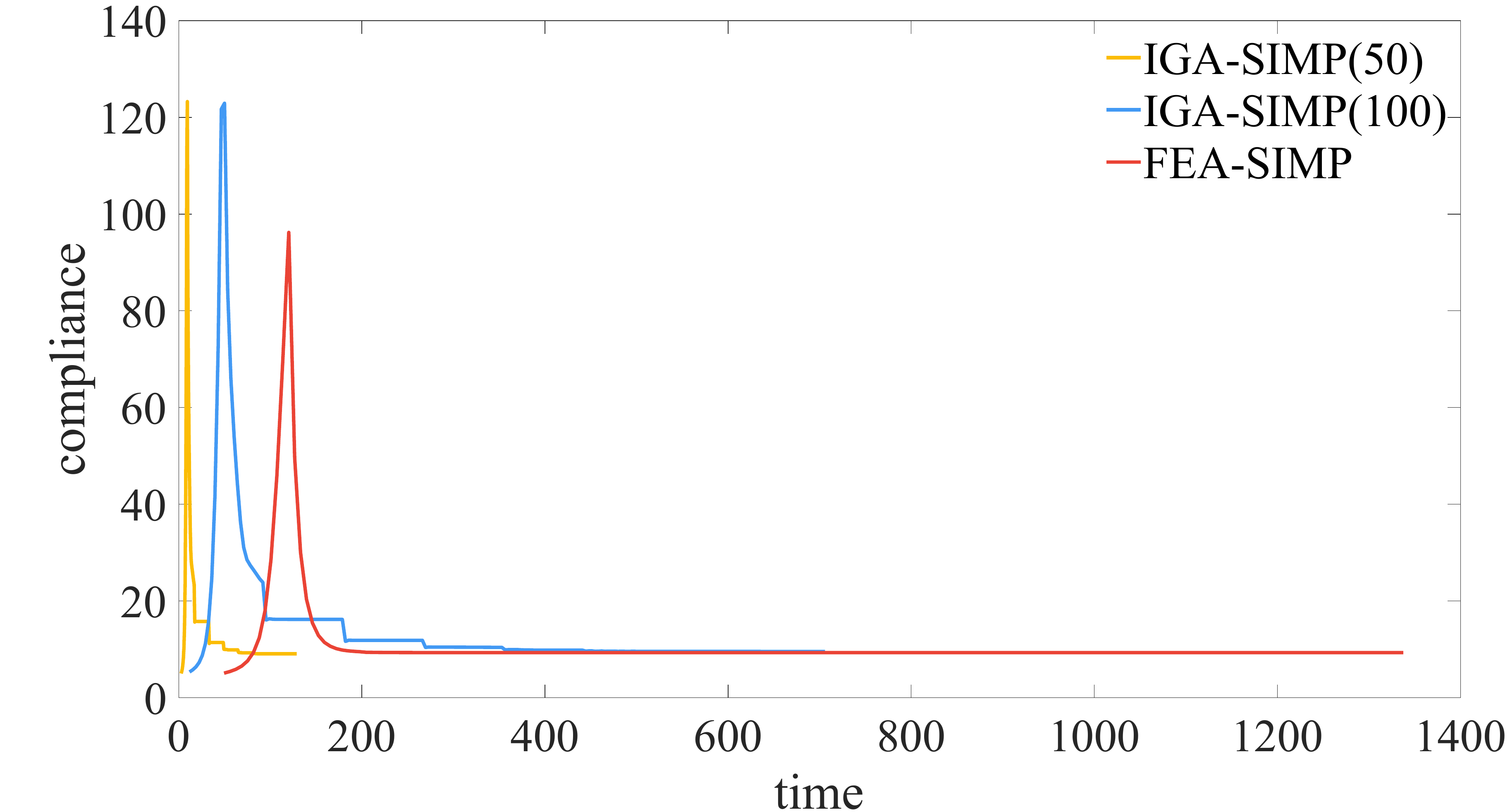}
			\end{overpic}
		}
		\subfigure[ case 2]
		{
			\begin{overpic}[width= 0.32\textwidth]{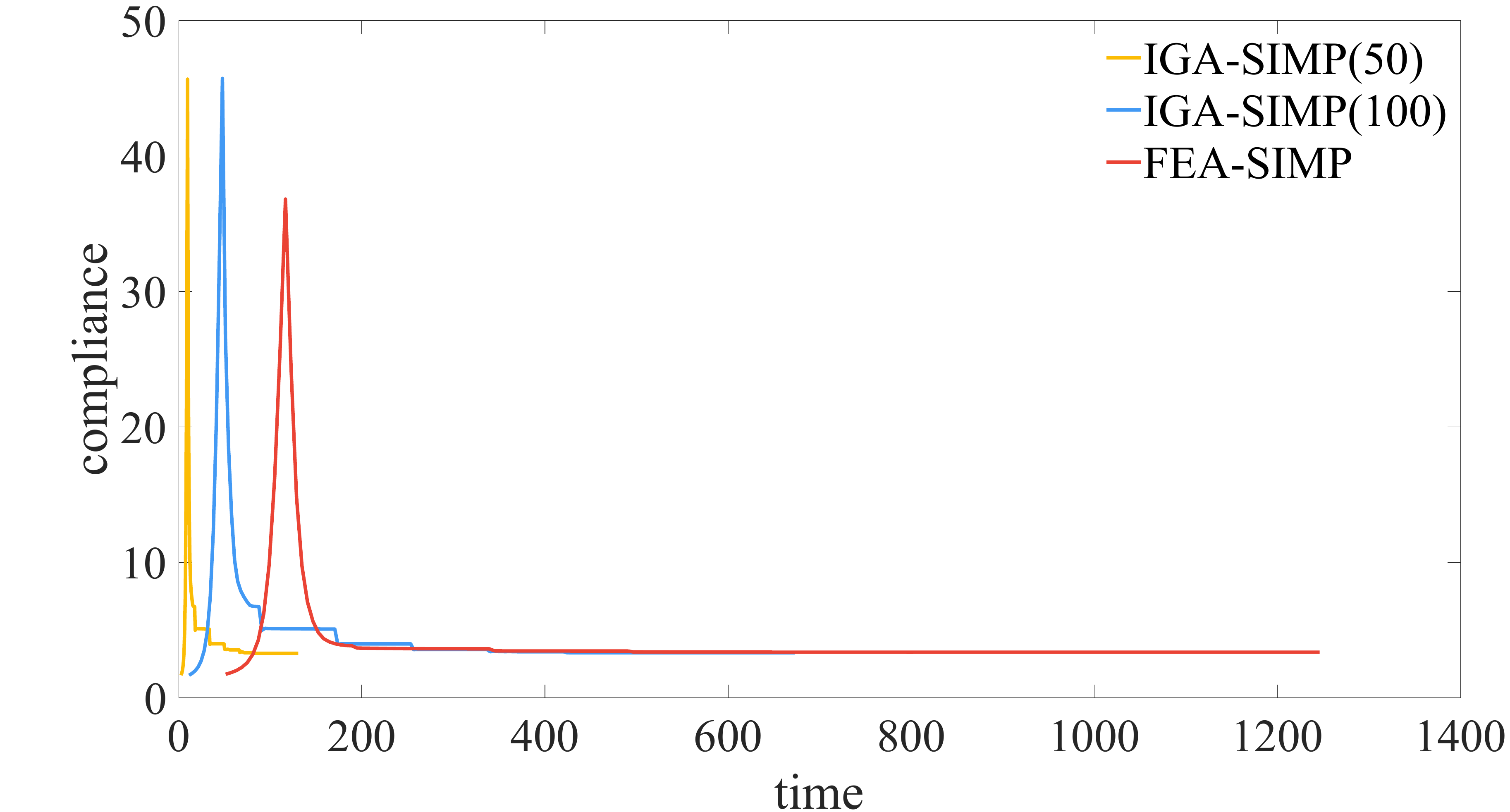}
			\end{overpic}
		}
		\subfigure[ case 3]
		{
			\begin{overpic}[width= 0.32\textwidth]{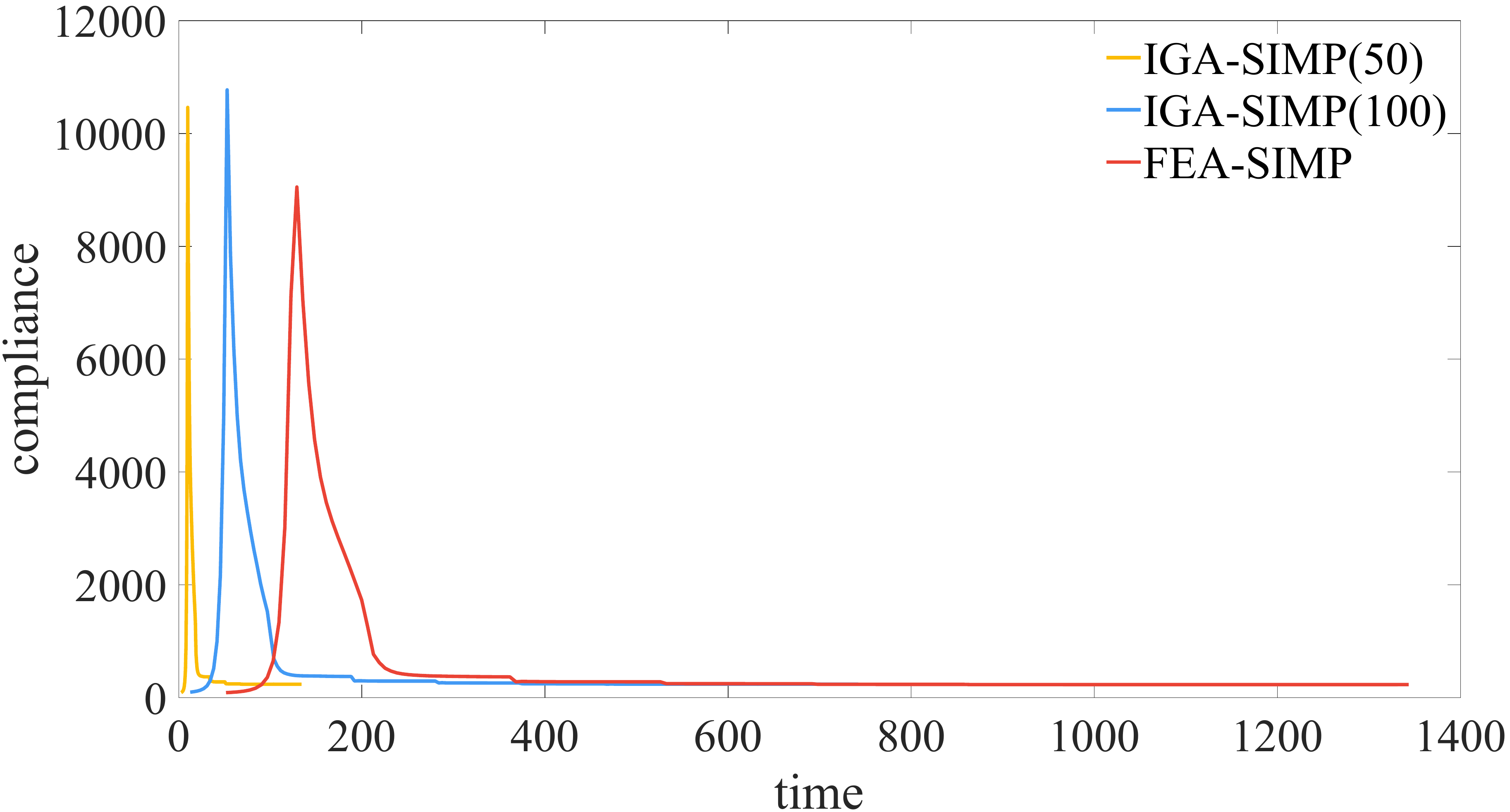}
			\end{overpic}
		}
		\caption{Iteration history for IGA-SIMP method and FEA-SIMP method.}
		\label{fig: iteration history.}
	\end{figure*}
	\begin{figure}
		\begin{center}
			\subfigure[case 1]{
				\includegraphics[width=0.31\textwidth]{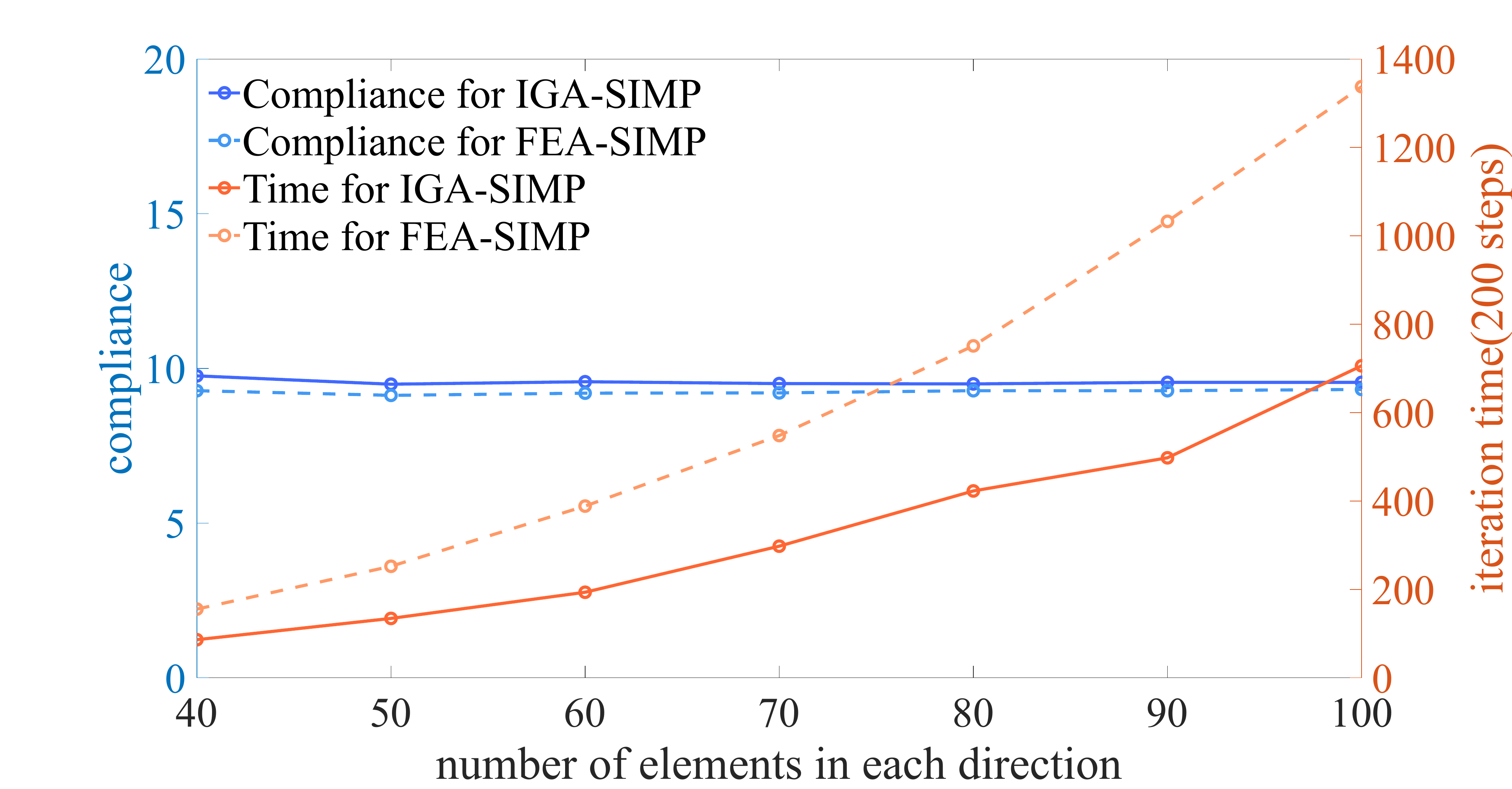}
				\subfigcapskip=-5pt	
			} 
			\subfigure[case 2]{
				\includegraphics[width=0.31\textwidth]{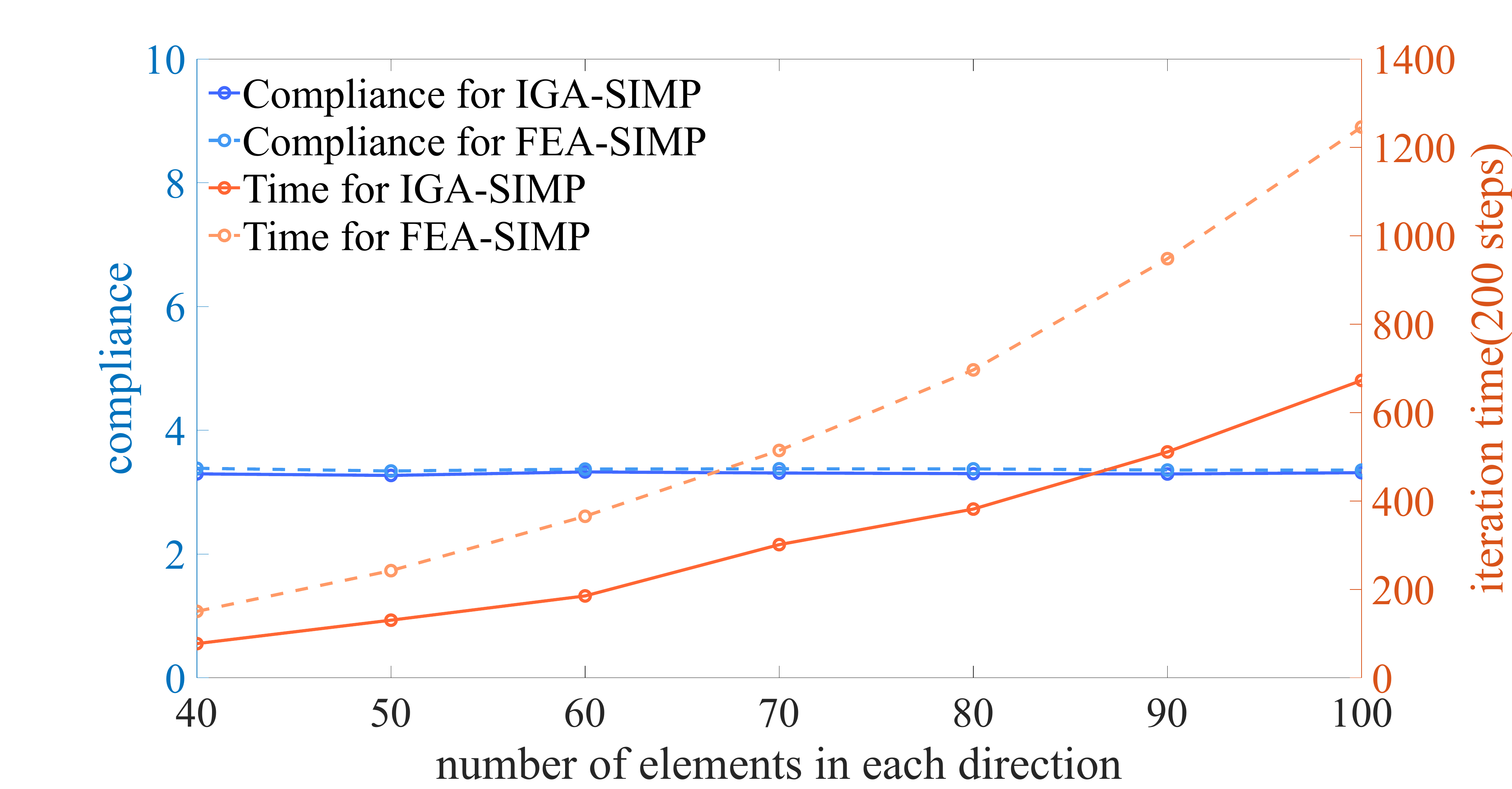}
				\subfigcapskip=-5pt	
			}
			\subfigure[case 3]{
				\includegraphics[width=0.31\textwidth]{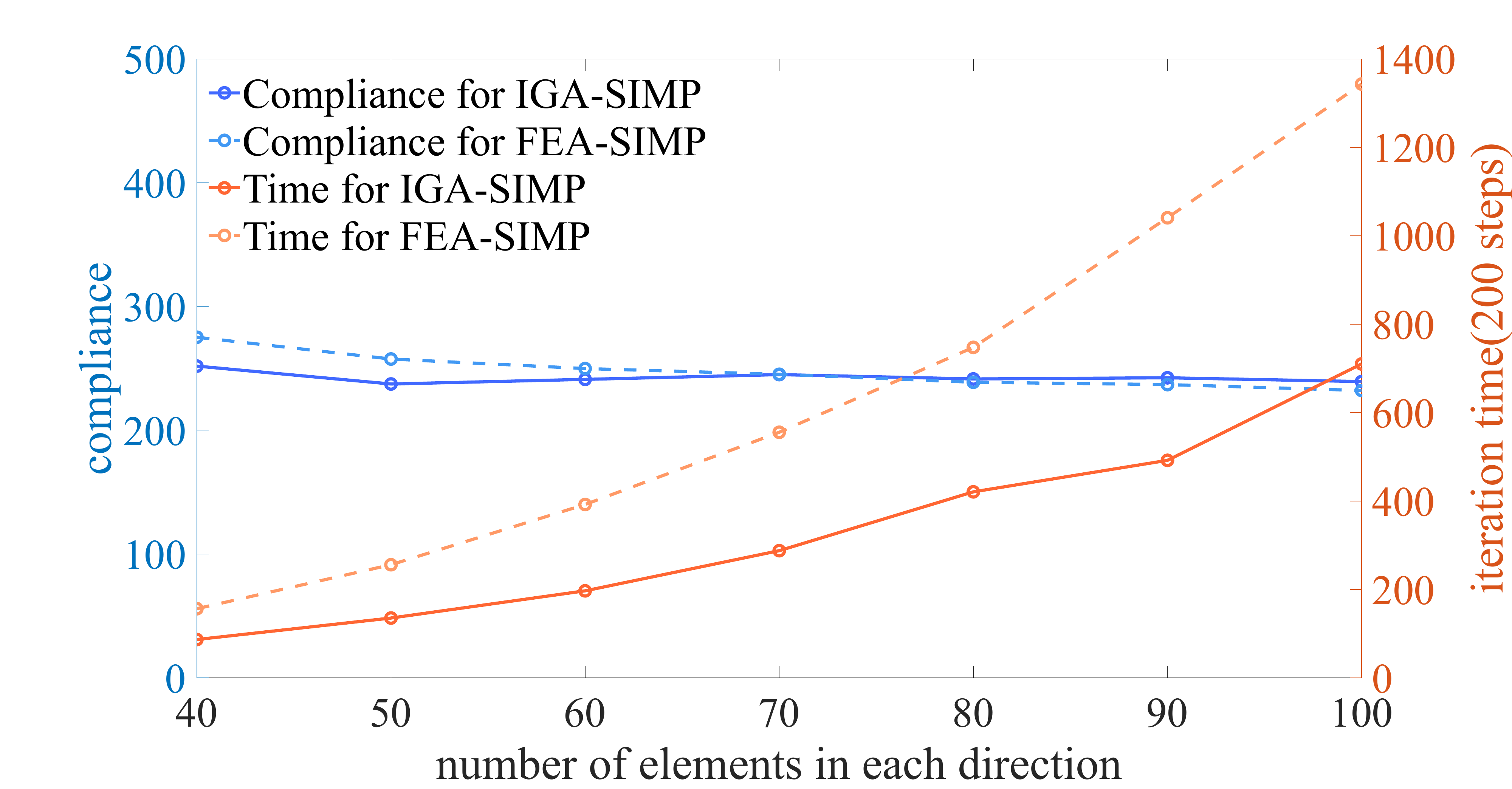}
				\subfigcapskip=-5pt	
			}
		\end{center}
		\vspace{-0.4cm}
		\caption{Numerical comparisons between IGA-SIMP method and FEA-SIMP method for shell topology optimization on three cases (a)-(c). The number of elements used in the analysis model varies from 40 to 100.}
		\label{fig: numerical comparison - IGA-SIMP VS FEM-SIMP}
	\end{figure}

	The IGA-SIMP method and FEA-SIMP method are compared on a hyperbolic shell subjected to different loads and boundary conditions. To minimize errors arising from the accuracy of geometric representation, we employ bi-quadratic Lagrange isoparametric elements to depict the geometric structure within the FEA-SIMP method. Three distinct scenarios, as illustrated in Fig.\ref{fig: geometrical comparison - IGA-SIMP VS FEM-SIMP}(a), are examined.
	The blue triangle denotes a fixed position, while multiple triangles on a side indicate that the entire side is fixed. The red arrow represents a force, with a single arrow representing the force applied to the top center and multiple arrows indicating forces along the isoparametric curves (e.g., s=0.5, t=0.5). These markings apply to the force and boundary conditions in all subsequent examples.
	The global volume constraint presented in Eq.~\eqref{eq: formulation 1} is taken into account, and the available volume of the solid material is denoted by $V^* = 0.3V_s$, where $V_s$ represents the volume of full solid shell structures.

	The results of the IGA-SIMP method and FEA-SIMP method are presented in Fig.~\ref{fig: geometrical comparison - IGA-SIMP VS FEM-SIMP}, while the iteration history is depicted in Fig.~\ref{fig: iteration history.} Both methods yield almost identical optimum models in terms of geometry. However, the shell structures obtained through the FEA-SIMP method exhibit rougher boundaries compared to those obtained through the IGA-SIMP method.

	Table~\ref{Table: IGA-SIMP VS FEM-SIMP} lists relevant parameters and numerical results. $DOF = m_2\times{n_2}$ and $ N_{DV} = m_1\times{n_1}$ represent the number of degrees of freedom and design variables, respectively.
	To evaluate the differences in the material distribution between the IGA-based method with $50\times 50$ and $100\times 100$ elements and the FEA-based method with $100\times 100$ elements, we define their respective density functions in the parameter space $[0,1]^2$  as $\rho_A(s,t)$, $\rho_B(s,t)$ and $\rho_C(s,t)$. 
	Through further experiments, we have observed the convergence of the IGA-SIMP method with an increasing number of elements. Therefore we consider the result obtained from the IGA-SIMP method with $200\times 200$ elements as the ground truth with density function $\rho_D(s,t)$.
	We compute the $L_1$ norms of the differences 
	between $\rho_A$, $\rho_B$, $\rho_C$ and $\rho_D$ as shown in Table~\ref{Table: IGA-SIMP VS FEM-SIMP}. 
	The results indicate that both methods exhibit nearly identical accuracy, with errors within $5\%$.
	Furthermore, the compliance of optimized results using the IGA-SIMP method with $50 \times 50$ elements is close to those obtained using $100 \times 100$ elements, but the former optimization method is approximately 5 times faster than the latter. Additionally, compared to the IGA-SIMP method, the FEA-SIMP method under the $100 \times 100$ element framework requires twice as much time. 
	As a result, the IGA-SIMP method is about 10 times faster than the FEA-SIMP method while producing nearly identical optimum shell structures with similar geometry, performance and accuracy, demonstrating its advantage over FEA-SIMP method.

	Moreover, there is a disparity between the optimal objective function values when varying the number of elements($N_E$) analyzed. To further investigate, tests were conducted as shown in Fig. \ref{fig: numerical comparison - IGA-SIMP VS FEM-SIMP} in which the number of elements $N_E$ for both the s- and t-directions ranged from 40 to 100 in increments of 10.
	From the figure, we can deduce that when using a coarse mesh, the optimal structures exhibit lower stability and higher compliance compared to when a finer mesh is used. The compliance of the optimal shell structures gradually converges as the number of elements $N_E$ increases for both methods.

	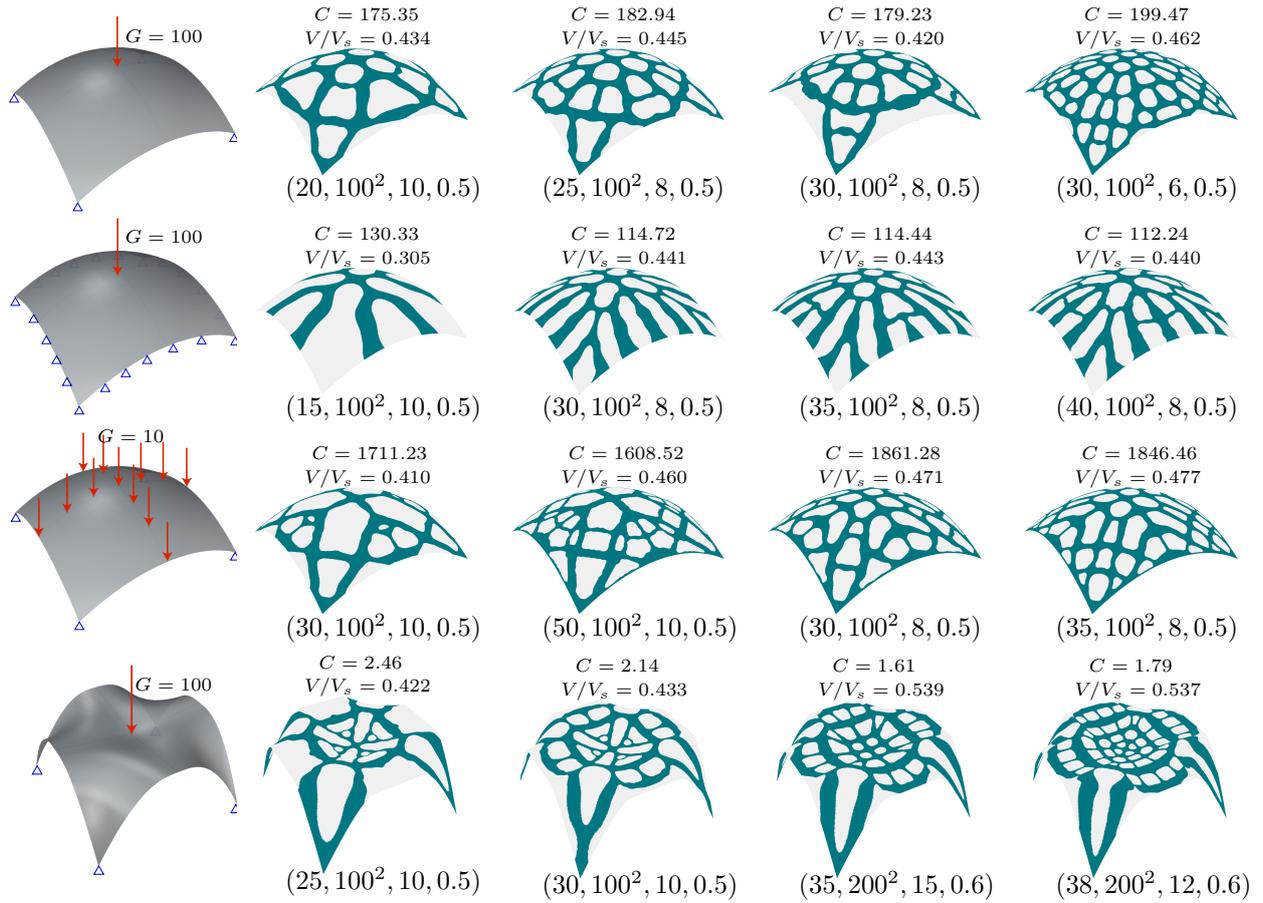
\begin{figure*}[htbp]
		\begin{center}
			\begin{minipage}[t]{0.13\linewidth}
				\centering
				\begin{overpic}[scale = 0.12]{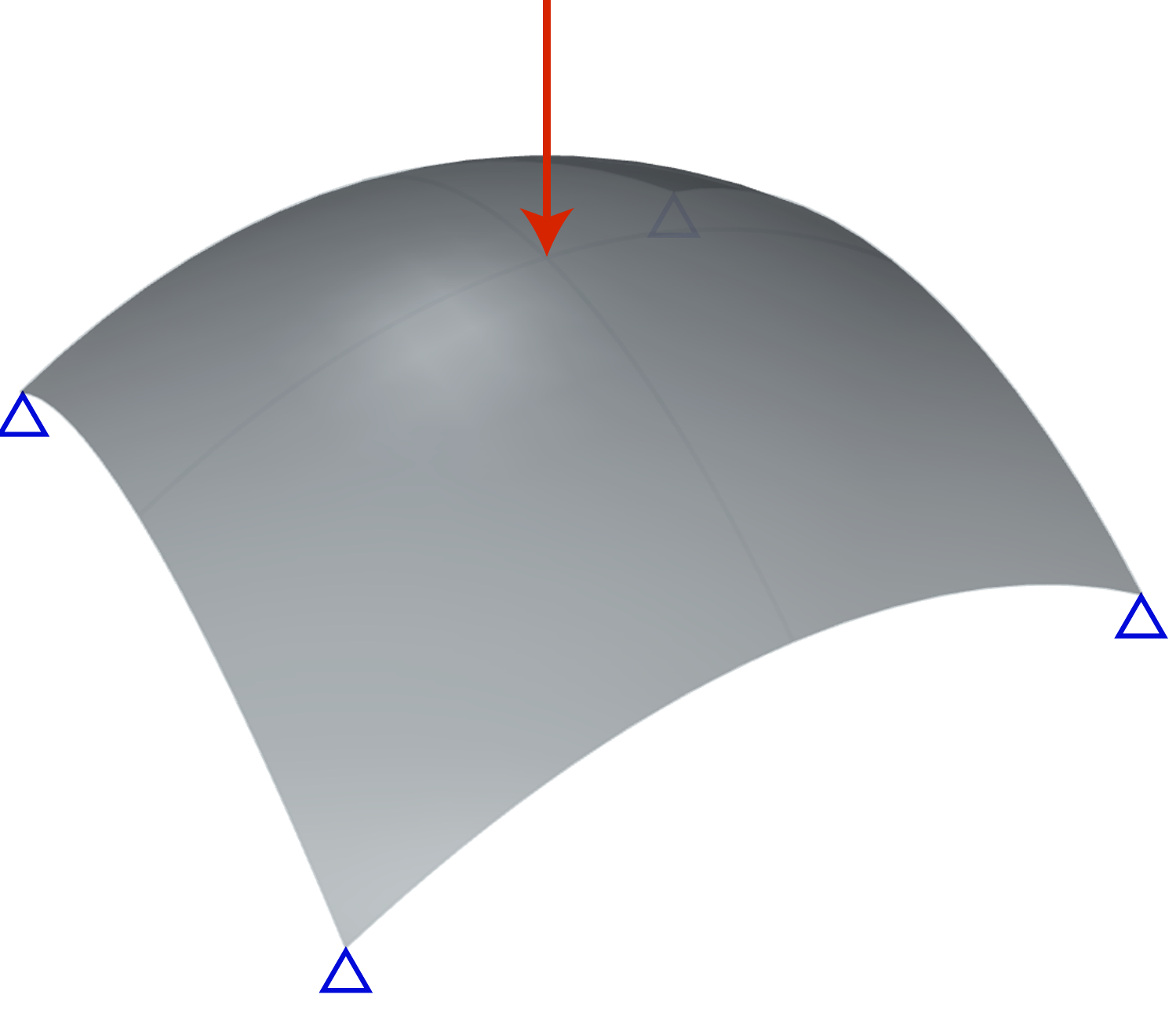}
					\put(50,75){\scriptsize{$G = 100$}}
				\end{overpic}
				\\
				\begin{overpic}[scale = 0.12]{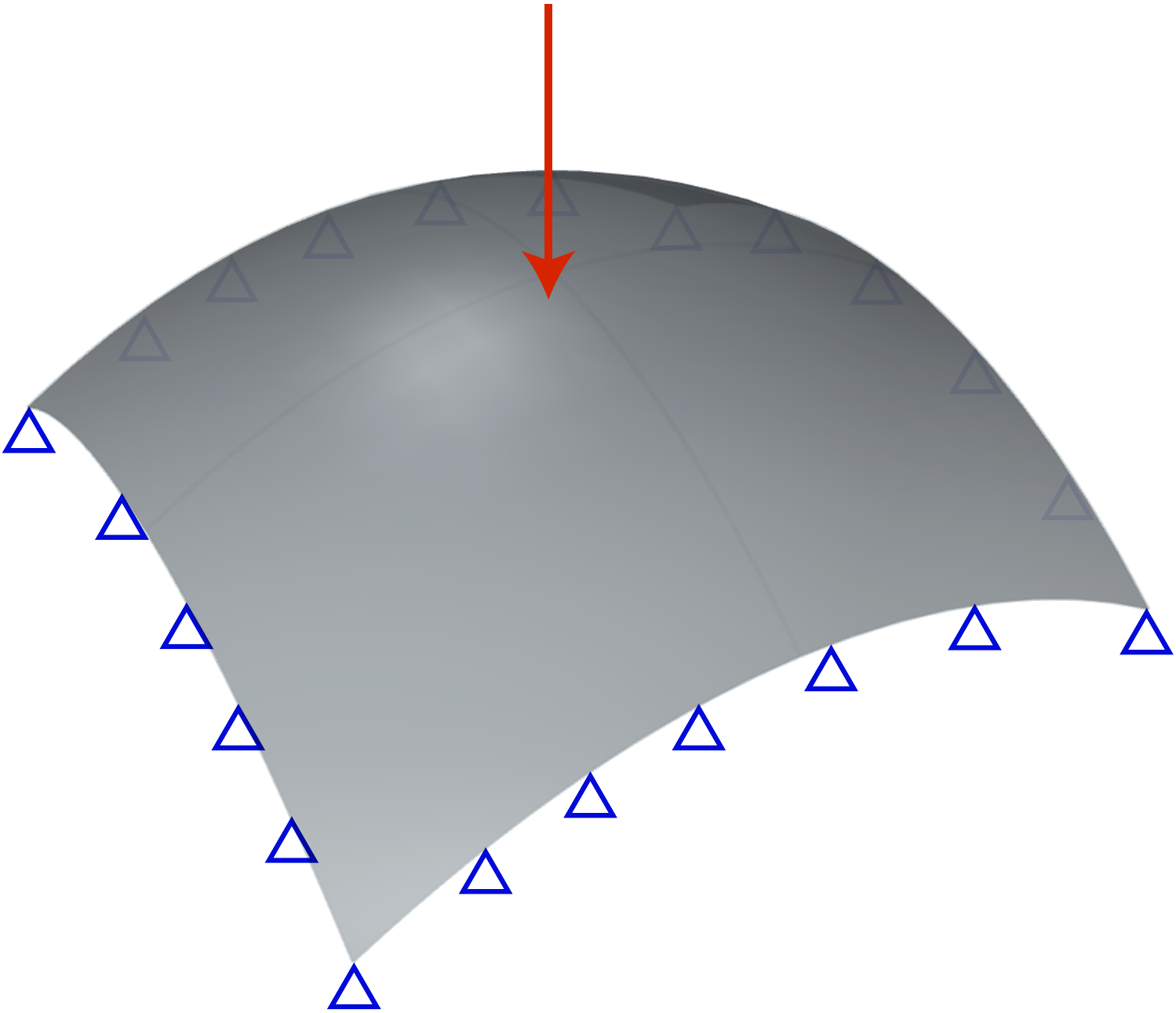}
					\put(50,75){\scriptsize{$G = 100$}}
				\end{overpic}
				\\
				\vspace{0.2cm}
				\begin{overpic}[scale = 0.12]{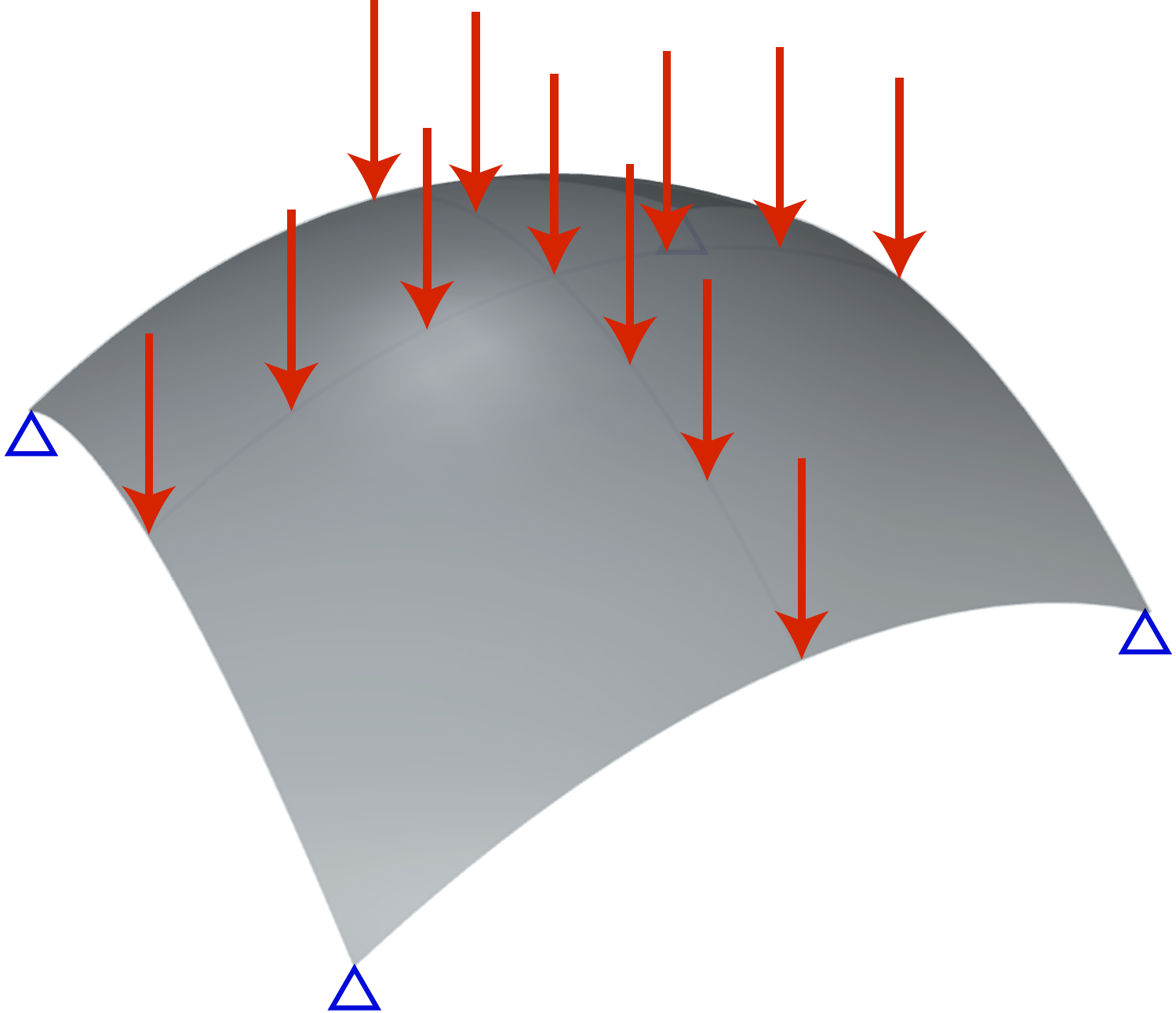}
					\put(38,82){\scriptsize{$G = 10$}}
				\end{overpic}
				\\
				\vspace{0.4cm}
				\begin{overpic}[scale = 0.13]{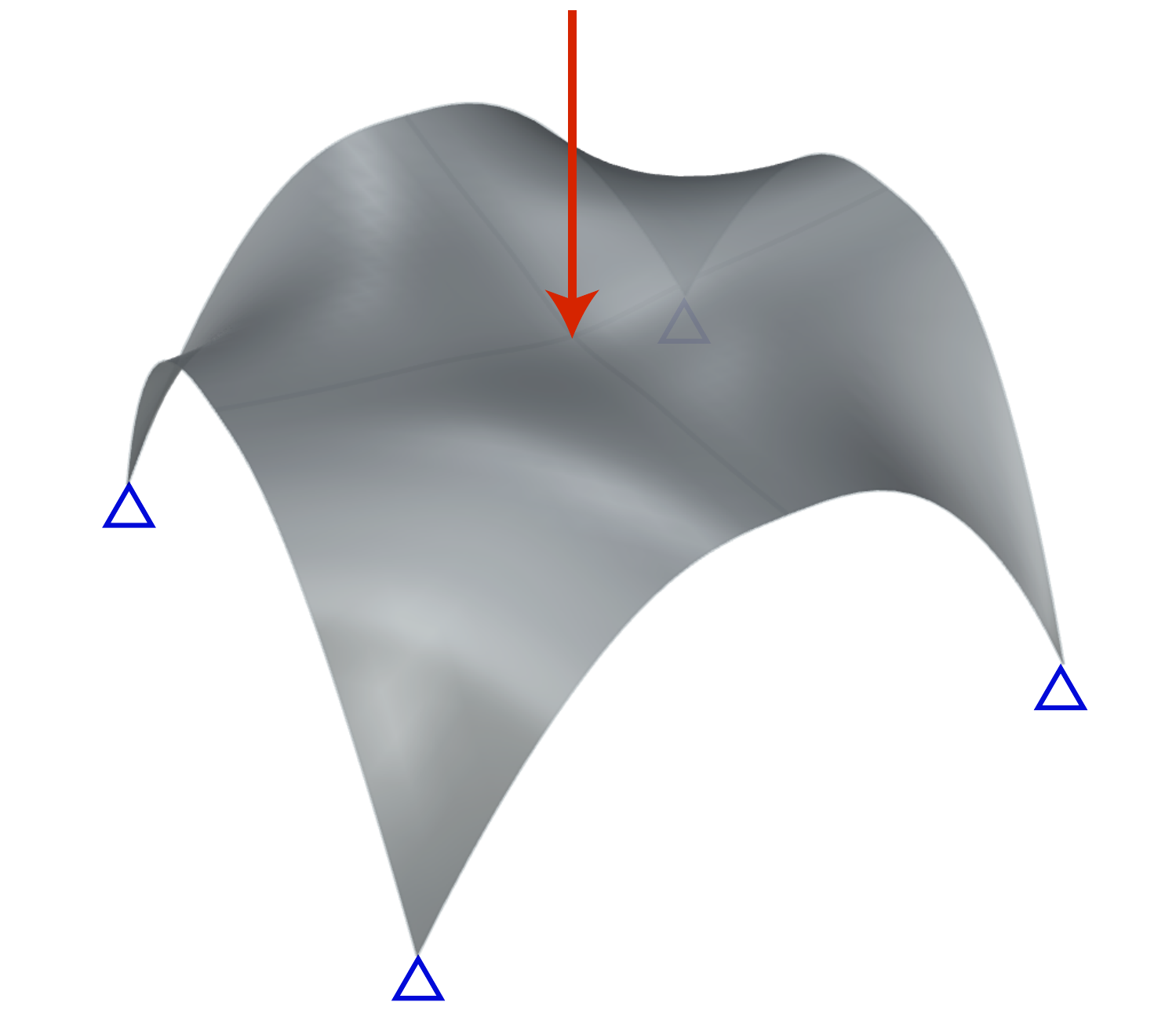}
					\put(50,75){\scriptsize{$G = 100$}}
				\end{overpic}\\
			\end{minipage}
			\qquad
			\begin{minipage}[t]{0.2\linewidth}
				\centering
				\begin{overpic}[scale = 0.14]{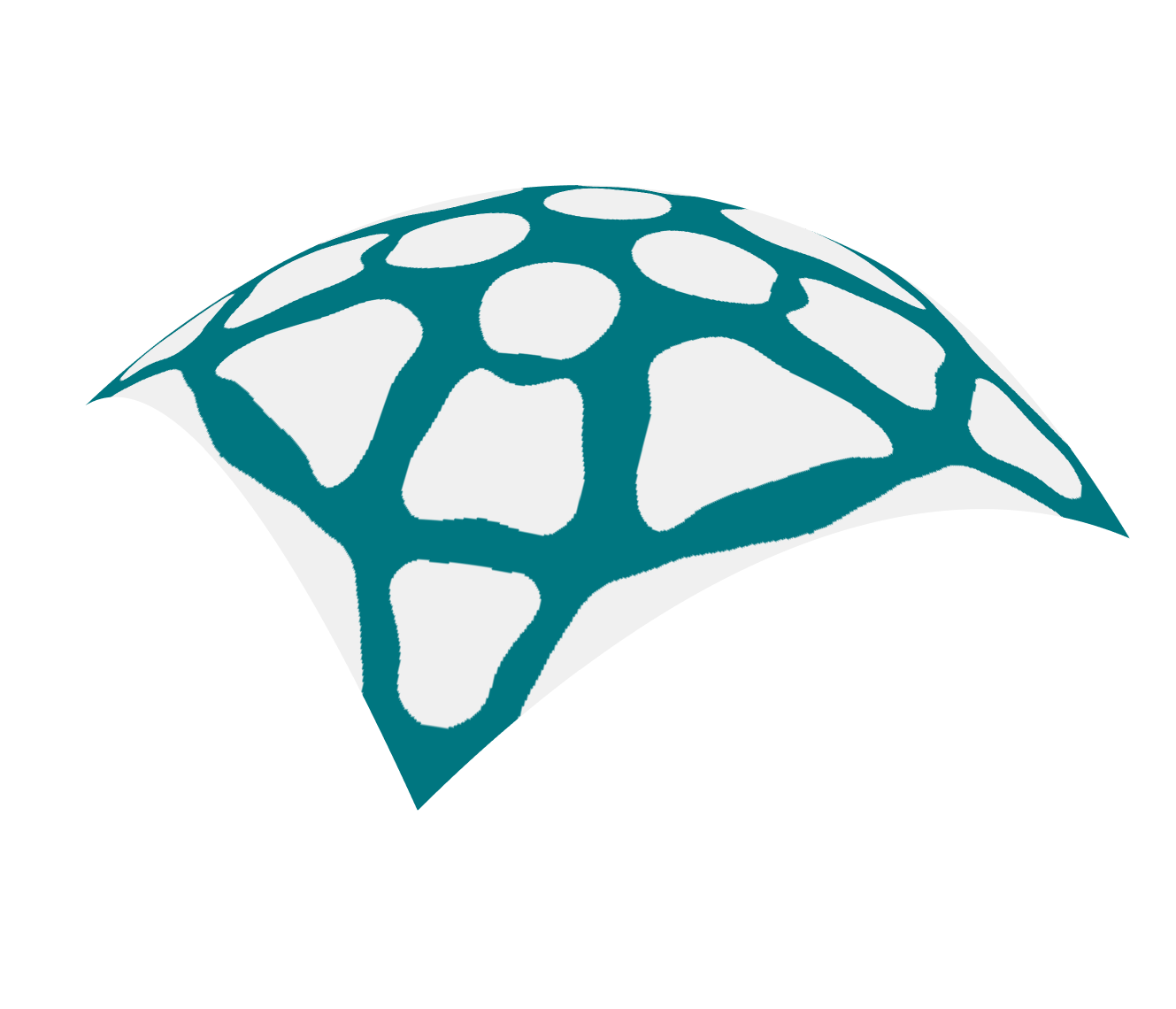}
					\put(20,8){$(20, 100^2, 10, 0.5)$}
					\put(32,82){\scriptsize{$C = 175.35$}}
					\put(28,72){\scriptsize{$V/V_s = 0.434$}}
				\end{overpic}\\
				\vspace{0.2cm}
				\begin{overpic}[scale = 0.14]{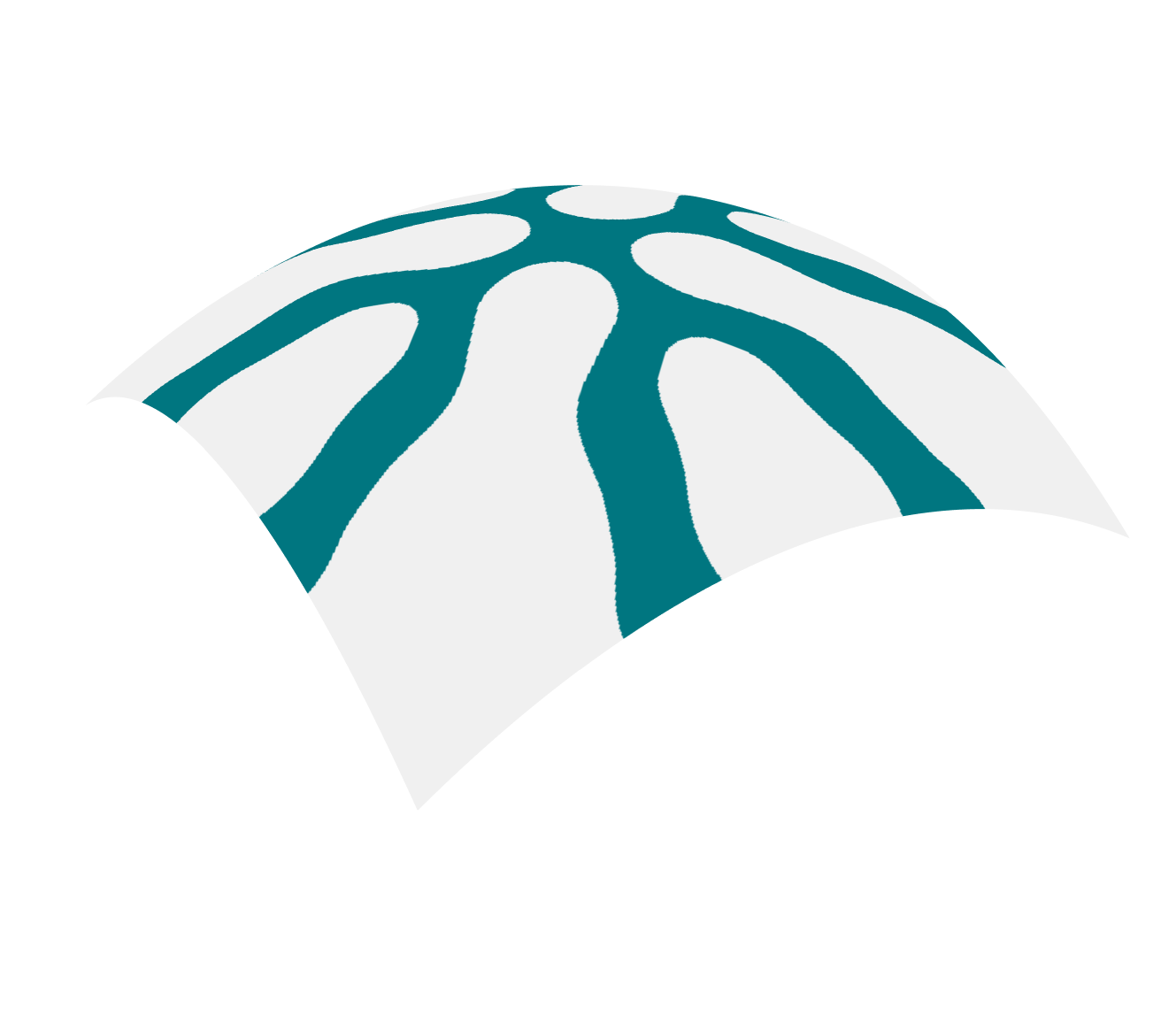}
					\put(20,8){$(15, 100^2, 10, 0.5)$}
					\put(32,82){\scriptsize{$C = 130.33$}}
					\put(28,72){\scriptsize{$V/V_s = 0.305$}}
				\end{overpic}\\
				\vspace{0.2cm}
				\begin{overpic}[scale = 0.14]{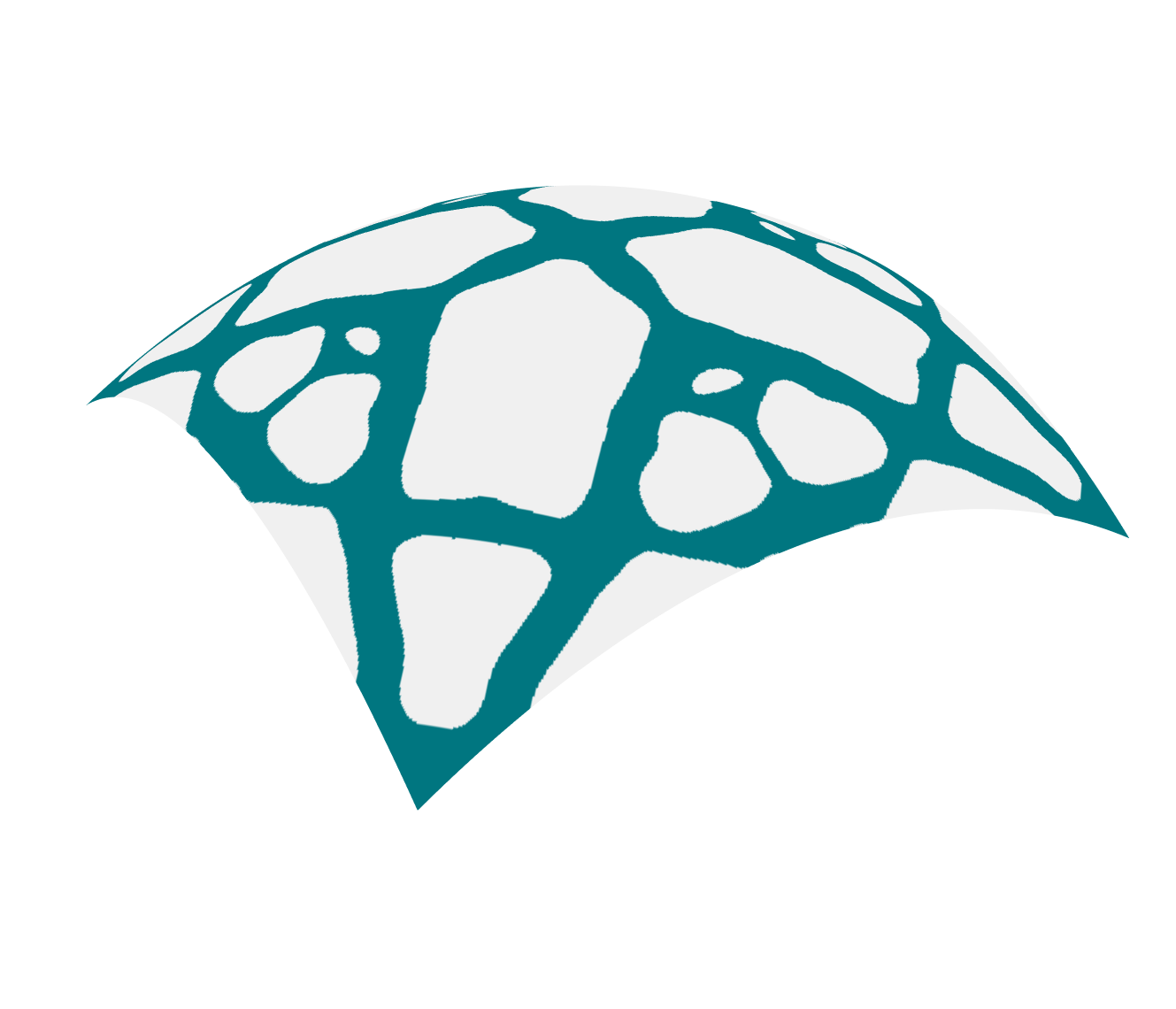}
					\put(20,8){$(30, 100^2, 10, 0.5)$}
					\put(31,82){\scriptsize{$C = 1711.23$}}
					\put(28,72){\scriptsize{$V/V_s = 0.410$}}
				\end{overpic}\\
				\vspace{0.4cm}
				\begin{overpic}[scale = 0.14]{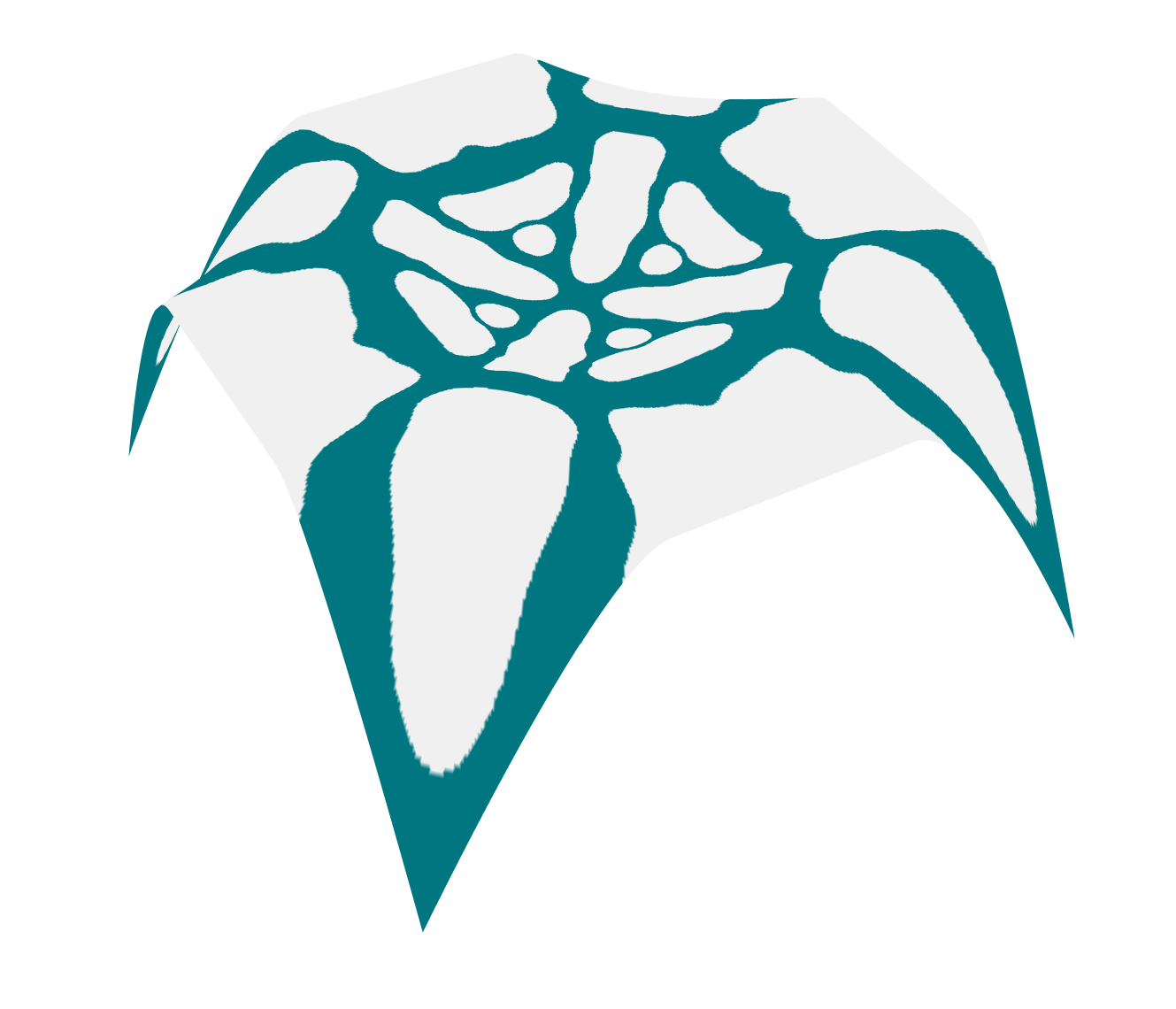}
					\put(20,0){$(25, 100^2, 10, 0.5)$}
					\put(34,93){\scriptsize{$C = 2.46$}}
					\put(28,83){\scriptsize{$V/V_s = 0.422$}}
				\end{overpic}\\
				\vspace{0.02cm}
			\end{minipage}
			\begin{minipage}[t]{0.2\linewidth}
				\centering
				\begin{overpic}[scale = 0.14]{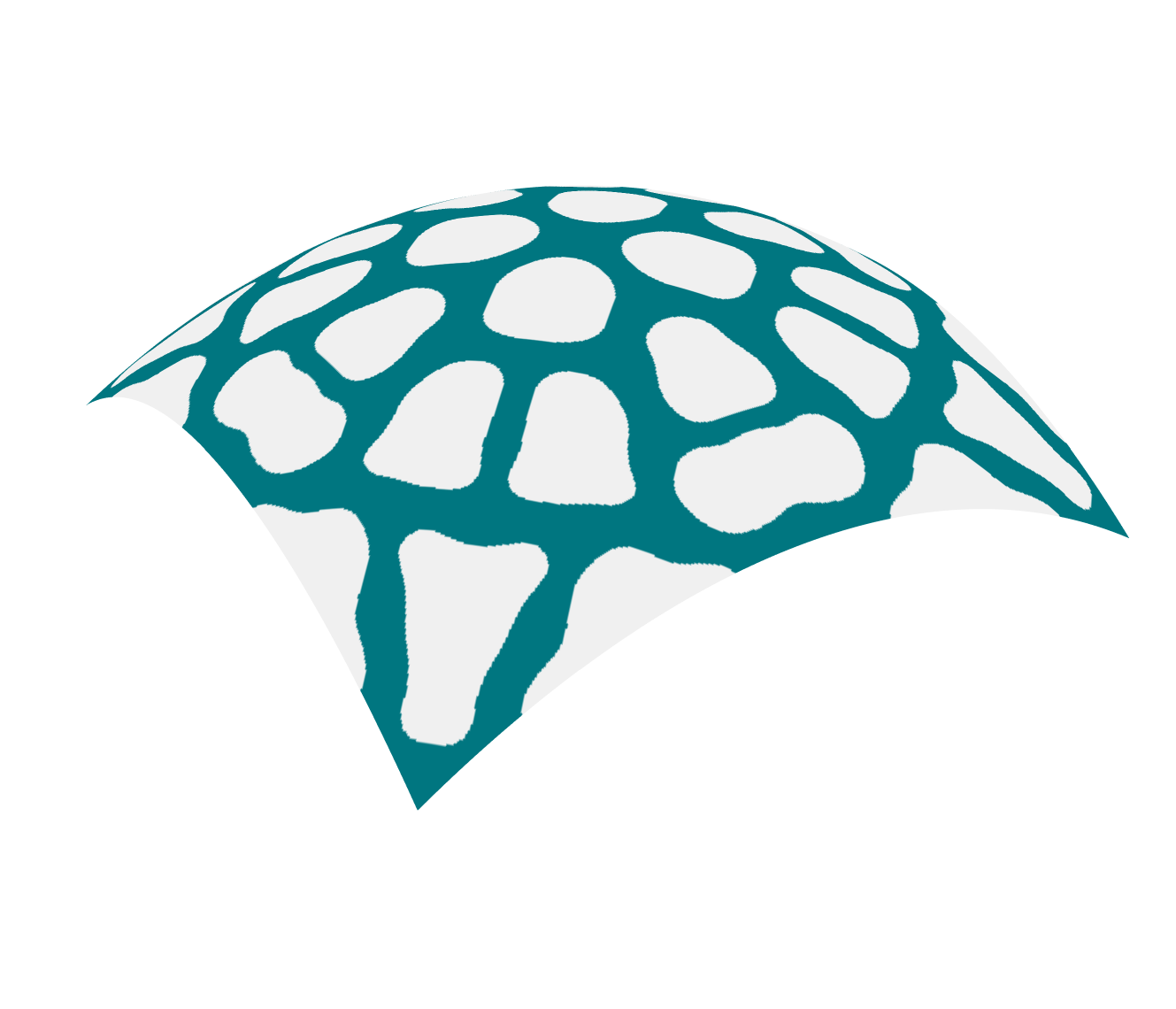}
					\put(20,8){$(25, 100^2, 8, 0.5)$}
					\put(32,82){\scriptsize{$C = 182.94$}}
					\put(28,72){\scriptsize{$V/V_s = 0.445$}}
				\end{overpic}\\
				\vspace{0.2cm}
				\begin{overpic}[scale = 0.14]{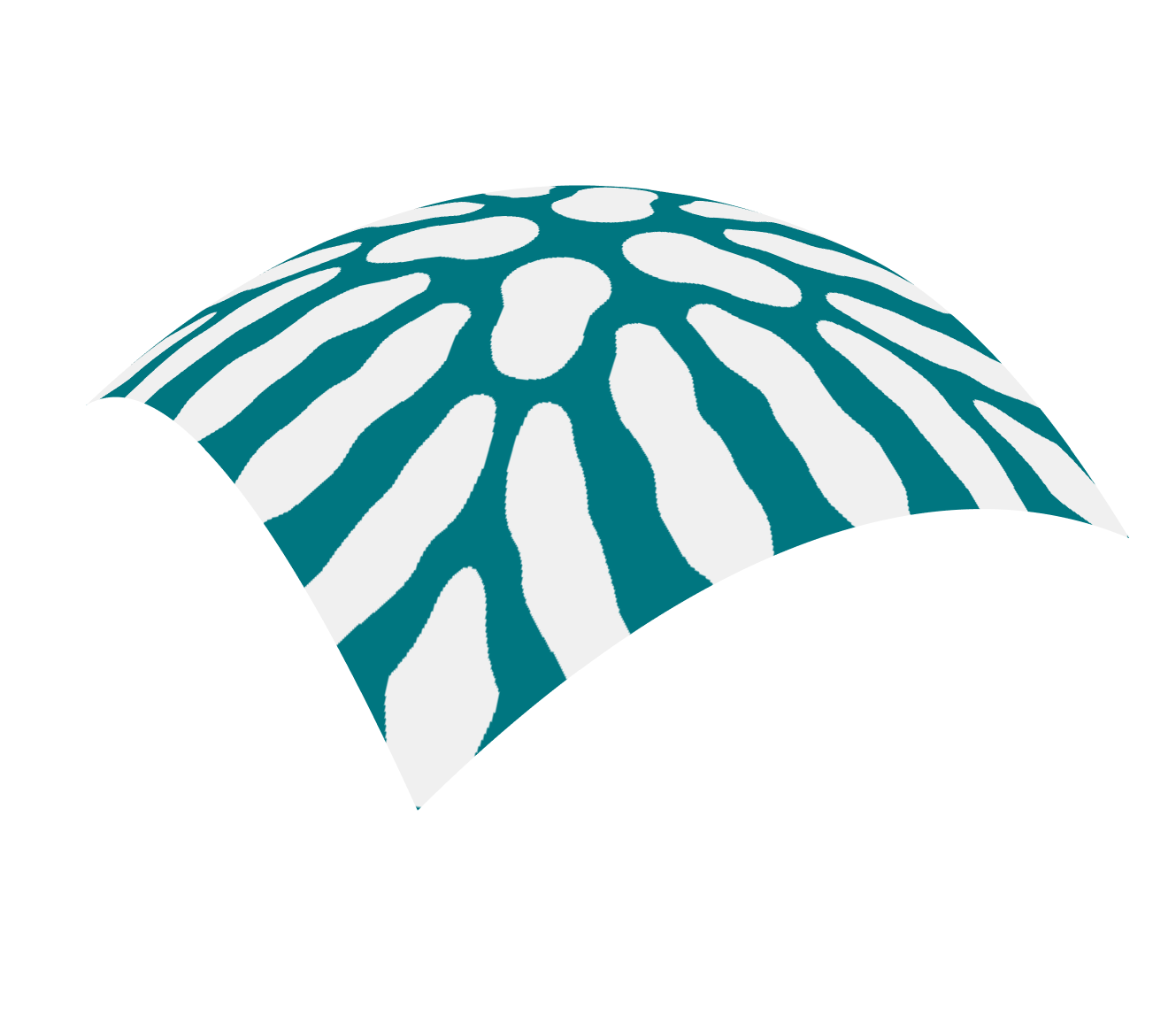}
					\put(20,8){$(30, 100^2, 8, 0.5)$}
					\put(32,82){\scriptsize{$C = 114.72$}}
					\put(28,72){\scriptsize{$V/V_s = 0.441$}}
				\end{overpic}\\
				\vspace{0.2cm}
				\begin{overpic}[scale = 0.14]{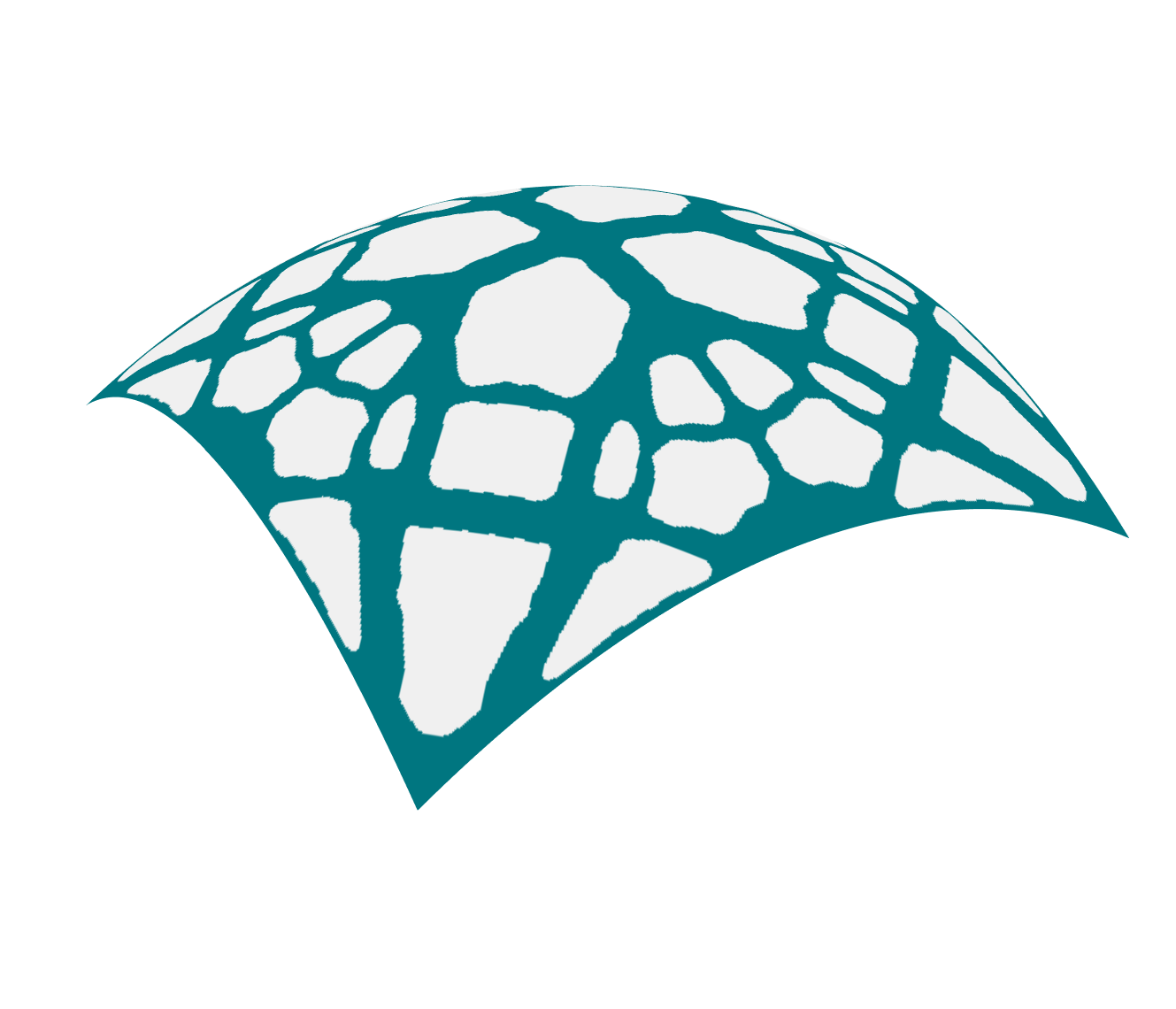}
					\put(20,8){$(50, 100^2, 10, 0.5)$}
					\put(31,82){\scriptsize{$C = 1608.52$}}
					\put(28,72){\scriptsize{$V/V_s = 0.460$}}
				\end{overpic}\\
				\vspace{0.45cm}
				\begin{overpic}[scale = 0.14]{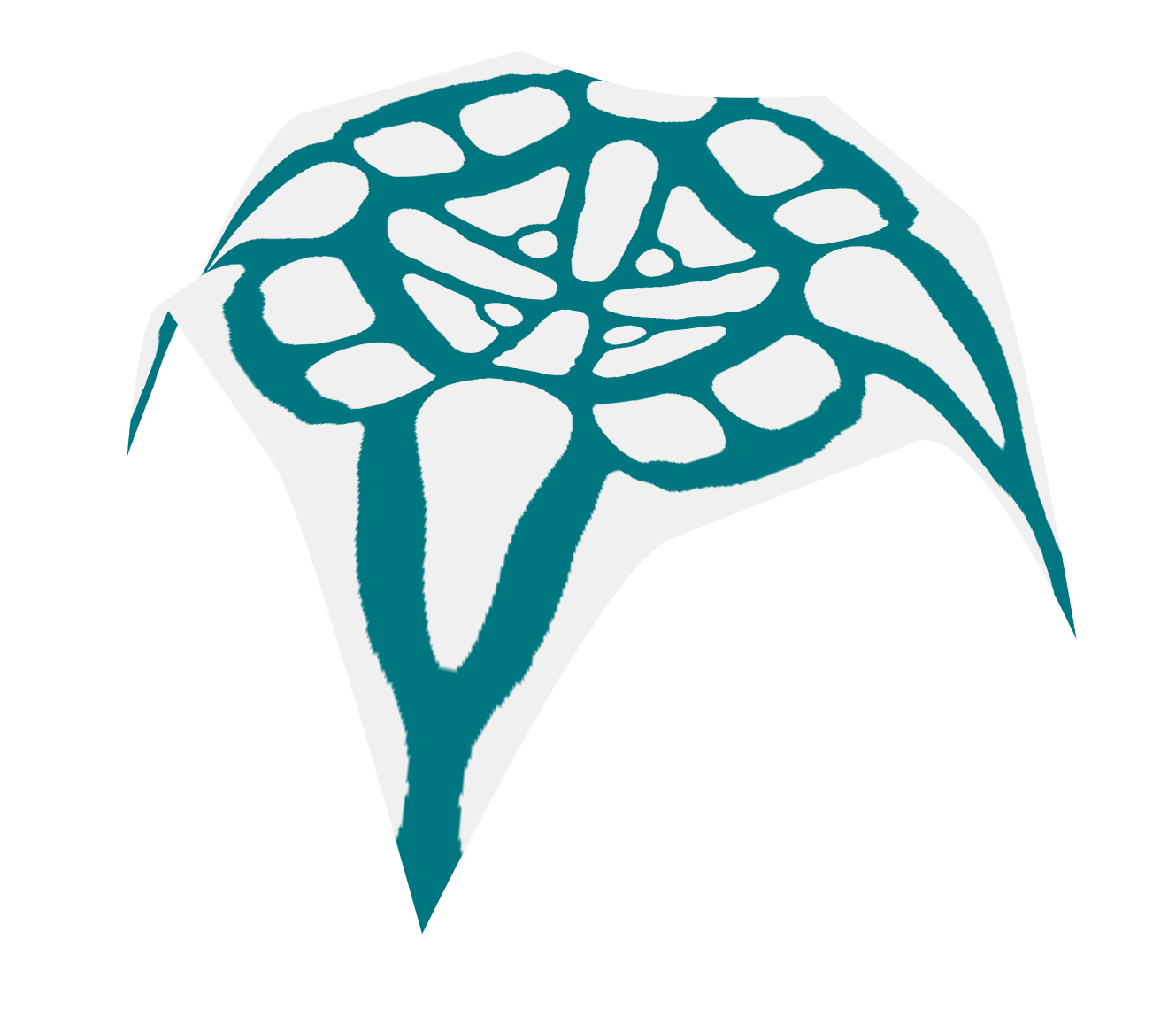}
					\put(20,0){$(30, 100^2, 10, 0.5)$}
					\put(34,93){\scriptsize{$C = 2.14$}}
					\put(28,83){\scriptsize{$V/V_s = 0.433$}}
				\end{overpic}\\
				\vspace{0.02cm}
			\end{minipage}
			\begin{minipage}[t]{0.2\linewidth}
				\centering
				\begin{overpic}[scale = 0.14]{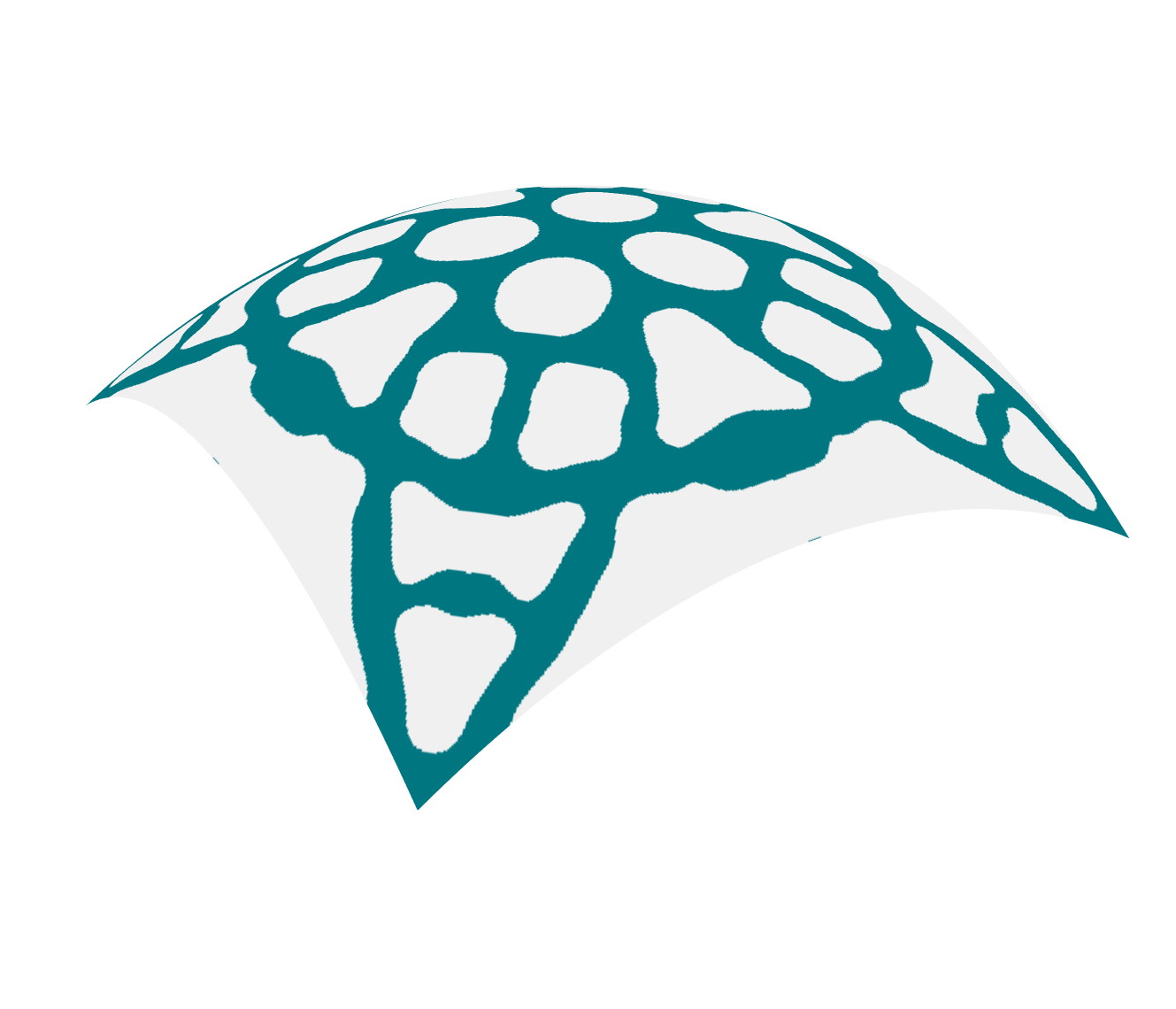}
					\put(20,8){$(30, 100^2, 8, 0.5)$}
					\put(32,82){\scriptsize{$C = 179.23$}}
					\put(28,72){\scriptsize{$V/V_s = 0.420$}}
				\end{overpic}\\
				\vspace{0.2cm}
				\begin{overpic}[scale = 0.14]{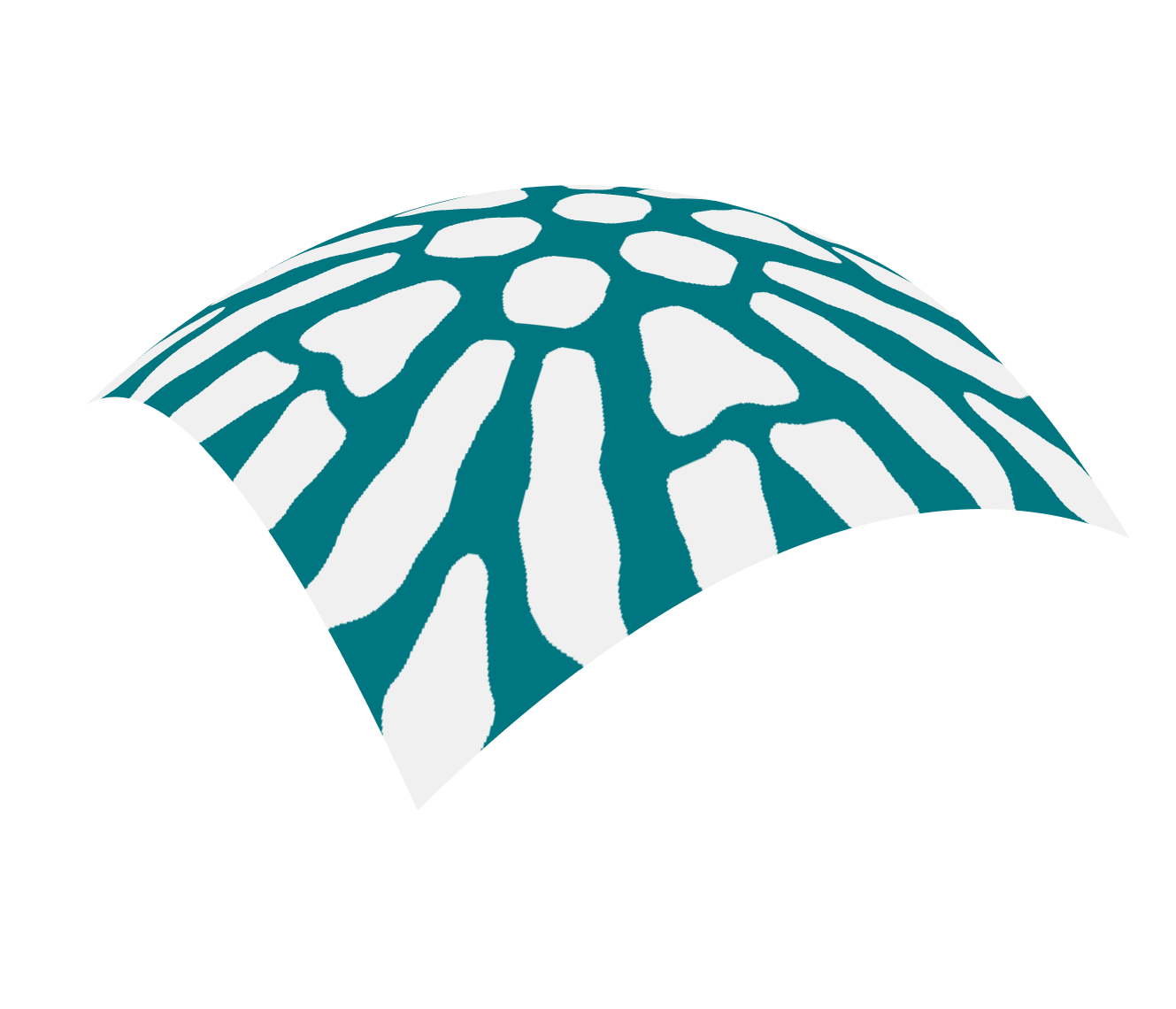}
					\put(20,8){$(35, 100^2, 8, 0.5)$}
					\put(32,82){\scriptsize{$C = 114.44$}}
					\put(28,72){\scriptsize{$V/V_s = 0.443$}}
				\end{overpic}\\
				\vspace{0.2cm}
				\begin{overpic}[scale = 0.14]{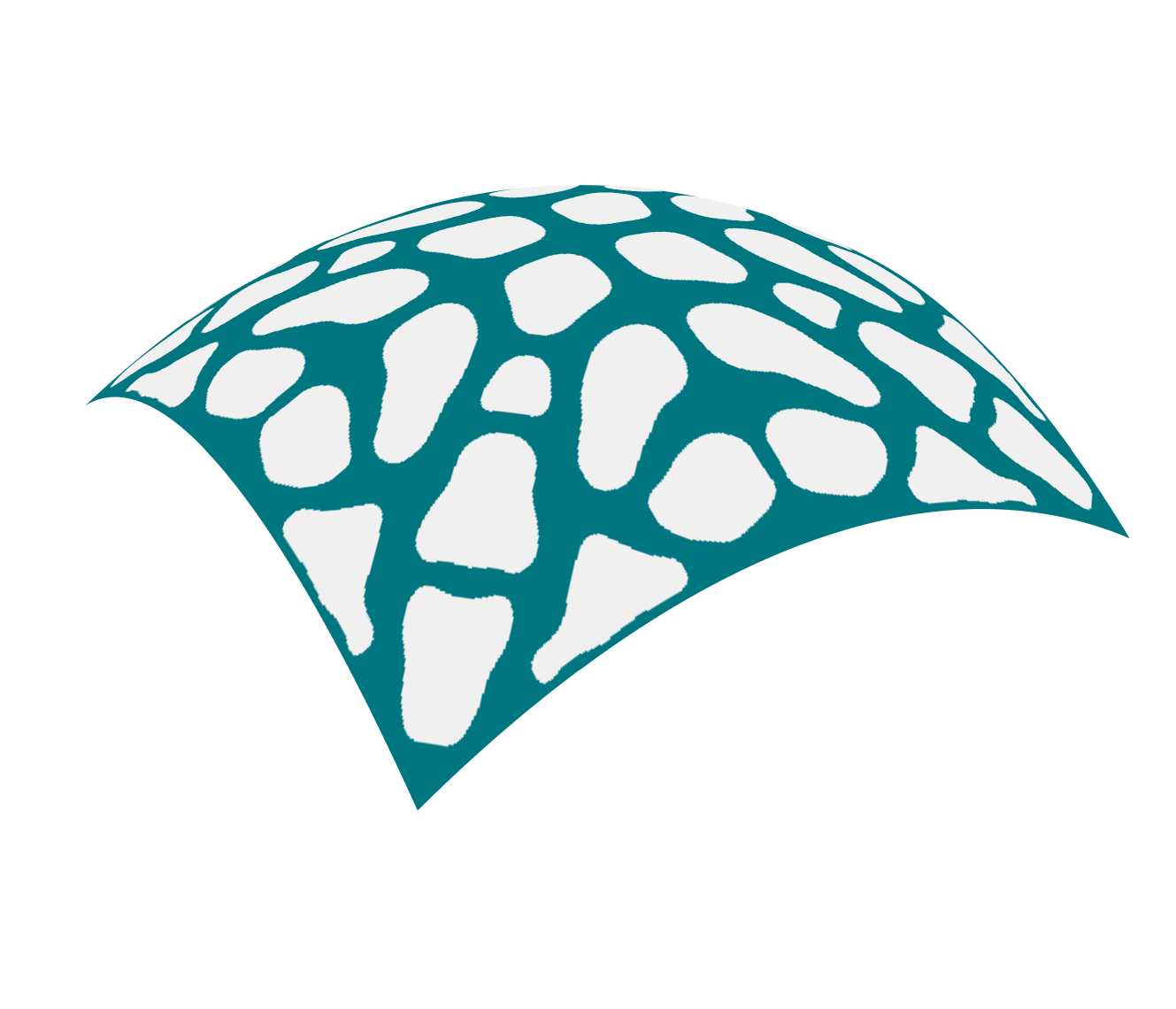}
					\put(20,8){$(30, 100^2, 8, 0.5)$}
					\put(31,82){\scriptsize{$C = 1861.28$}}
					\put(28,72){\scriptsize{$V/V_s = 0.471$}}
				\end{overpic}\\
				\vspace{0.45cm}
				\begin{overpic}[scale = 0.14]{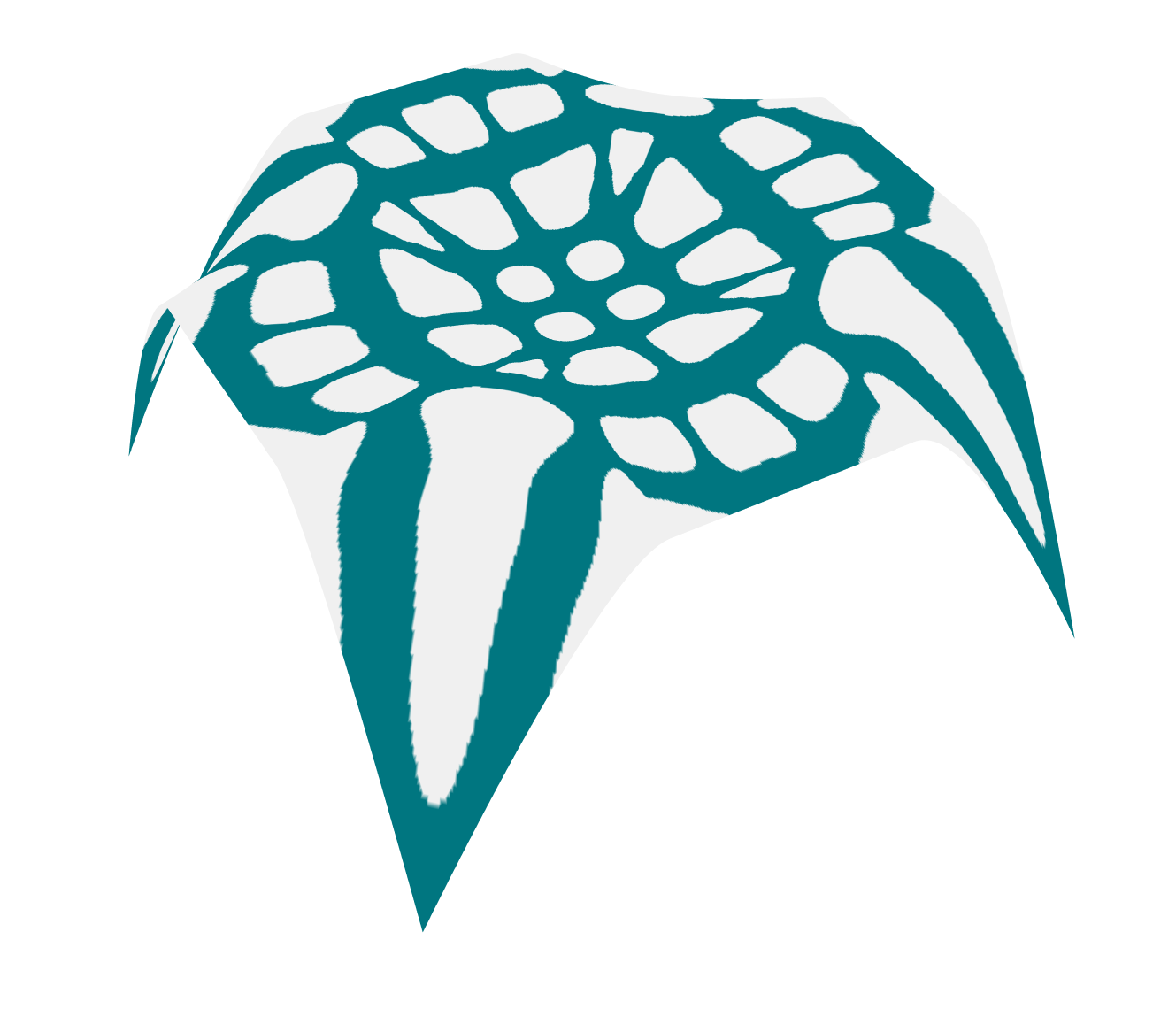}
					\put(20,0){$(35, 200^2, 15, 0.6)$}
					\put(34,93){\scriptsize{$C = 1.61$}}
					\put(28,83){\scriptsize{$V/V_s = 0.539$}}
				\end{overpic}\\
				\vspace{0.02cm}
			\end{minipage}
			\begin{minipage}[t]{0.2\linewidth}
				\centering
				\begin{overpic}[scale = 0.14]{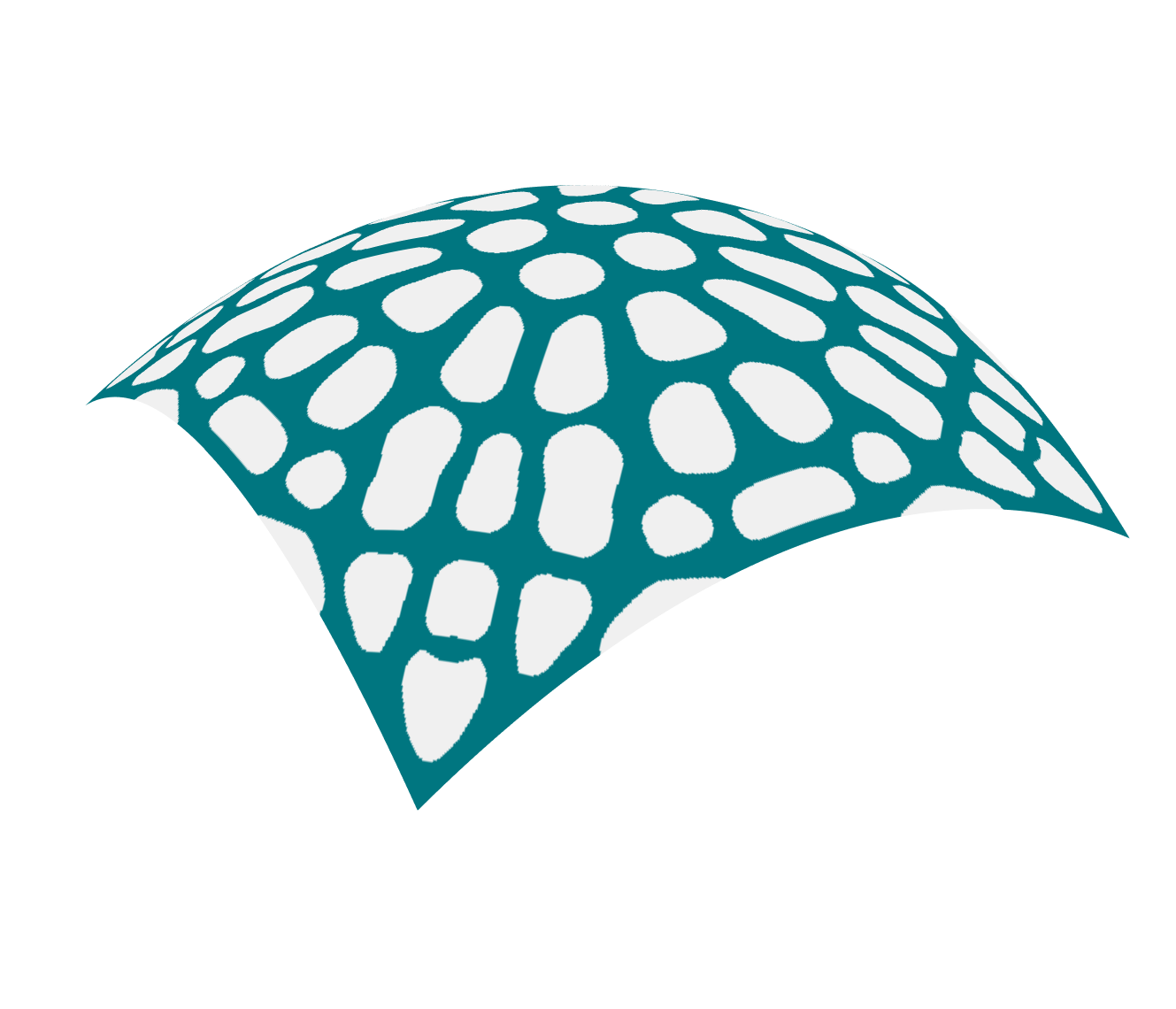}
					\put(20,8){$(30, 100^2, 6, 0.5)$}
					\put(32,82){\scriptsize{$C = 199.47$}}
					\put(28,72){\scriptsize{$V/V_s = 0.462$}}
				\end{overpic}\\
				\vspace{0.2cm}
				\begin{overpic}[scale = 0.14]{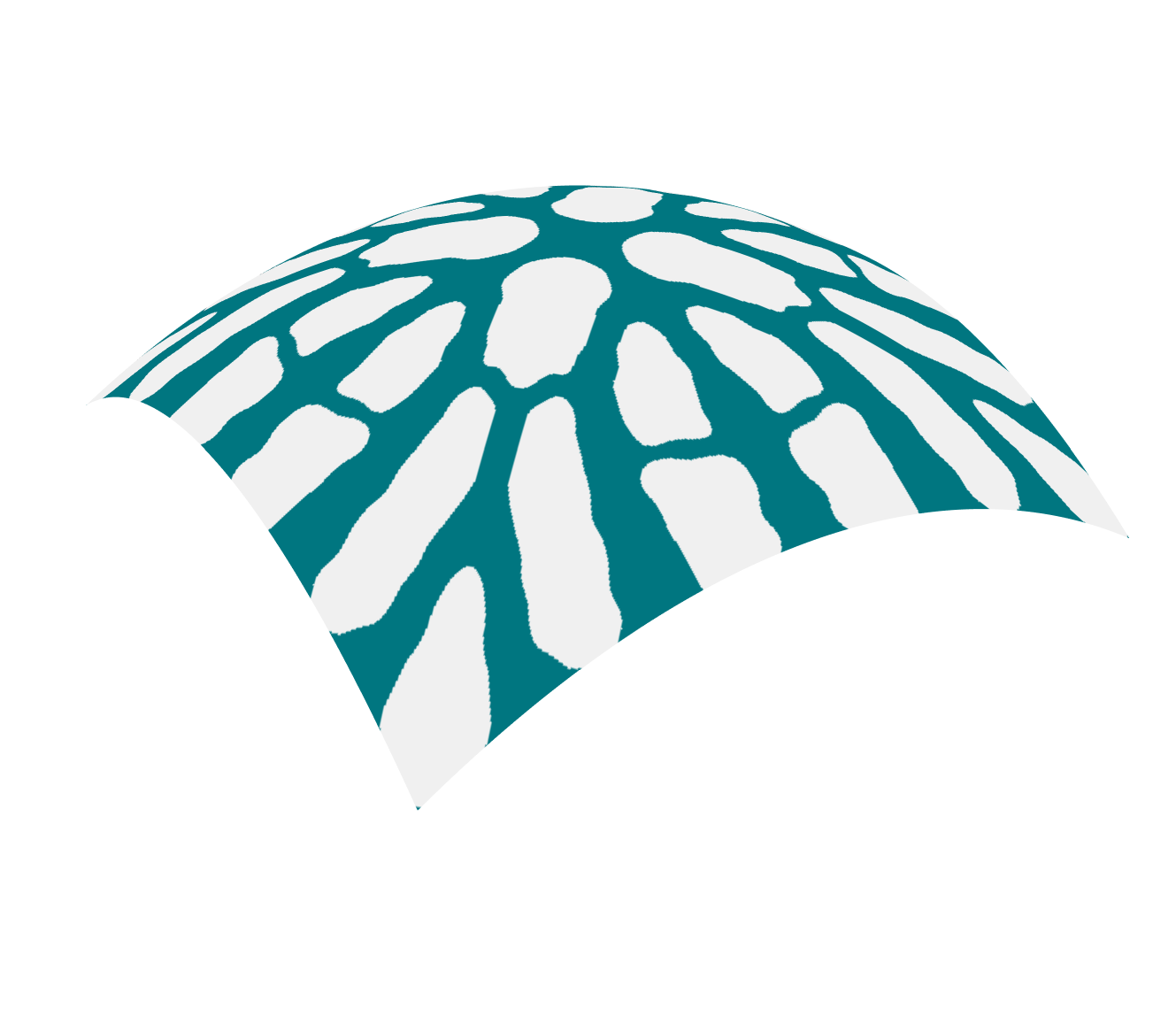}
					\put(20,8){$(40, 100^2, 8, 0.5)$}
					\put(32,82){\scriptsize{$C = 112.24$}}
					\put(28,72){\scriptsize{$V/V_s = 0.440$}}
				\end{overpic}\\
				\vspace{0.2cm}
				\begin{overpic}[scale = 0.14]{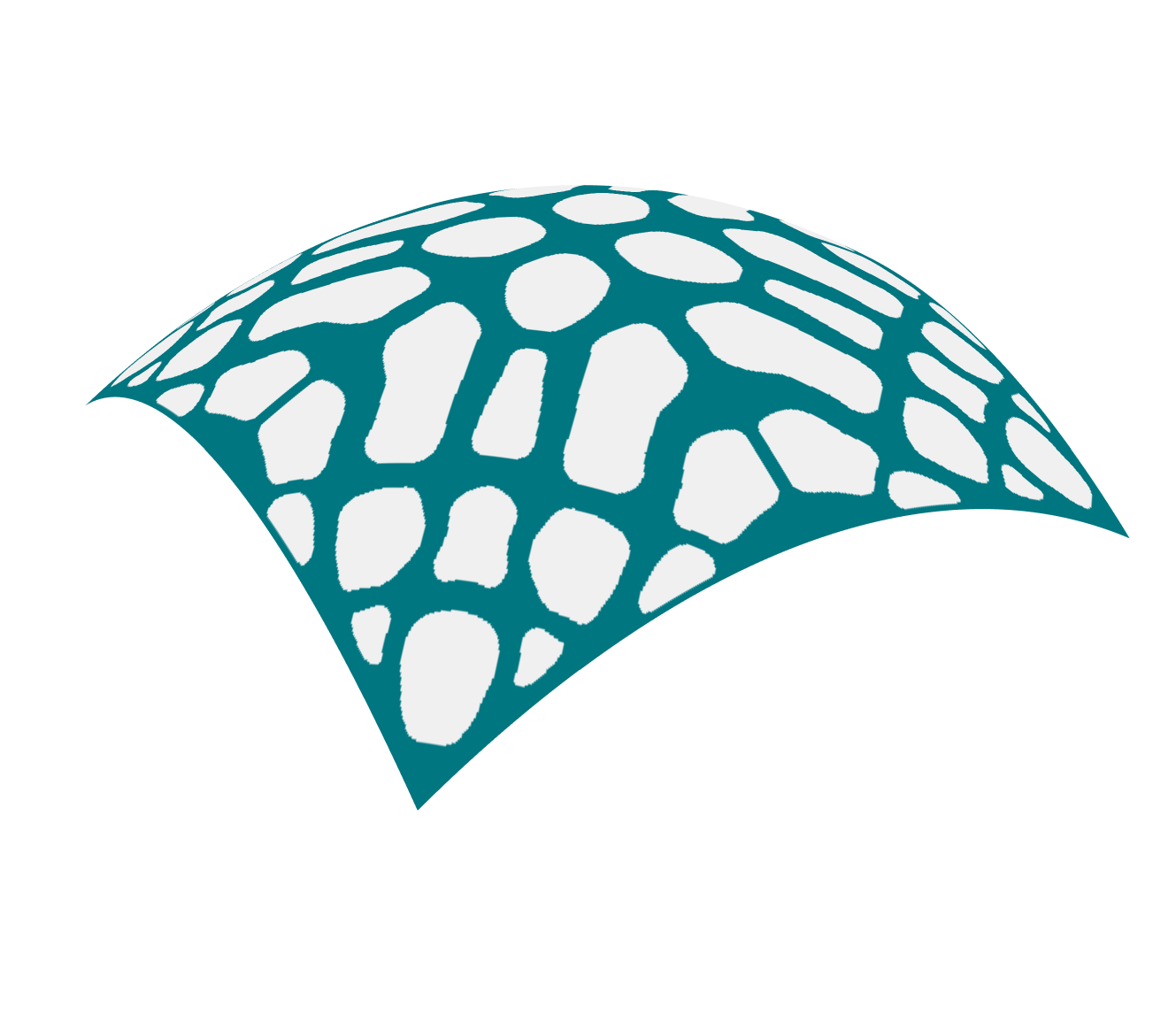}
					\put(20,8){$(35, 100^2, 8, 0.5)$}
					\put(31,82){\scriptsize{$C = 1846.46$}}
					\put(28,72){\scriptsize{$V/V_s = 0.477$}}
				\end{overpic}\\
				\vspace{0.45cm}
				\begin{overpic}[scale = 0.14]{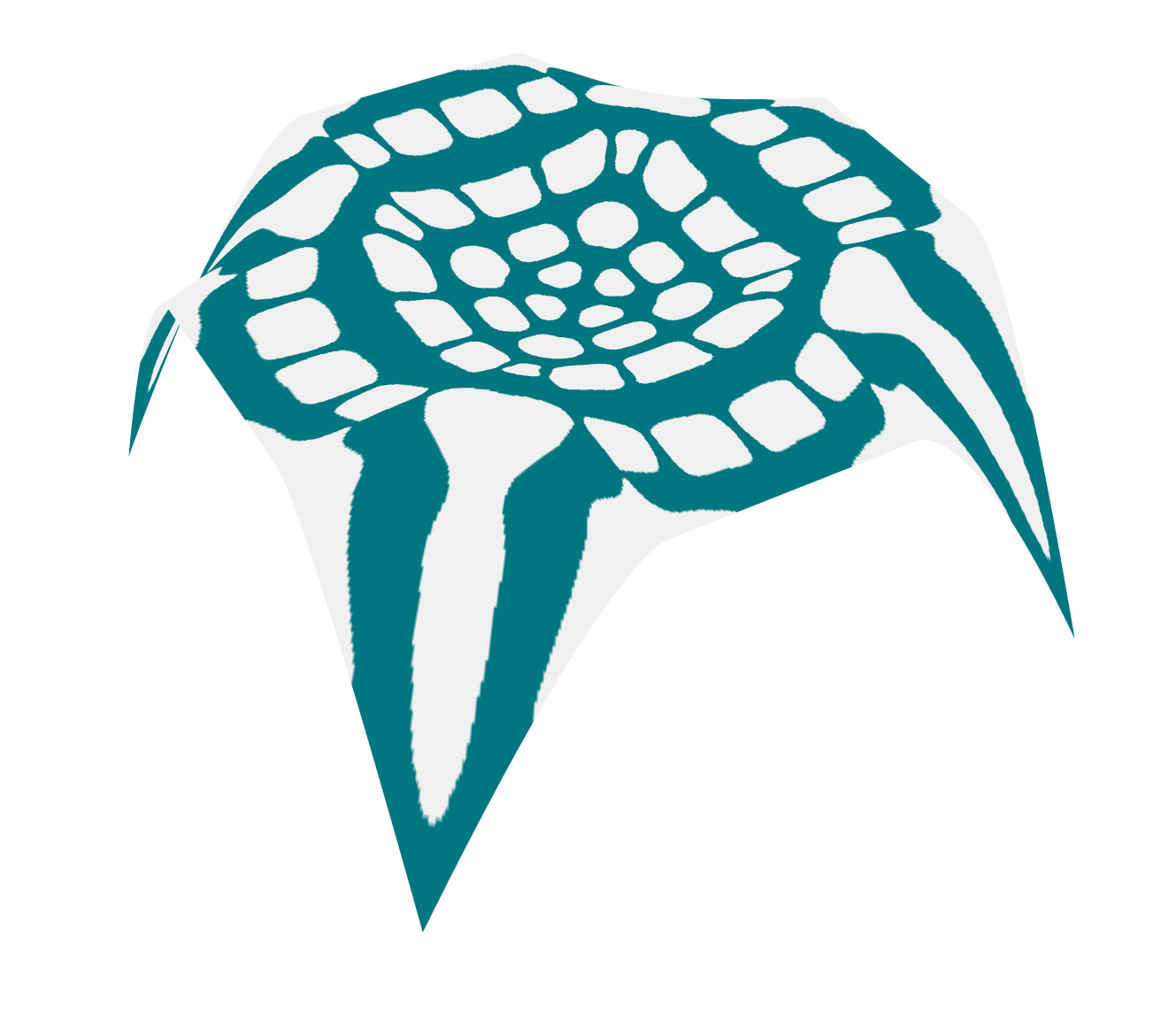}
					\put(20,0){$(38, 200^2, 12, 0.6)$}
					\put(34,93){\scriptsize{$C = 1.79$}}
					\put(28,83){\scriptsize{$V/V_s = 0.537$}}
				\end{overpic}\\
				\vspace{0.02cm}
			\end{minipage}
			\caption{Porous shells obtained using our proposed IGA-SIMP method. The first column depicts the shell models with load and boundary conditions. Other figures represent the optimum structures with control coefficients $N_{DV}$, different numbers of elements $N_E$, different local radius $\mathcal{R}$ and local volume constraint upper bound $\alpha$. Different results are obtained by different combinations of parameters $(N_{DV}, N_E, \mathcal{R},\alpha)$.}
			\label{fig: porous shell examples}
		\end{center}
	\end{figure*}  

	\subsection{IGA-SIMP design: Porous shell structures}
	A variety of porous shell structures are generated based on Model~\eqref{eq: formulation 2}.
	The first column of Fig.\ref{fig: porous shell examples} portrays the shell models with their respective load and boundary conditions. 
	Three examples as before in Fig~\ref{fig: geometrical comparison - IGA-SIMP VS FEM-SIMP} are considered, and the last one is a twisted shell with a thickness $h = 1$.

	By varying the number of control coefficients $N_{DV}$, the number of elements $N_E$, local radius $\mathcal{R}$, and the upper bound on different local volume fractions $\alpha$, we obtain the optimized results including geometries, volume $V$ and compliance $C$ displayed in Fig.~\ref{fig: porous shell examples}. All the final optimal results satisfy the corresponding local volume constraint.
	Increasing the number of design variables $N_{DV}$ tends to reduce compliance, create more holes and disperse the material geometrically, as observed in the last two examples (second and third columns). 
	Smaller $\mathcal{R}$ values produce more and smaller holes but with stricter constraints leading to higher compliance results, which can be observed in the first two examples. The local filter radius $\mathcal{R}$ should be adjusted according to the number of elements $N_E$ in the analysis model, and more elements lead to more uniformly distributed results under the same ratio of $\mathcal{R}/N_E$ due to improved analysis accuracy.
	Larger $\alpha$ values yield stiffer shells with more material. $\alpha$ is a comparable decimal to the volume fraction in the global volume constraint. 
	To design a porous shell structure, $\alpha$ can be selected as a decimal greater than or equal to the volume fraction to achieve specific requirements, with smaller values saving material and larger values producing denser and stiffer structures.

	\subsection{Fairing boundaries discussion}
	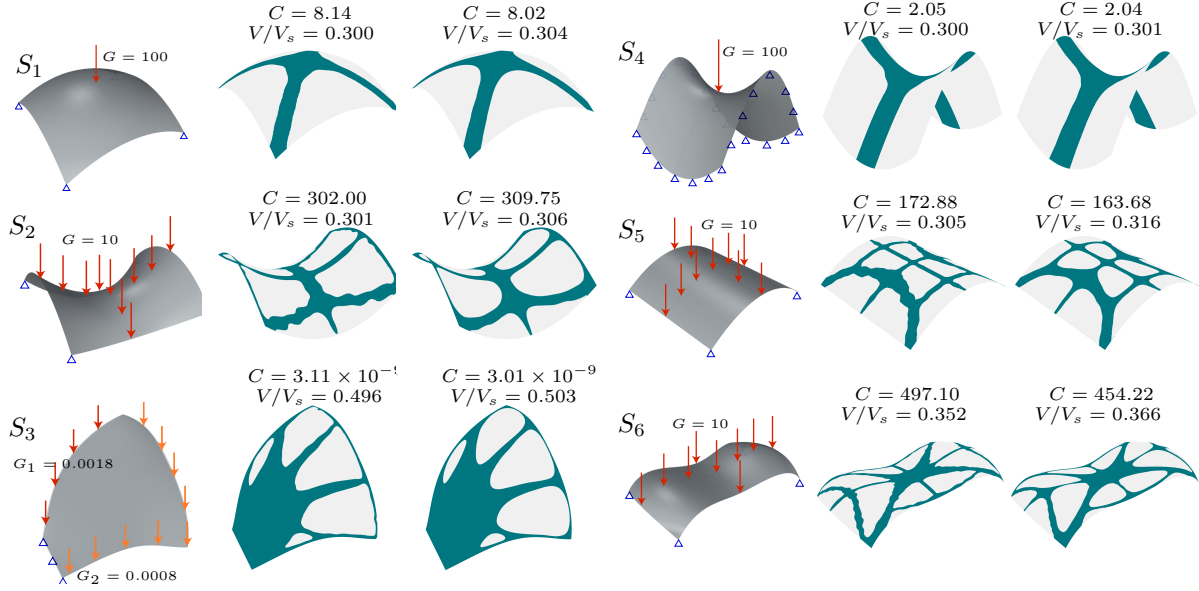
\begin{figure*}[htbp]
		\begin{center}
			\begin{minipage}[t]{0.15\linewidth}
				\centering
				\begin{overpic}[scale = 0.09]{Fig.9-a_s2b1f1.pdf}
					\put(50,77){\tiny{$G = 100$}}
					\put(0,70){$S_1$}
				\end{overpic}\\
				\vspace{0.05cm}
				\begin{overpic}[scale = 0.11]{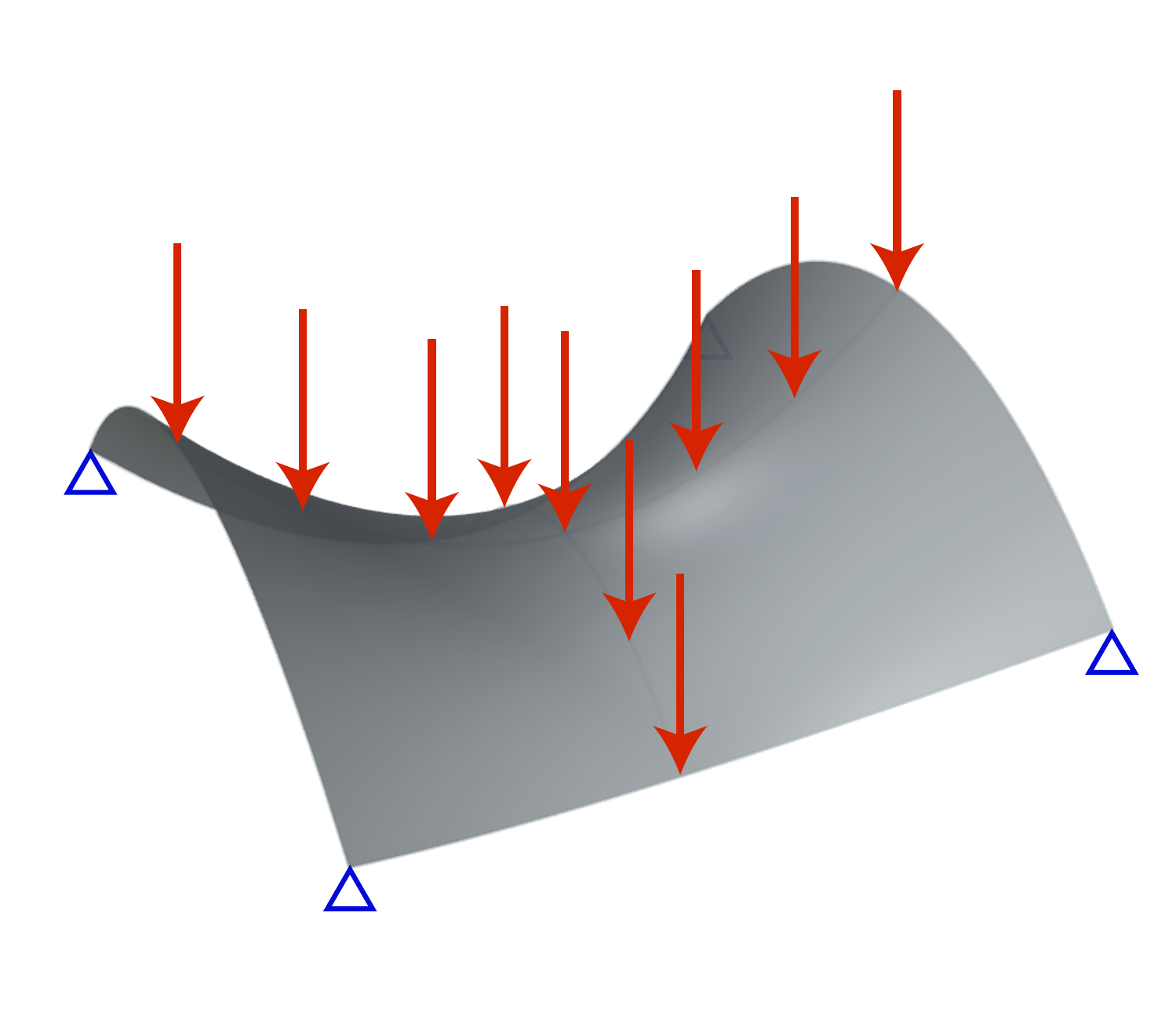}
					\put(25,65){\tiny{$G = 10$}}
					\put(0,70){$S_2$}
				\end{overpic}\\
				\vspace{0.25cm}
				\begin{overpic}[scale = 0.11]{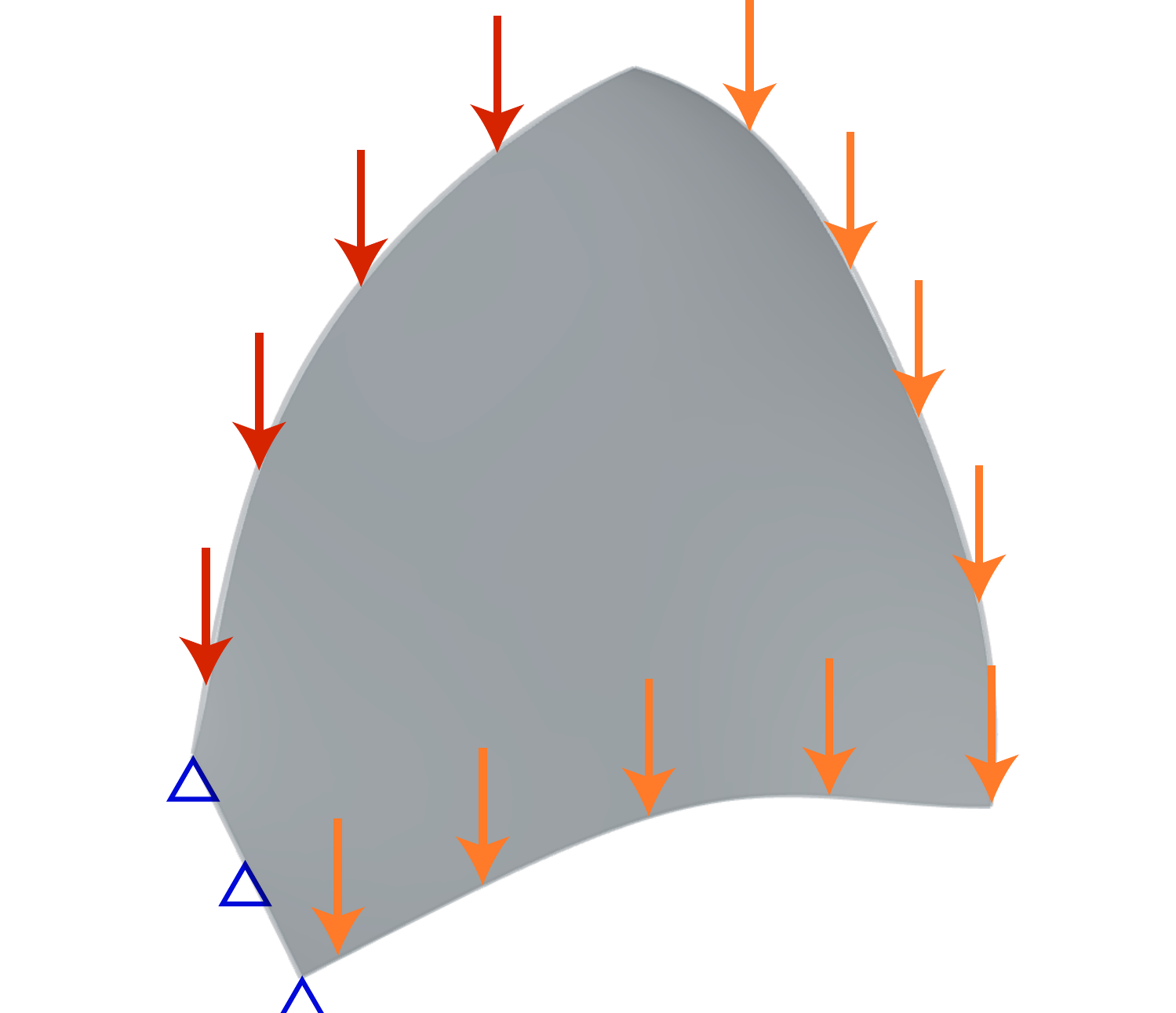}
					\put(2,55){\tiny{$G_1 = 0.0018$}}
					\put(32,2){\tiny{$G_2 = 0.0008$}}
					\put(0,70){$S_3$}
				\end{overpic}
				\\
			\end{minipage}
			%
			\begin{minipage}[t]{0.15\linewidth}
				\centering
				\begin{overpic}[scale = 0.12]{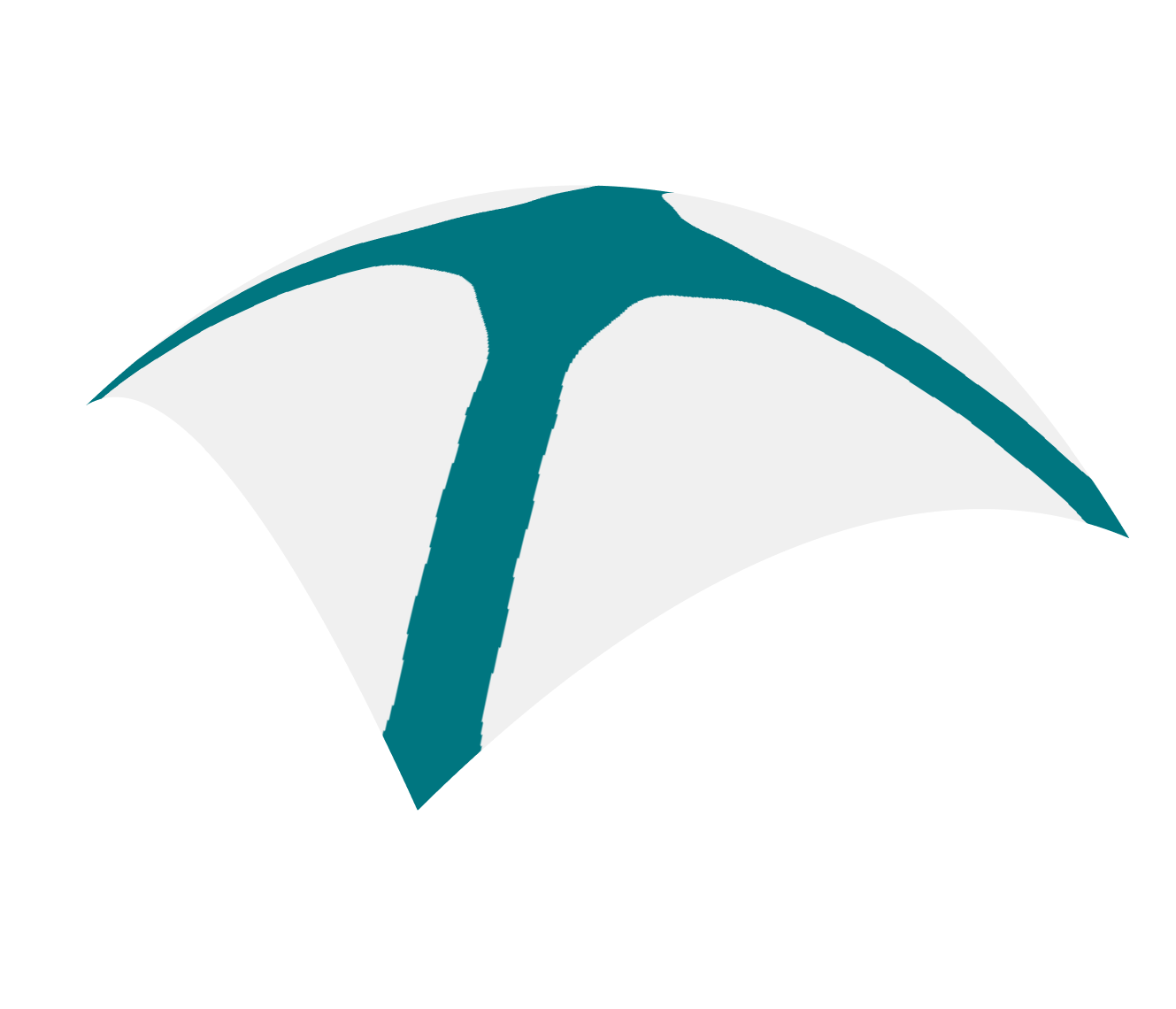}
					\put(22,76){\scriptsize{$V/V_s = 0.300$}}
					\put(32,86){\scriptsize{$C = 8.14$}}
				\end{overpic}\\
				\vspace{0.05cm}
				\begin{overpic}[scale = 0.12]{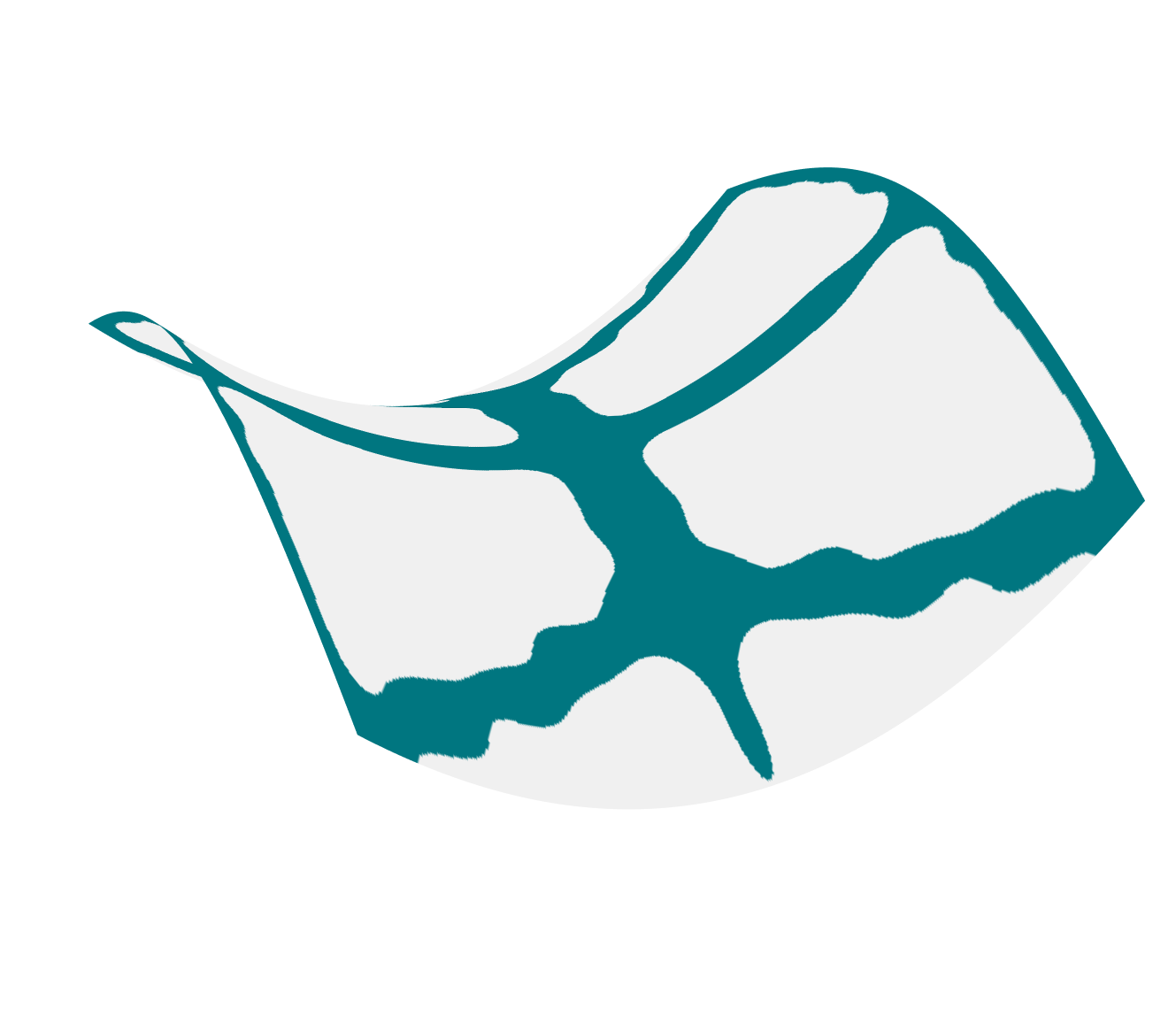}
					\put(22,73){\scriptsize{$V/V_s = 0.301$}}
					\put(29,83){\scriptsize{$C = 302.00$}}
				\end{overpic}\\
				\vspace{0.35cm}
				\begin{overpic}[scale = 0.12]{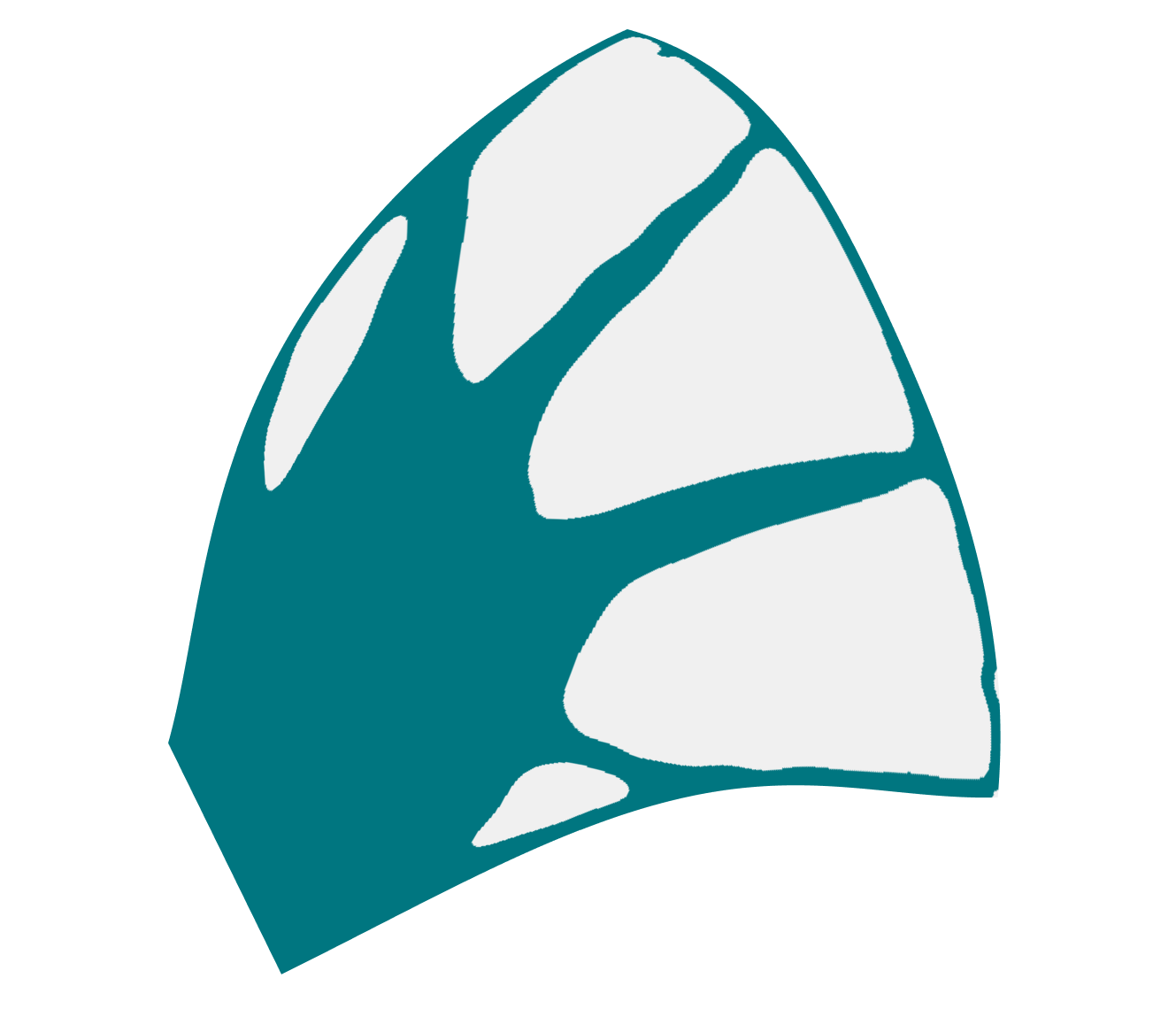}
					\put(26,85){\scriptsize{$V/V_s = 0.496$}}
					\put(22,95){\scriptsize{$C = 3.11\times 10^{-9}$}}
				\end{overpic}\\
			\end{minipage}
			\begin{minipage}[t]{0.15\linewidth}
				\centering
				\begin{overpic}[scale = 0.12]{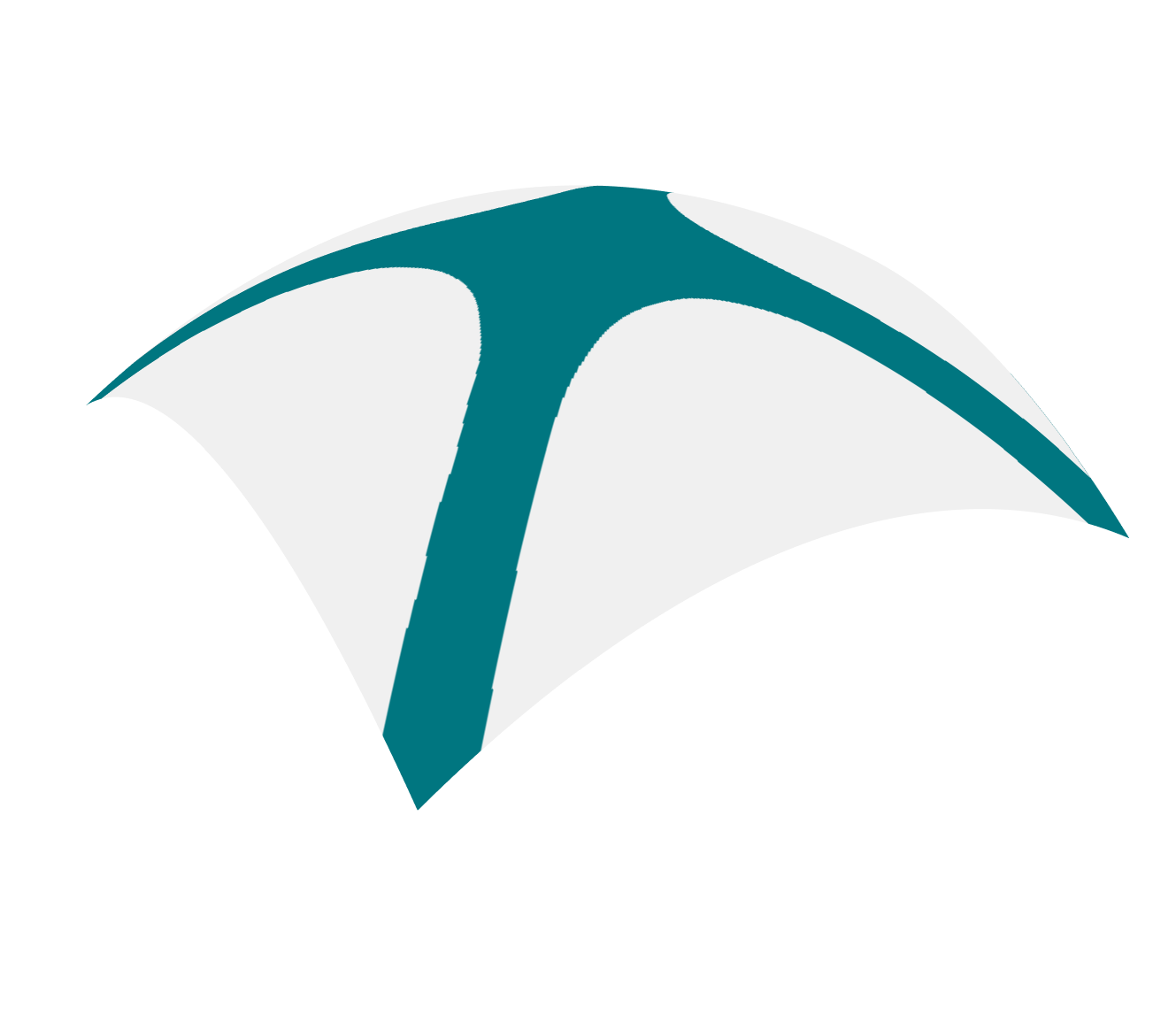}
					\put(22,76){\scriptsize{$V/V_s = 0.304$}}
					\put(32,86){\scriptsize{$C = 8.02$}}
				\end{overpic}\\
				\vspace{0.05cm}
				\begin{overpic}[scale = 0.12]{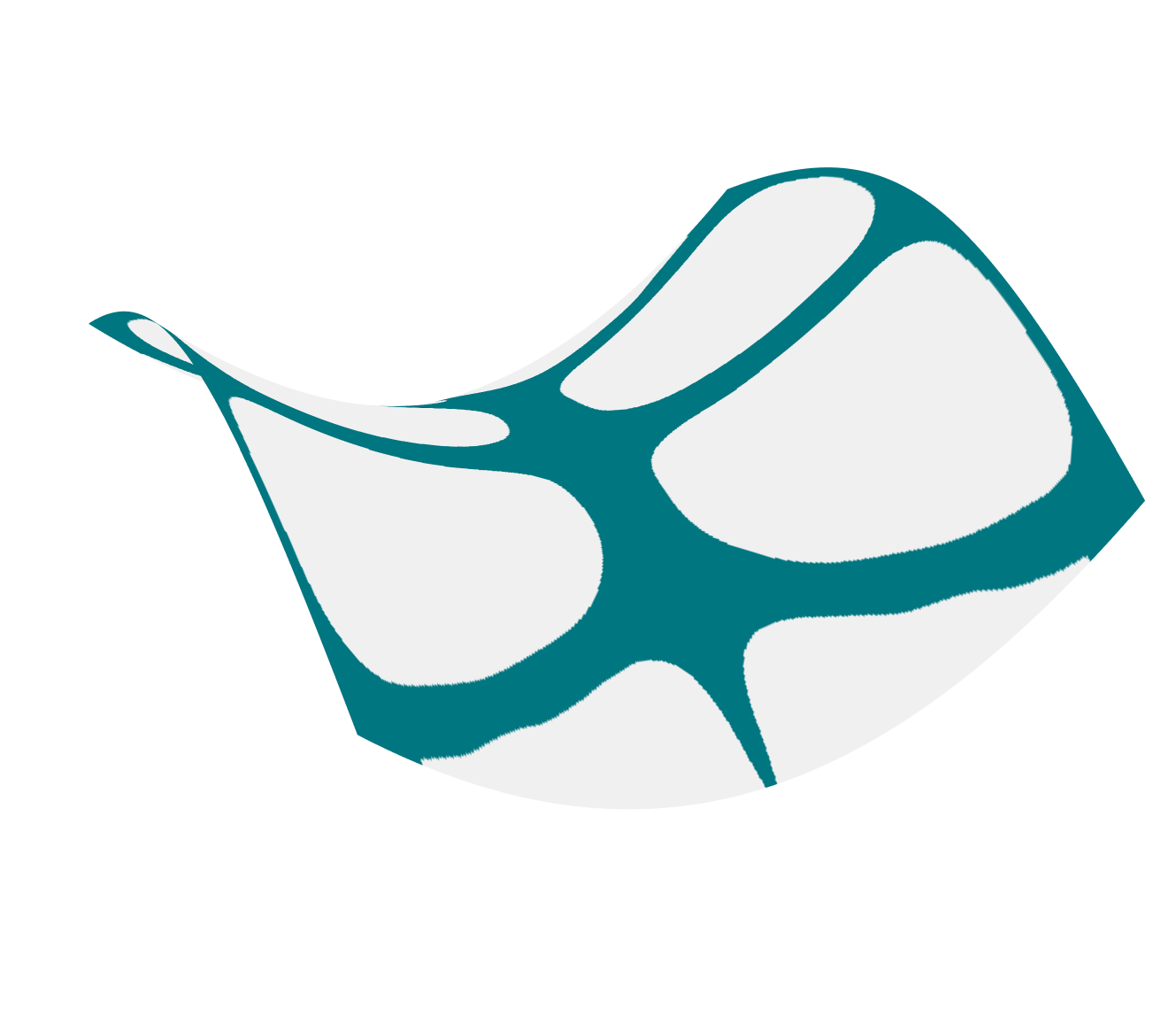}
					\put(22,73){\scriptsize{$V/V_s = 0.306$}}
					\put(29,83){\scriptsize{$C = 309.75$}}
				\end{overpic}\\
				\vspace{0.35cm}
				\begin{overpic}[scale = 0.12]{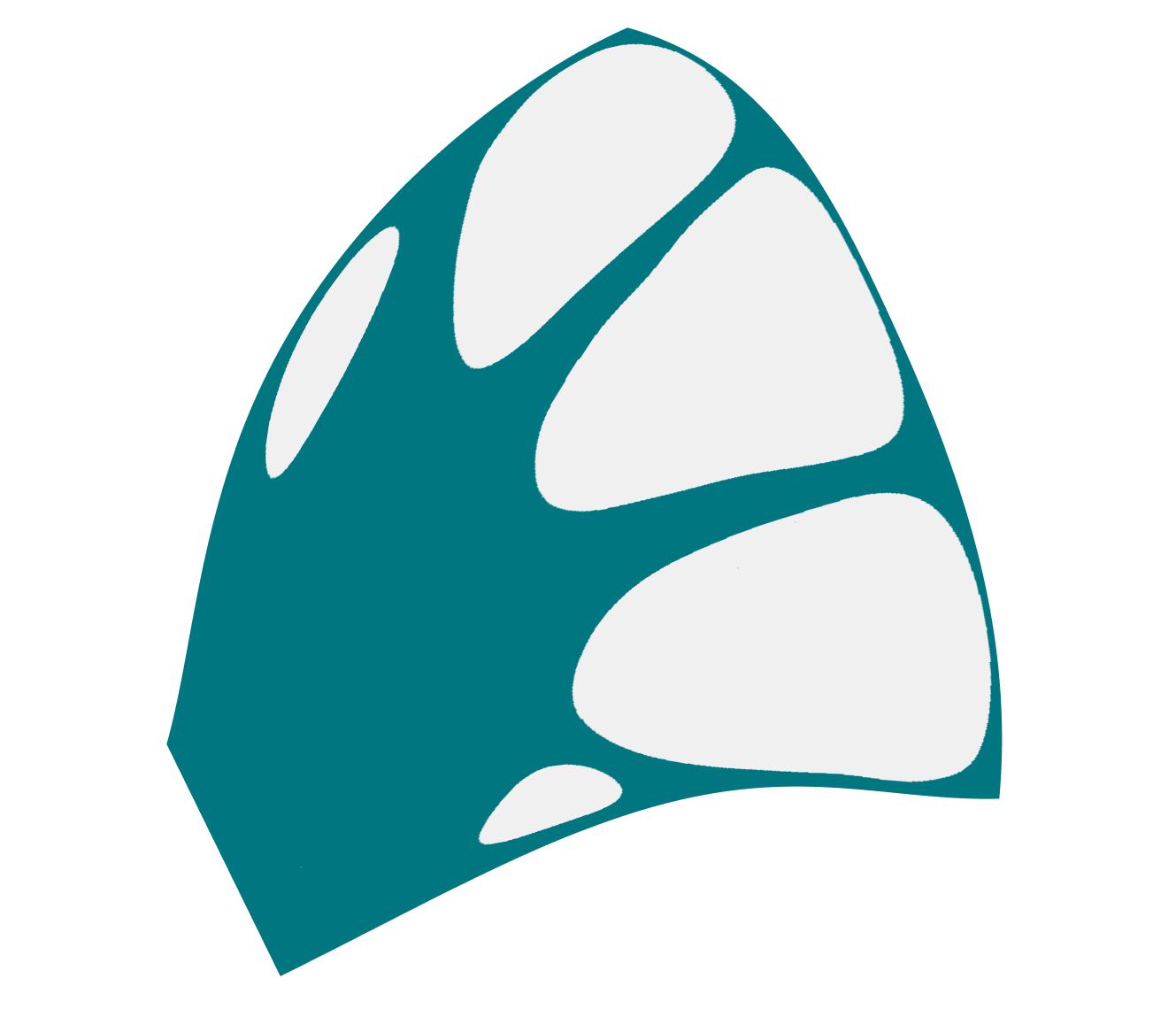}
					\put(26,85){\scriptsize{$V/V_s = 0.503$}}
					\put(22,95){\scriptsize{$C = 3.01\times 10^{-9}$}}
				\end{overpic}\\
			\end{minipage}
			\quad 
			\begin{minipage}[t]{0.15\linewidth}
				\centering
				\begin{overpic}[scale = 0.1]{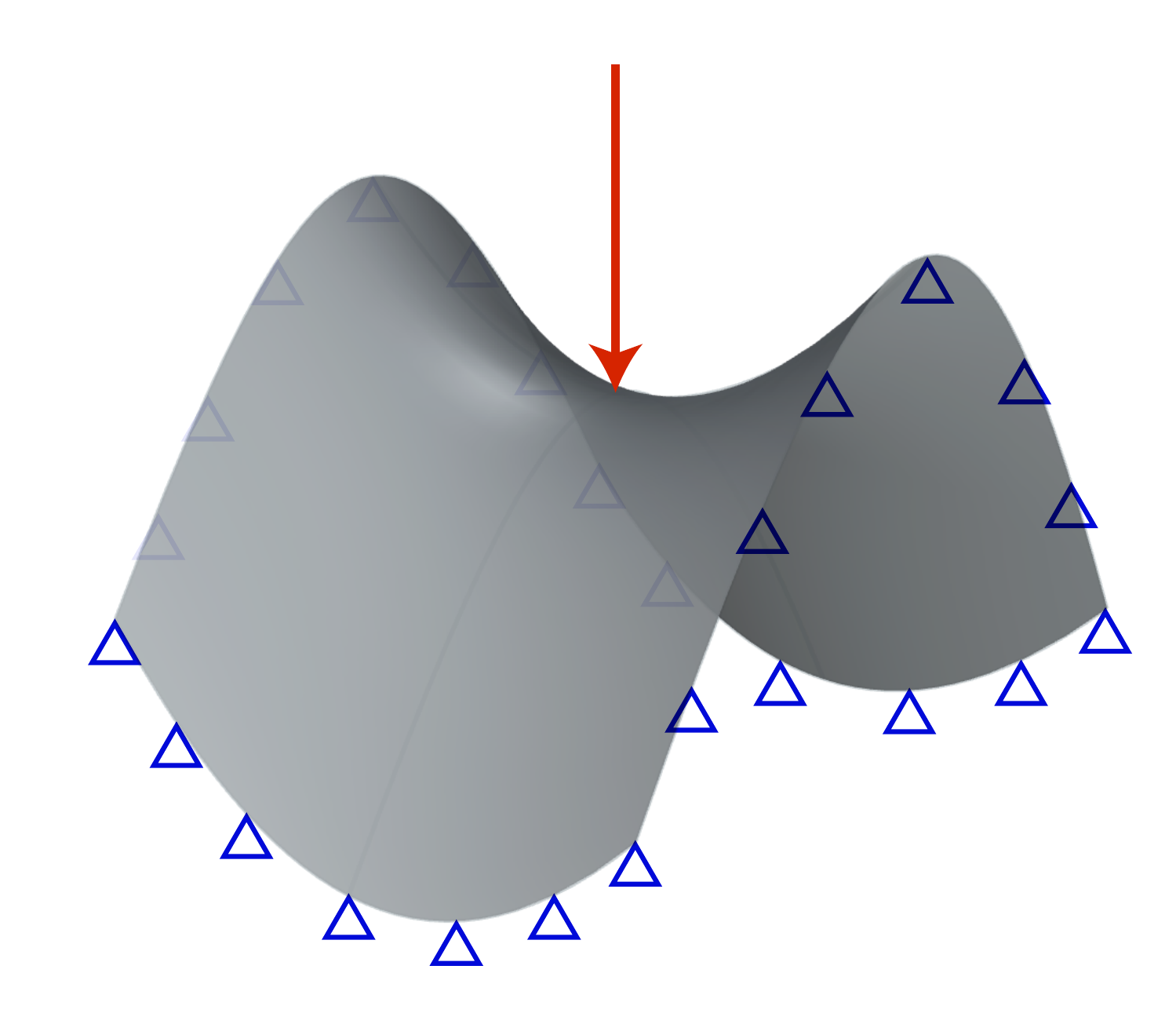}
					\put(54,72){\tiny{$G = 100$}}
					\put(0,70){$S_4$}
				\end{overpic}
				\\
				\vspace{0.15cm}
				\begin{overpic}[scale = 0.1]{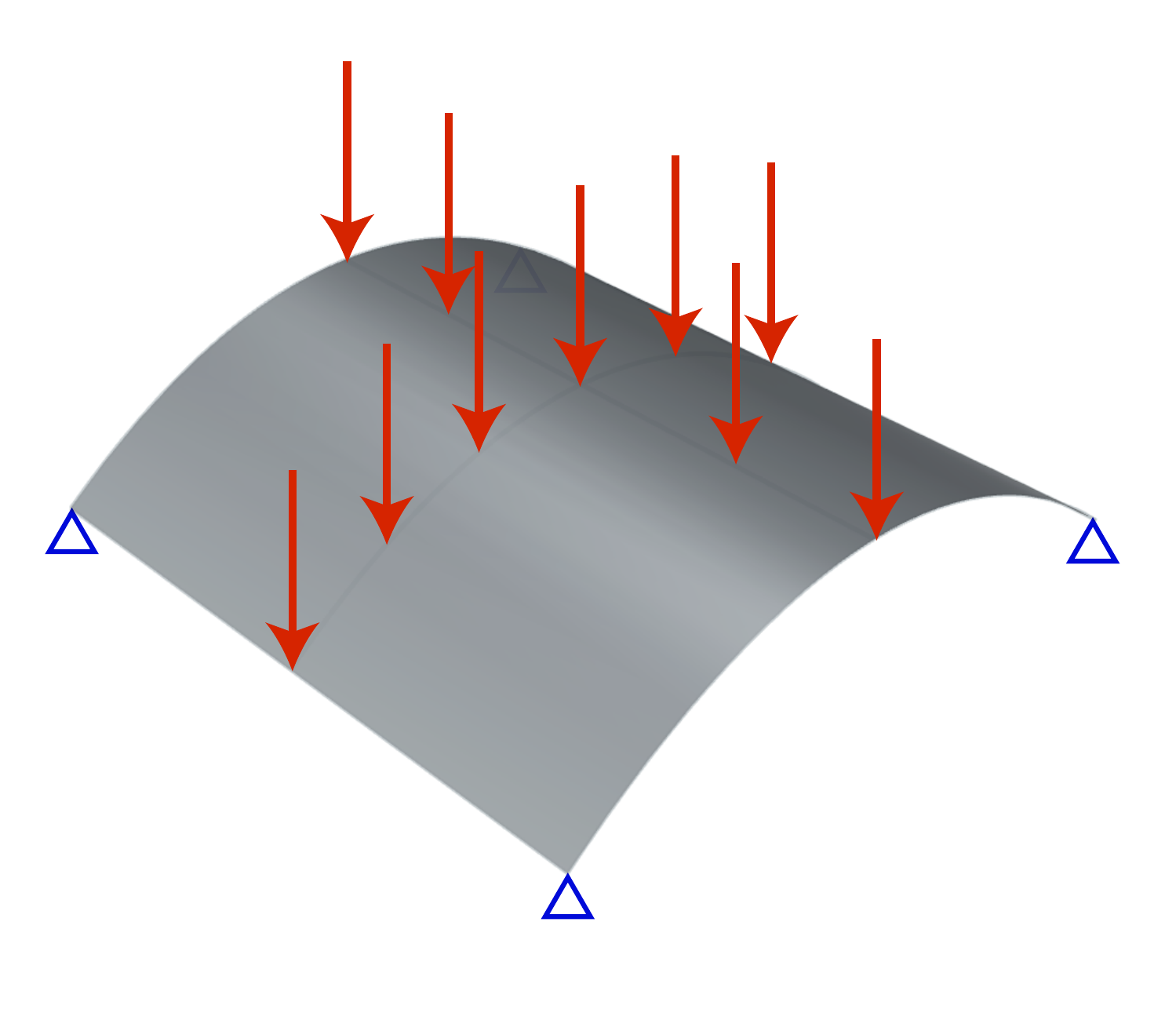}
					\put(43,75){\tiny{$G = 10$}}
					\put(0,70){$S_5$}
				\end{overpic}
				\\
				\vspace{0.35cm}
				\begin{overpic}[scale = 0.1]{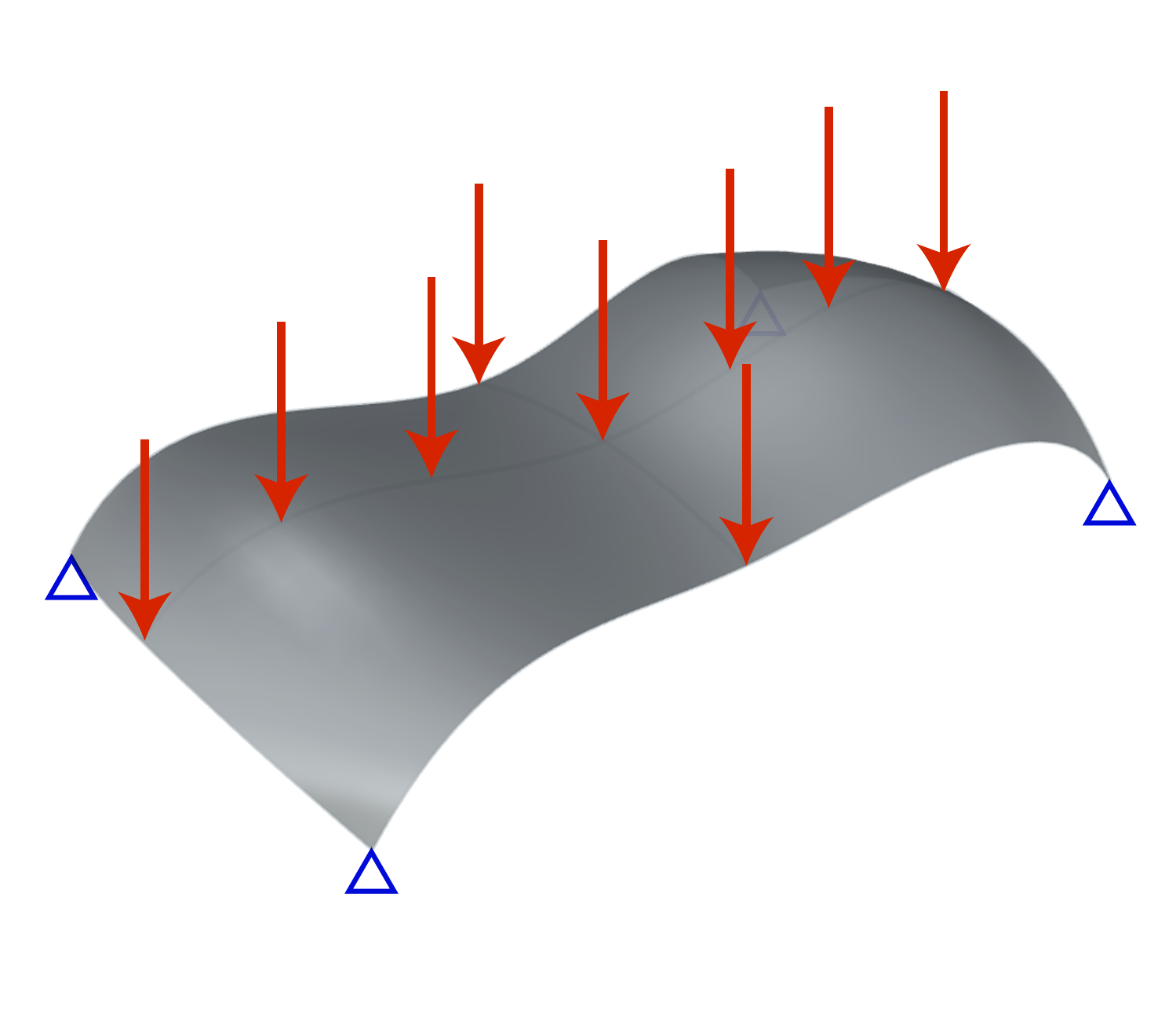}
					\put(28,72){\tiny{$G = 10$}}
					\put(0,70){$S_6$}
				\end{overpic}
				\\
			\end{minipage}
			%
			\begin{minipage}[t]{0.15\linewidth}
				\centering
				\begin{overpic}[scale = 0.12]{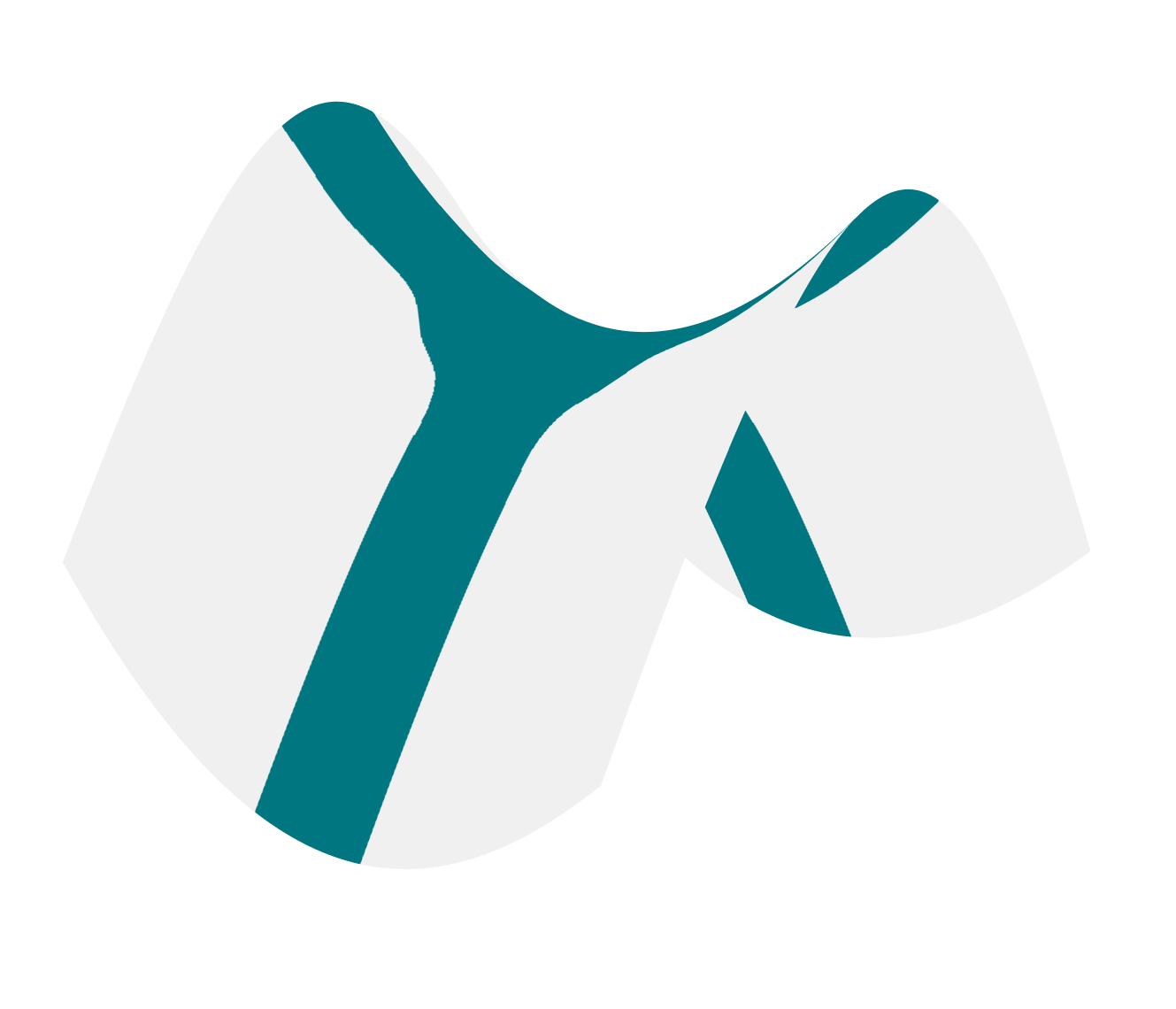}
					\put(15,78){\scriptsize{$V/V_s = 0.300$}}
					\put(25,88){\scriptsize{$C = 2.05$}}
				\end{overpic}\\
				\vspace{0.15cm}
				\begin{overpic}[scale = 0.12]{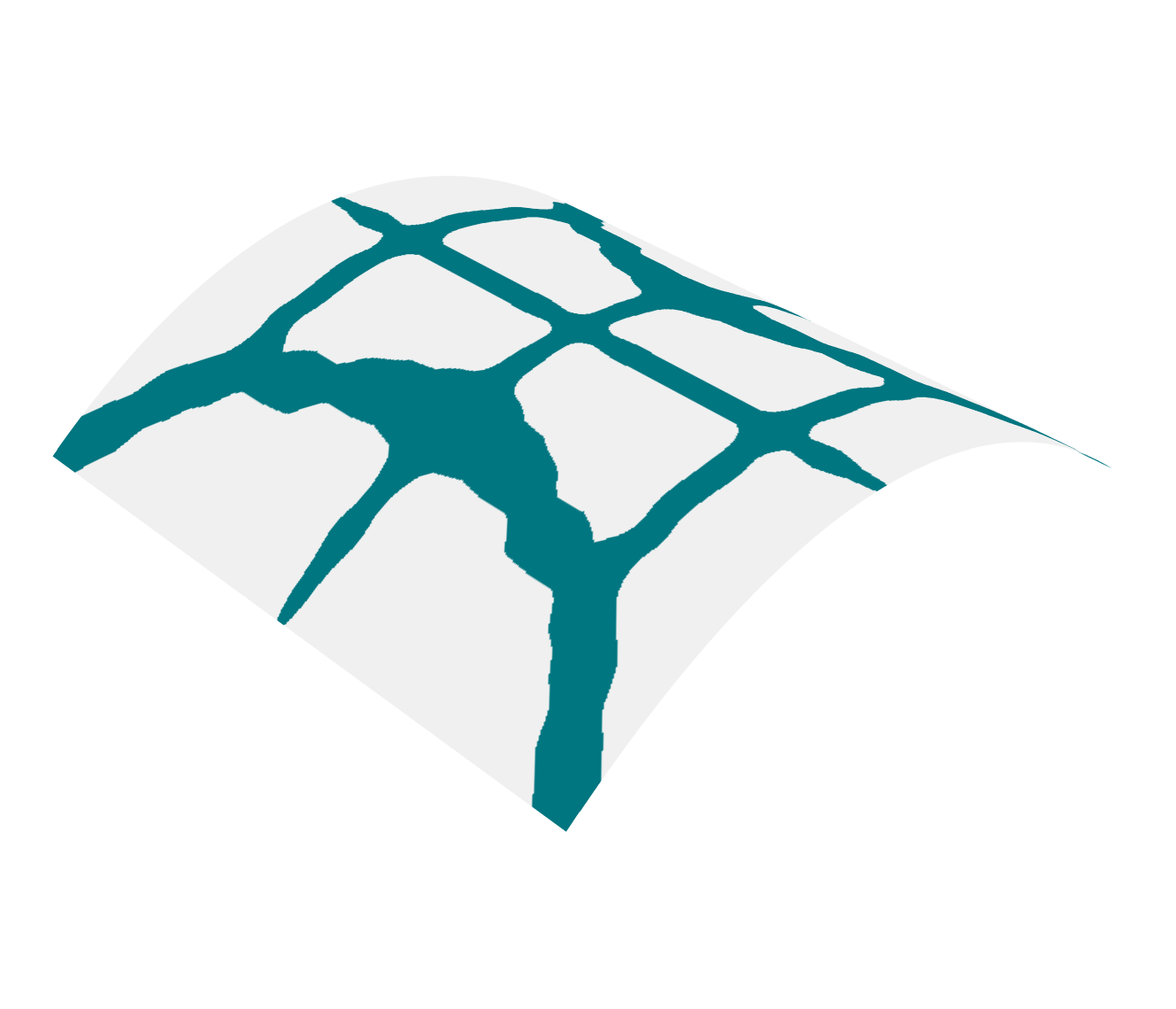}
					\put(14,75){\scriptsize{$V/V_s = 0.305$}}
					\put(20,85){\scriptsize{$C = 172.88$}}
				\end{overpic}\\
				\vspace{0.35cm}
				\begin{overpic}[scale = 0.12]{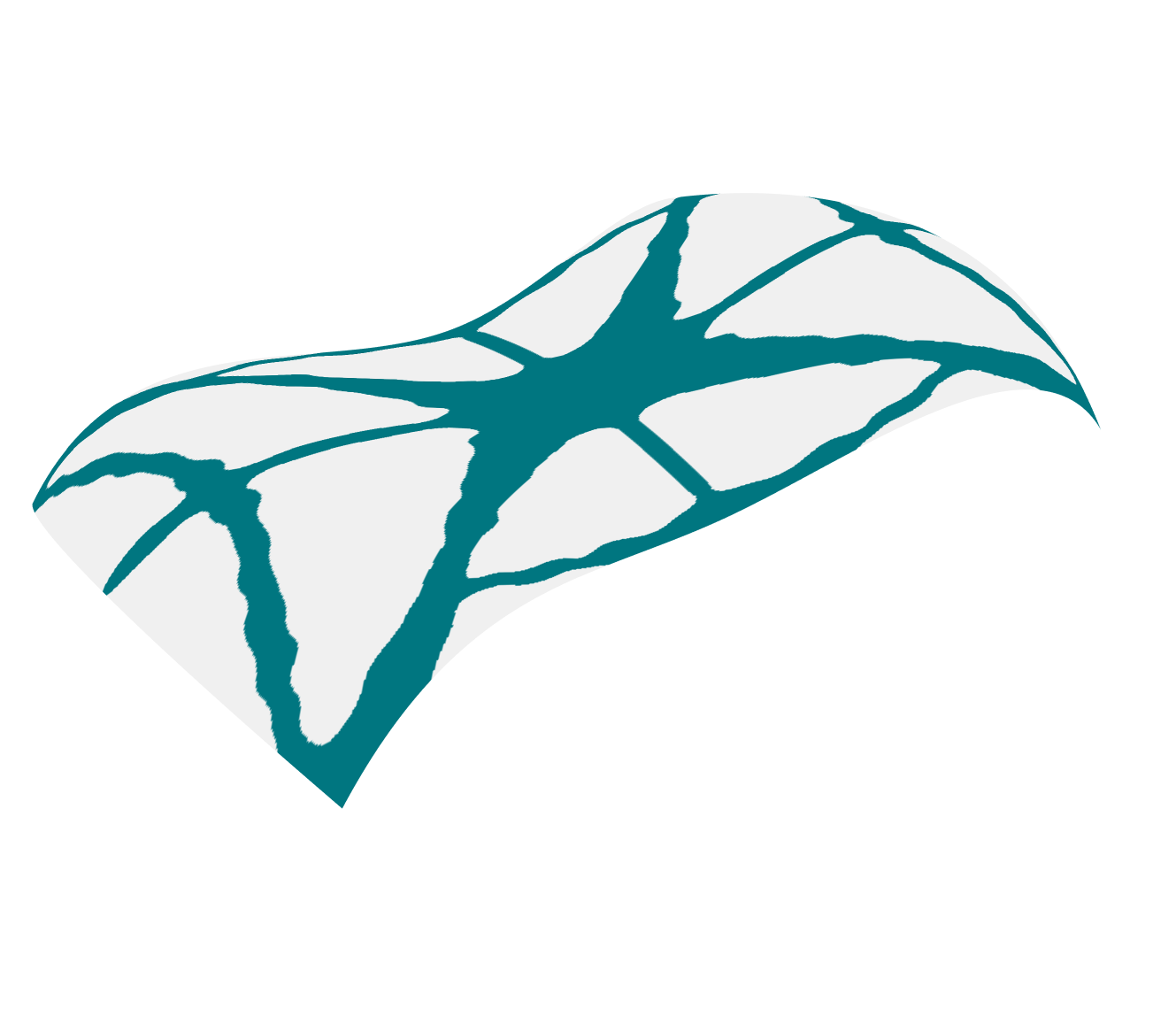}
					\put(14,80){\scriptsize{$V/V_s = 0.352$}}
					\put(21,90){\scriptsize{$C = 497.10$}}
				\end{overpic}\\
			\end{minipage}
			\begin{minipage}[t]{0.15\linewidth}
				\centering
				\begin{overpic}[scale = 0.12]{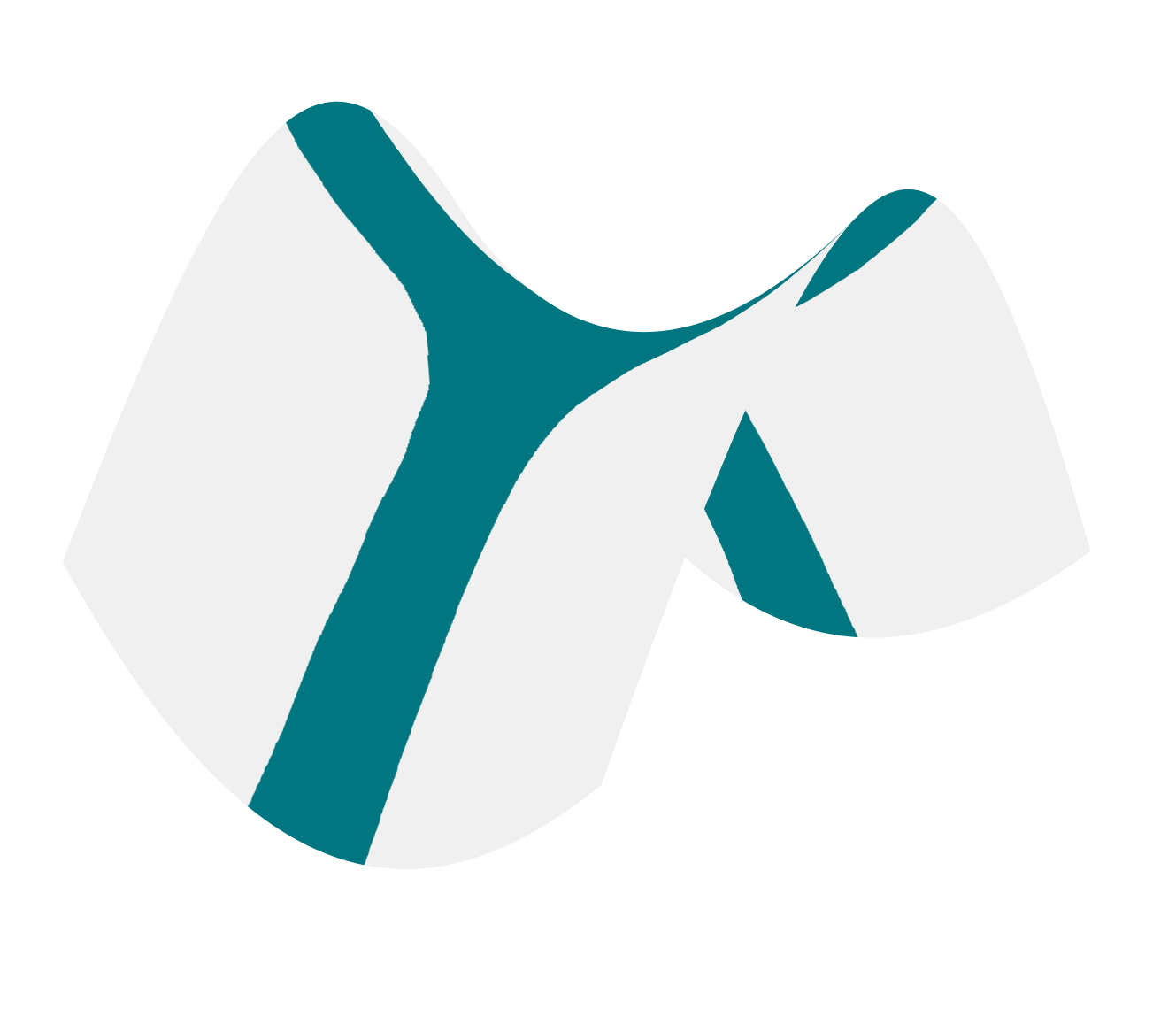}
					\put(15,78){\scriptsize{$V/V_s = 0.301$}}
					\put(26,88){\scriptsize{$C = 2.04$}}
				\end{overpic}\\
				\vspace{0.15cm}
				\begin{overpic}[scale = 0.12]{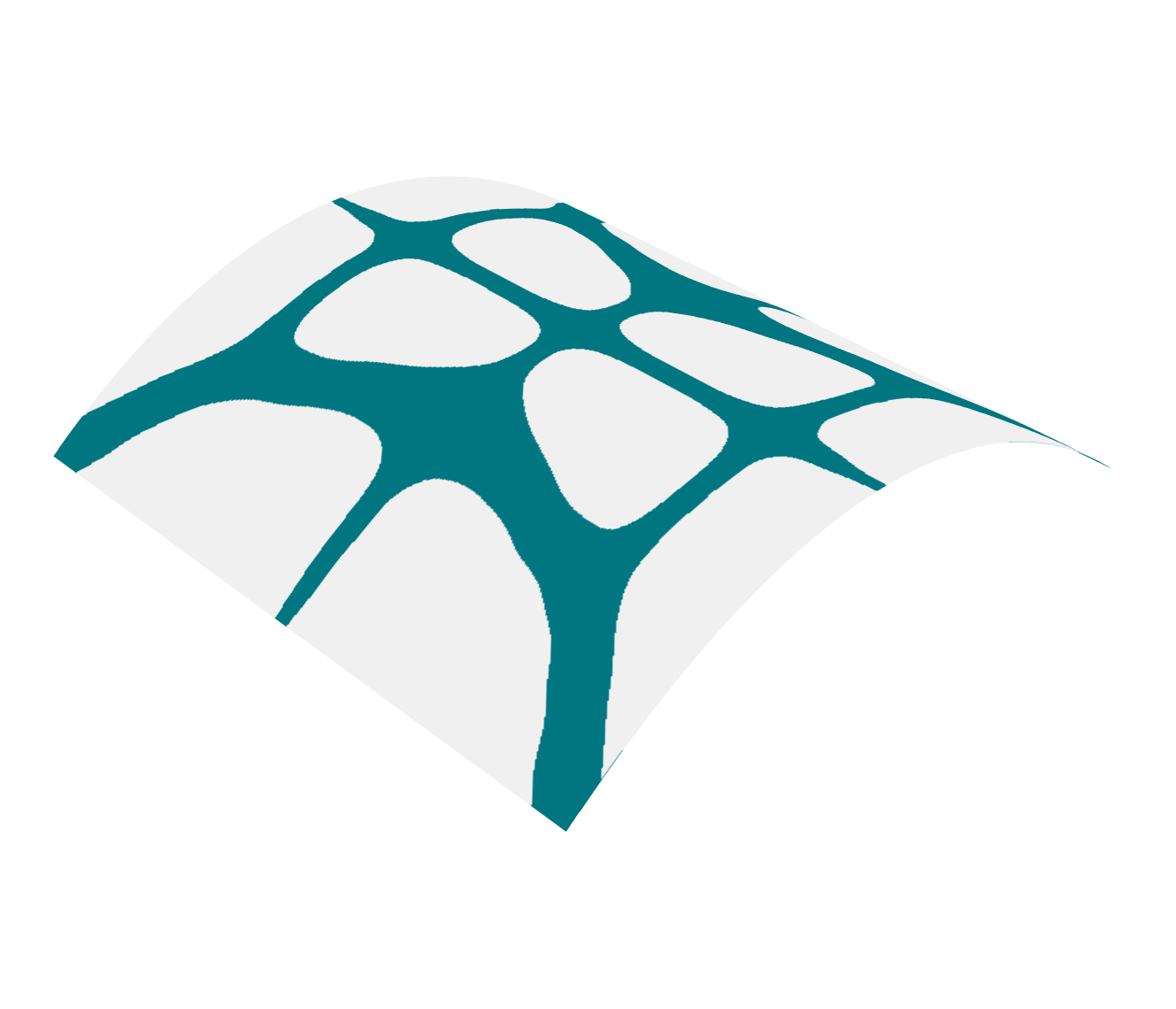}
					\put(14,75){\scriptsize{$V/V_s = 0.316$}}
					\put(20,85){\scriptsize{$C = 163.68$}}
				\end{overpic}\\
				\vspace{0.35cm}
				\begin{overpic}[scale = 0.12]{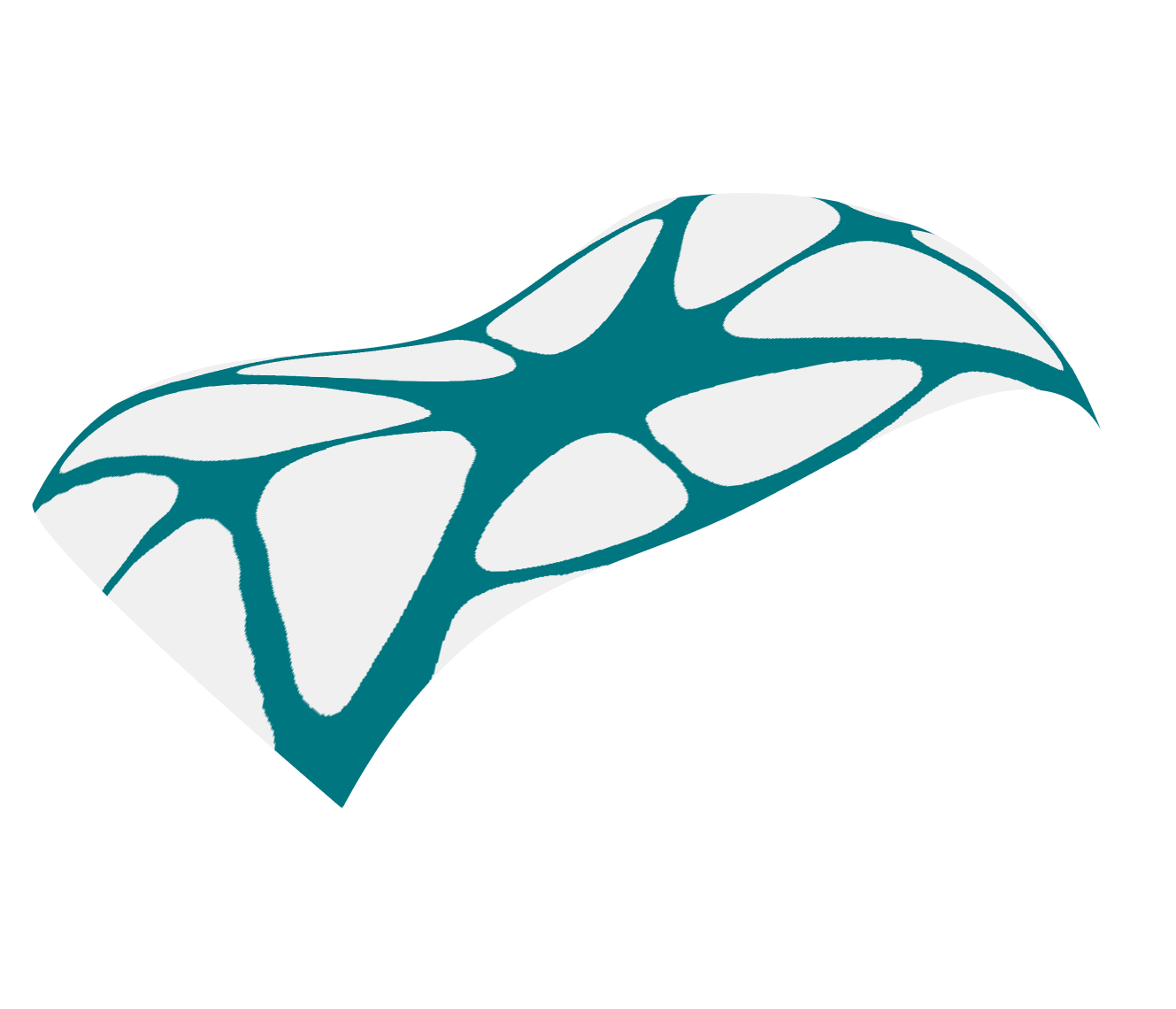}
					\put(14,80){\scriptsize{$V/V_s = 0.366$}}
					\put(21,90){\scriptsize{$C = 454.22$}}
				\end{overpic}\\
			\end{minipage}
			\caption{Fairing results of different shell structures. The upper bound of volume fraction is 0.3 for $S_1,S_2,S_4,S_5$, 0.5 for $S_3$, and 0.35 for $S_6.$}
			\label{fig: smoothing results}
		\end{center}
	\end{figure*}

	Six shell structures are designed using the IGA-SIMP method, resulting in optimal results with wavy boundaries.
	To enhance the quality of the designs for  aesthetics and manufacture purpose, we apply an automatic fairing process, which leads to clearer, explicit and fair boundaries results, as depicted in Fig.~\ref{fig: smoothing results}.
	The compliance and volume of the shell structures before and after the fairing process are also presented in the figure.

	The compliance is comparable in these six cases while the volume of the post-processed structures is slightly larger than the original structures.
	To achieve a structure that closely resembles the original design, the fairing weight $\lambda$ is decreased, thereby minimizing the impact of the fairing process and vice versa.
		This fairing approach can be easily employed with other density-based design representations.

	\begin{table*}[htbp]
		\centering
		\setlength{\abovecaptionskip}{0.35cm} 
		\setlength{\belowcaptionskip}{0.35cm} 
		\caption{Data statistics for 3 cases with different numbers of design variables $N_{DV}$ where $N_E = 50\times 50, V^*=0.3V_s$}
		\begin{tabular*}{\hsize}{@{}@{\extracolsep{\fill}}|c|cccccccccc|@{}}
			\hline
			$N_{DV}$ & $ 5\times 5$ & $ 10\times 10$ & $ 15\times 15$ & $ 20\times 20$ & $ 25\times 25$ & $ 30\times 30$ & $ 35\times 35$ & $ 40\times 40$ & $ 45\times 45$ & $ 50\times 50$\\
			\hline
			case 1 & 15.14 & 9.90 & 9.55 & 9.47 & 9.47 
			& 9.46  & 9.48 & 9.43 & 9.48 & 9.52 \\
			case 2 & 4.31 & 3.31 & 3.27 & 3.36 & 3.36                                      & 3.28 & 3.27 & 3.27 & 3.27 & 3.27 \\
			case 3 & 10838.4 & 740.2 & 266.3 & 256.7 & 241.5 
			& 237.9   & 238.1 & 255.6 & 259.0 & 248.3 \\
			\hline
		\end{tabular*}
		\label{table: 3 cases - numerical results}
	\end{table*}
	\begin{figure*}[htbp]
		\begin{center}
			\begin{overpic}[width=1\textwidth]{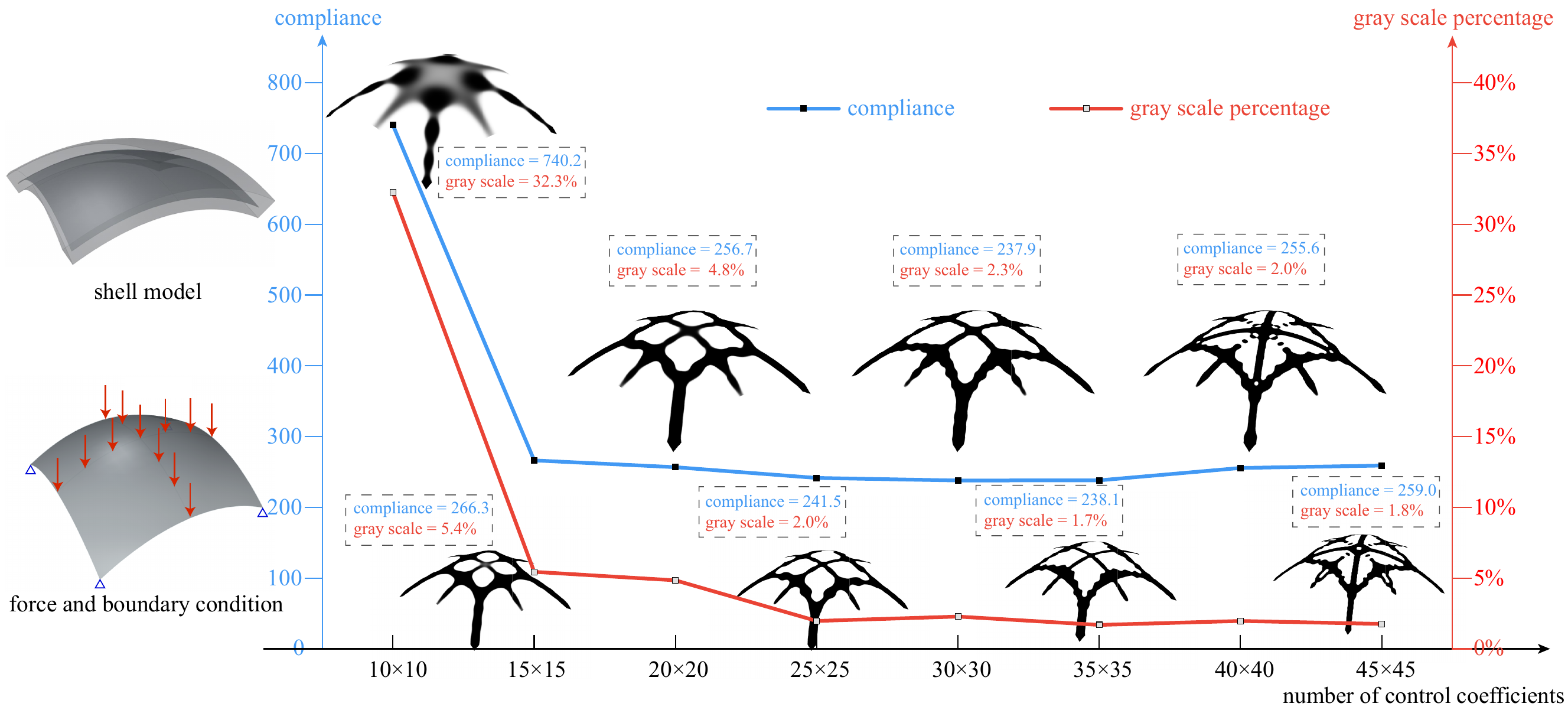}
				\put(8,20){\tiny{$G = 10$}}
			\end{overpic}
		\end{center}
		\caption{Results for case 3, the number of design variables in each direction varied from 10 to 45 with $N_E = 50\times 50, V^*=0.3V_s$. The top left part depicts the force and boundary condition and the bottom left shows the initial structure with a density equal to $0.3$ everywhere. The compliance and grayscale percentage are $(740.2, 32.3\%)$, $(266.3, 5.4\%)$, $(256.7, 4.8\%)$, $(241.5, 2.0\%)$, $(237.9, 2.3\%)$, $(255.6, 2.0\%)$, $(259, 1.8\%)$ respectively, and the grayscale is defined as the elements whose density is between 0.1 and 0.9.}
		\label{fig: discussion of number of control variables}
	\end{figure*}

	\subsection{Parameter selection}
	The choice of the number of design variables $N_{DV}$ significantly impacts the optimization results. In this subsection, we provide a recommended range by testing different $N_{DV}$ through the refinement operation in three cases. The obtained results are presented in Table \ref{table: 3 cases - numerical results}.  It is evident that as $N_{DV}$ increases, the compliance becomes more stable for all three cases.

	The detailed geometric results of case 3 illustrating the significant impact of the number of design variables on the results are shown in Fig.~\ref{fig: discussion of number of control variables}. It should be noted that the fairing process is not applied here and the grayscale is included in the volume. Insufficient design variables lead to high compliance and a large percentage of grayscale. As we can see, the optimization results with $N_{DV} = 10\times 10$ produce an unstable structure with compliance $C = 740.2$ and $32.3\%$ grayscale area.
	Conversely, excessive control coefficients may also lead to unstable structures.  In case 3, for example, using $40\times40$ and $45\times45$ control coefficients results in isolating phenomenon, where isolated blocks appear and additional processing skills are required to address it 
	\citep{gao2019isogeometric}. 

		The occurrence of isolated blocks shares a resemblance with the checkerboard patterns often encountered in standard FEA, primarily attributed to numerical instability with artificially high stiffness~\citep{diaz1995checkerboard}. 
		This can be improved by elevating the degree of basis functions in the design model or applying filter operation for the density coefficients ~\citep{bendsoe2003topology}.

	To achieve more reasonable results and save computational storage, we recommend restricting the number of design variables per parameter direction to between 15 and 35. For more complex structures that require more control points to represent in initial CAD models, we suggest increasing the recommended limit interval accordingly.

	\subsection{Comparison with other IGA-based methods}
	The IGA-SIMP method proposed in this study is distinct from other existing IGA-based methods in its approach to representing optimal structure. 
	While IGA-TSA methods and IGA-MMC/MMV methods use a boundary description model to display the structure topology, the IGA-SIMP method represents the optimal structure with a material description model 
	\citep{gao2020comprehensive}. 
	These two description models are different branches of design methods 
	each with their own advantages and disadvantages
	\citep{kang2013integrated, sigmund2013topology}. 

	The IGA-TSA method utilizes a spline surface and trimming curves to represent the outer and inner boundaries of a structure, improving the accuracy of the analysis and the optimization results. However, the merging and splitting of holes can be intricate, potentially increasing the computational burden. On the other hand, the IGA-MMC/MMV method use moving morphable components or voids to represent material and voids, respectively, allowing for flexible geometry control and overcoming the jagged boundary defect of SIMP methods. However, this approach may come at a significant cost, especially for structures with complex topologies. 
	Furthermore, the MMC/MMV method relies on predefined morphable components/voids to evolve the design, which could limit the final designs to a predetermined set of shapes.
	In comparison, our IGA-SIMP method offers conceptual clarity and is simple to execute in numerical implementation, making it well-suited for complex and porous shell structures.
	
	
        \section{Conclusion and Future Work}
	\label{sec:conclu}
	Shell structures have attracted significant attention in both academic research and industrial production. To achieve a clear and editable boundary representation, we propose an isogeometric topology optimization method based on the SIMP method and automatic post-processing operations. NURBS (Non-Uniform Rational B-Splines) are utilized to represent the shell surface, shape functions for finite element analysis, and material distribution based on the SIMP method.
	The density function establishes a connection between the design space and the geometry space, resulting in a structure that is free from the checkerboard problem, with minimal grayscale and blurred regions. To create porous structures, a local volume constraint is introduced. However, after the topology optimization process, indistinct and wavy boundaries may appear.
	To address this issue, we employ an automatic boundary fairing operation that utilizes NURBS as well. This operation aims to obtain explicit and fairing boundaries. The final results are suitable for easy editing and subsequent procedures in computer-aided design (CAD), computer-aided engineering (CAE), and computer-aided manufacturing (CAM). Several numerical experiments are conducted to demonstrate the features and effectiveness of the proposed method, and the results are thoroughly discussed. Additionally, we provide a reasonable range of optimization parameters. A detailed comparison and analysis of different shell topology optimization algorithms (FEA-SIMP method) are also introduced; It includes a geometric and numerical comparison with the FEA-SIMP method and qualitative comparisons with other IGA-based methods, which can more directly illustrate the advantages of IGA. In this paper, we mainly consider shells with a mid-surface represented by a single patch of NURBS. More complex shell structures for multi-patches optimization problems will be discussed in future work. The optimization process will place greater emphasis on continuity between patch to patch.

	\section*{Declarations}
	\paragraph{Conflict of interest }
	The authors declare that they have no conflict of interest.
	\paragraph{Replication of results}
	Important details for the replication of results have been described in the manuscript. 
	
	\section*{Acknowledgment}
	The authors would like to acknowledge the financial support from  the National Natural Science Foundation of China (61972368, 12371383), the Provincial Natural Science Foundation of Anhui (2208085QA01), the Fundamental Research Funds for the Central Universities (WK0010000075), the Open Project Program of the State Key Laboratory of CAD \& CG (Grant No. A2303), Zhejiang University.

\bibliography{mybibfile}

\appendix
\section{Appendix}
	\label{sec:supplementary}
	
	\subsection{Shell structures analysis}
	\label{subsec:analysis}
	As presented in Eqs.\eqref{eq: shell}-\eqref{eq: displacement}, the geometry and displacement of a shell structure are defined by parameters $(s,t,\zeta)$ and are represented by a conjunction of two parts, one for the mid-surface and the other for the thickness direction. The displacement is determined by 5 displacement variables $(u_{ij},v_{ij},w_{ij},\alpha_{ij},\beta_{ij})$. In this part, our study will focus on deducing the static equilibrium equation and shell properties including strain, stress, stiffness matrix and force.

	The strain is calculated by taking the partial derivatives of displacement $\boldsymbol{u}$ with respective to $x,y,z$. To achieve this, we employ the chain rule and link it to $5$ displacement variables for each control point. The resulting representation takes the form of 
	\begin{equation}
		\boldsymbol{\epsilon}_{global} =  \sum\limits_{ij}  \boldsymbol{H} \boldsymbol{\Gamma \Phi}_{ij}  \left( u_{ij}, v_{ij}, w_{ij}, \alpha_{ij}, \beta_{ij} \right)^\T,
	\end{equation}
	where $\boldsymbol{H}$ is a constant $6 \times 9$  matrix that establishes the connection between the strain and the vector
	$$\{\frac{\partial{u}}{\partial{x}},\frac{\partial{u}}{\partial{y}},\frac{\partial{u}}{\partial{z}},\frac{\partial{v}}{\partial{x}},\frac{\partial{v}}{\partial{y}},\frac{\partial{y}}{\partial{z}},\frac{\partial{w}}{\partial{x}},\frac{\partial{w}}{\partial{y}},\frac{\partial{w}}{\partial{z}}\}.$$ $\boldsymbol{\Gamma}$, a $9 \times 9$ matrix, further transforms the partial derivatives to those with respect to parameters $(s,t,\zeta)$ and $\boldsymbol{\Phi}_{ij}$ correlate these partial derivatives with displacement variables. To comply with the assumption of zero normal stress in the local system ${\boldsymbol{v}_1,\boldsymbol{v}_2,\boldsymbol{v}_3}$ in the Reissner-Mindlin theory, we compute strain and stress in the local system. This results in a local strain
	\begin{equation}
		\boldsymbol{\epsilon}_{local} = \sum\limits_{ij} \boldsymbol{T H \Gamma} \boldsymbol{\Phi}_{ij} \left( u_{ij}, v_{ij}, w_{ij}, \alpha_{ij}, \beta_{ij} \right)^\T,
	\end{equation}
	which is determined through the use of $\boldsymbol{T}$, a $6\times 6$ orthogonal matrix.  $\boldsymbol{T}$ is obtained through the chain rule of multivariate function derivation and depicts the transformation between local and global strain.

	The local stress is calculated by
	\begin{equation}
		\boldsymbol{\sigma}_{\text{local}} = \boldsymbol{D} \boldsymbol{\epsilon}_{local},
	\end{equation}
	where $\boldsymbol{D}$ is a stain-stress matrix related to Young's modulus $E$ of the isotropic material according to \cite{bendsoe2003topology}, which is defined as  
	\begin{equation}
		\label{eq: SIMP Young's modulus}
		E(\rho) = E_{min} + \rho^{\iota} (E_0-E_{min}).
	\end{equation}
	where $E_{min}$ is a minimum limit. $E_0$ denotes the exact
	Young’ s modulus for the isotropic material and $\iota = 5$. 
	Then the global stress is represented as
	\begin{equation}
		\boldsymbol{\sigma}_{global} 
		= \boldsymbol{T}^{-1} \boldsymbol{\sigma}_{local} 
		= \boldsymbol{T}^{-1} \boldsymbol{D} \epsilon_{local} 
		= \boldsymbol{T}^{\T} \boldsymbol{D T} \boldsymbol{\epsilon}_{global}.
	\end{equation}

	The stiffness matrix of an solid element is built as 
	\begin{equation}
		\begin{aligned}
			\boldsymbol{K}_e^0  &= \int\limits_{s}  \int\limits_{t} \int\limits_{\zeta} \boldsymbol{B}^\T \boldsymbol{D B} \left|det(\boldsymbol{J}) \right| \  ds dt d\zeta, \\
			&\approx \sum\limits_{k=1}^{N_G} \boldsymbol{B}^\T_k \boldsymbol{D} \boldsymbol{B}_k \left|det(\boldsymbol{J}_k)\right| W_k,
		\end{aligned}
	\end{equation}
	where $N_G$ denotes the number of Gauss integration points $\left\{(s_k,t_k,\zeta_k)\right\}_{k=1}^{N_G}$ and $\left\{W_k\right\}_{k=1}^{N_G}$ are the corresponding Gauss weights. 
	At each Gauss point, the Jacobi matrix $J_k$ is calculated, which depicts the transformation from the physical coordinate system ${x,y,z}$ to the parametric space ${s,t,\zeta}$. 
	The matrix $\boldsymbol{B}$ is obtained as $\boldsymbol{B}=\boldsymbol{TH\Gamma \boldsymbol{\Phi}}$,  where $\boldsymbol{\Phi}$ is a $9 \times 5(p+1)(q+1)$ matrix composed of a group of $\left\{\boldsymbol{\Phi}_{ij}\right\}$ corresponding to those NURBS basis functions $N_{ij}(s,t) \neq 0$, in other words, those NURBS basis functions whose support interact with element $e$. 
	Note that the solid shell elements are defined based on the curved elements of the mid-surface divided by the isoparametric curves corresponding to knot spans of $\mathcal{S}$ and $\mathcal{T}$. 
	The calculation of the stiffness matrix $\boldsymbol{K}_{e}$ for element $e$ is similar to that of the modulus in Eq.\eqref{eq: SIMP Young's modulus}, with $\boldsymbol{K}_e^0$ and $\boldsymbol{K}_e$ corresponding to $E_0$ and $E$, respectively.
	
	Additionally, the expression of load vector $\boldsymbol{F}$ is 
	\begin{equation}
		\label{eq: force}
		\boldsymbol{F} =\int\limits_{s}  \int\limits_{t} \int\limits_{\zeta} \boldsymbol{M}^\T \boldsymbol G |det(\boldsymbol J)| ds dt d\zeta
	\end{equation}
	where G is the force onto the shell structure and $\boldsymbol{J}$ is the Jacobi matrix. $\boldsymbol{M}$ is the basis function matrix which connect the displacement in Eq.\eqref{eq: displacement} with the displacement coefficients vector $\boldsymbol{U}$, that is, $\boldsymbol{u}(s,t) = \boldsymbol{M} \boldsymbol{U}$.

	\subsection{Sensitivity analysis}
	\label{subsec:sensitivity analysis}
	With the representation of compliance in Eq.~\eqref{eq: compliance} and assuming that the load is independent of density,  we can compute the first-order partial derivative with respect to density coefficients as follows,
	\begin{equation}
		\label{eq: partial dirivative of C with respect to rho}
		\frac{\partial C}{\partial \rho_{ij}} = \sum\limits_{e}  \frac{\partial {(\boldsymbol{U}_e^{\T} \boldsymbol{F}_e)}}{\partial \rho_{ij}} .
	\end{equation}
	Using the adjoint method, we set $\hat{C} = C + \boldsymbol{\lambda}^{\T}(\boldsymbol{KU}-\boldsymbol{F})$ and obtain that 
	\begin{equation}
		\frac{\partial{\hat{C}}}{\partial{\rho_{ij}}} 
		= (\boldsymbol{K}^{\T} \boldsymbol{\lambda} + \boldsymbol{F})^{\T} \frac{\partial \boldsymbol{U}}{\partial \rho_{ij}} + \boldsymbol{\lambda} \frac{\partial{\boldsymbol{K}}}{\partial{\rho_{ij}}} \boldsymbol{U}.
	\end{equation}
	We take $\boldsymbol{\lambda} = -\boldsymbol{U}$ to cancel out terms with $\frac{\partial{\boldsymbol{U}}}{\partial{\rho_{ij}}} $.
	Then Eq.~\eqref{eq: partial dirivative of C with respect to rho} becomes 
	\begin{equation}
		\begin{aligned}
			\frac{\partial{C}}{\partial{\rho_{ij}}}
			&= \frac{\partial{\hat{C}}}{\partial{\rho_{ij}}} \\
			&= -\sum\limits_{e} \boldsymbol{U}_e^\T \frac{\partial \boldsymbol{K}_e}{\partial \rho_{ij}} \boldsymbol{U}_e \\
			&= - \sum\limits_{e}  \boldsymbol{U}_e^\T \cdot \tilde\rho_e^{\tilde{p}-1} (\boldsymbol{K_e}^0 - \boldsymbol{K}_{min}) \cdot \boldsymbol{U}_e \cdot  \tilde{p} \frac{\partial \tilde\rho_e}{\partial \rho_{ij}}.
		\end{aligned}
	\end{equation}
	where $\tilde\rho_e = H(\rho_e)$ in Eq.\eqref{eq: heaviside function} and
	$$\frac{\partial \tilde\rho_e}{\partial \rho_{ij}} = H^{\prime}({\rho}_e) \cdot \frac{\partial {\rho}_e}{\partial \rho_{ij}} = H^{\prime}(\tilde{\rho}_e) \cdot N_{ij}$$
	
	The first-order partial derivative of volume in Eq.~\eqref{eq: volume}  with respect to density coefficient $\rho_{ij}$ is 
	\begin{equation}
		\label{eq: partial dirivative of V with respect to rho}
		\frac{\partial V}{\partial \rho_{ij}} = \sum\limits_{e} V_e^0 \cdot  \frac{\partial \tilde\rho_e}{\partial \rho_{ij}} = \sum\limits_{e} V_e^0 \cdot H^{\prime}(\tilde{\rho}_e) \cdot N_{ij}, 
	\end{equation}
	
	Note that the sensitivity analysis discussed above pertains only to the first formulation in Eq.~\eqref{eq: formulation 1}. And a similar discussion can be conducted for the second one in Eq.~\eqref{eq: formulation 2} as
	\begin{equation}\nonumber
		\begin{aligned}
			\frac{\partial{C}}{\partial{\rho_{ij}}}
			&= -\sum\limits_{e} \boldsymbol{U}_e^\T \frac{\partial \boldsymbol{K}_e}{\partial \rho_{ij}} \boldsymbol{u}_e \\
			&= - \sum\limits_{e}  \boldsymbol{U}_e^\T \cdot \tilde\rho_e^{\tilde{p}-1} (\boldsymbol{K_e}^0 - \boldsymbol{K}_{min}) \cdot \boldsymbol{U}_e \cdot  \tilde{p} \cdot H^{\prime}({\rho}_e) \cdot N_{ij},
		\end{aligned}
	\end{equation}
	and
	\begin{equation}\nonumber
		\frac{\partial \Bar{V}}{\partial \rho_{ij}} =  \frac{1}{\gamma}\left(\frac{1}{N_E}\sum_{e=1}^{N_E} (\Bar{\rho}_e)^{\gamma}\right)^{\frac{1}{\gamma}-1}
		\cdot
		\frac{1}{N_E} \sum_{e=1}^{N_E} \gamma \left(\Bar{\rho}_e\right)^{\gamma-1}
		\cdot
		\frac{\partial \Bar\rho_e}{\partial \rho_{ij}}, 
	\end{equation}
	where $\Bar{\rho}_e = \frac{\sum\limits_{f \in \mathcal{N}_e} \tilde\rho_f V_f }{\sum\limits_{f \in \mathcal{N}_e} V_f} $ and
	$$\frac{\partial \Bar\rho_e}{\partial \rho_{ij}} = \frac{\sum\limits_{f \in \mathcal{N}_e} V_f H^{\prime}(\rho_f) N_{ij}}{\sum\limits_{f \in \mathcal{N}_e} V_f}.$$
	

\end{document}